\documentstyle[amssymb,12pt]{amsart}
\newtheorem{prop}{Proposition}[section]
\newtheorem{prop:def}{Proposition-Definition}[section]

\newtheorem{lemma}{Lemma}[section]
\newtheorem{thm}{Theorem}[section]
\newtheorem{cor}{Corollary}[section]

\theoremstyle{remark}
\newtheorem{remark}{Remark}

\begin{document}
\newcommand{\nc}{\newcommand} \nc{\on}{\operatorname}
\nc{\pa}{\partial}
\nc{\cA}{{\cal A}}\nc{\cB}{{\cal B}}\nc{\cC}{{\cal C}}
\nc{\cE}{{\cal E}}\nc{\cG}{{\cal G}}\nc{\cH}{{\cal H}}
\nc{\cX}{{\cal X}}\nc{\cR}{{\cal R}}\nc{\cL}{{\cal L}}
\nc{\cK}{{\cal K}}\nc{\cV}{{\cal V}}\nc{\cI}{{\cal I}} 
\nc{\cU}{{\cal U}}
\nc{\sh}{\on{sh}}\nc{\Id}{\on{Id}}\nc{\Diff}{\on{Diff}}
\nc{\ad}{\on{ad}}\nc{\Der}{\on{Der}}\nc{\End}{\on{End}} 
\nc{\Hom}{\on{Hom}} 
\nc{\Alt}{\on{Alt}}\nc{\res}{\on{res}}\nc{\ddiv}{\on{div}}
\nc{\card}{\on{card}}\nc{\dimm}{\on{dim}}
\nc{\Jac}{\on{Jac}}\nc{\Ker}{\on{Ker}}\nc{\Sh}{\on{Sh}}
\nc{\Imm}{\on{Im}}\nc{\limm}{\on{lim}}\nc{\Ad}{\on{Ad}}
\nc{\ev}{\on{ev}}\nc{\Sym}{\on{Sym}} \nc{\ol}{\overline}
\nc{\Hol}{\on{Hol}}\nc{\Det}{\on{Det}}\nc{\Ham}{\on{Ham}}
\nc{\de}{\delta}\nc{\si}{\sigma}\nc{\ve}{\varepsilon}
\nc{\al}{\alpha}\nc{\vp}{\varphi}
\nc{\CC}{{\mathbb C}}\nc{\ZZ}{{\mathbb Z}}
\nc{\NN}{{\mathbb N}}\nc{\zz}{{\mathbf z}}\nc{\kk}{{{\mathbf k}}}
\nc{\bl}{{{\mathbf l}}}\nc{\pp}{{{\mathbf p}}} 
\nc{\ff}{{{\mathbf f}}}
\nc{\AAA}{{\mathbb A}}\nc{\BB}{{\mathbb B}}
\nc{\cO}{{\cal O}} \nc{\cF}{{\cal F}}
\nc{\cS}{{\cal S}}\nc{\cW}{{\cal W}}\nc{\cJ}{{\cal J}} 
\nc{\cZ}{{\cal Z}}
\nc{\la}{{\lambda}}\nc{\G}{{\mathfrak g}}\nc{\mm}{{\mathfrak m}}
\nc{\J}{{\mathfrak j}}
\nc{\A}{{\mathfrak a}}\nc{\gotS}{{\mathfrak S}}
\nc{\HH}{{\mathfrak h}} \nc{\T}{{\mathfrak t}}
\nc{\N}{{\mathfrak n}}\nc{\B}{{\mathfrak b}}\nc{\LL}{{\mathfrak l}}
\nc{\La}{\Lambda}
\nc{\g}{\gamma}\nc{\eps}{\epsilon}\nc{\wt}{\widetilde}
\nc{\wh}{\widehat}
\nc{\bn}{\begin{equation}}\nc{\en}{\end{equation}}
\nc{\SL}{{\mathfrak{sl}}}\nc{\ttt}{{\mathfrak{t}}}
\nc{\s}{{\mathfrak{s}}}

%
%
%

\newcommand{\ldar}[1]{\begin{picture}(10,50)(-5,-25)
\put(0,25){\vector(0,-1){50}}
\put(5,0){\mbox{$#1$}} 
\end{picture}}

\newcommand{\lrar}[1]{\begin{picture}(50,10)(-25,-5)
\put(-25,0){\vector(1,0){50}}
\put(0,5){\makebox(0,0)[b]{\mbox{$#1$}}}
\end{picture}}

\newcommand{\luar}[1]{\begin{picture}(10,50)(-5,-25)
\put(0,-25){\vector(0,1){50}}
\put(5,0){\mbox{$#1$}}
\end{picture}}

\title[PBW and duality for quantum groups and quantum current algebras]
{PBW and duality theorems for quantum groups and quantum current algebras}

\author{B. Enriquez}

\address{Centre de Math\'ematiques, Ecole Polytechnique, 
UMR 7640 du CNRS, 91128 Palaiseau, France}



\date{April 1999}

\begin{abstract}
We give proofs of the PBW and duality theorems for the quantum 
Kac-Moody algebras and quantum current algebras, relying on Lie bialgebra 
duality. We also show that the classical limit of the quantum current 
algebras associated with an untwisted affine Cartan matrix is the enveloping 
algebra of a quotient of the corresponding toroidal algebra; this quotient is 
trivial in all cases except the $A_1^{(1)}$ case. 
\end{abstract}

\maketitle

\section{Outline of results}

\subsection{Quantum Kac-Moody algebras} \label{KM}

Let $A = (a_{ij})_{1\leq i,j\leq n}$ be a symmetrizable Cartan matrix. 
Let $(d_i)_{1\leq i \leq n}$ be the coprime positive integers such that
the matrix  $(d_i a_{ij})_{1\leq i,j\leq n}$ is symmetric. Let $r$ be
the  rank of $A$; we assume that the matrix $(a_{ij})_{n-r+1\leq i,j\leq
n}$ is nondegenerate. 

Let $\G$ be the Kac-Moody Lie algebra associated with $A$; let $\N_+$ be
its positive pro-nilpotent subalgebra and $(\bar e_i)_{1\leq i\leq n}$
be the generators of $\N_+$ corresponding to the simple roots of $\G$.  

Let $\CC[[\hbar]]$ be the formal series ring in $\hbar$.  Let
$U_\hbar\N_+$ be the quotient of the free algebra with $n$ generators
$\CC[[\hbar]]\langle e_i, i = 1,\ldots,n\rangle$ by the two-sided ideal
generated by the quantum Serre relations  \begin{equation}
\label{quantum:serre}  \sum_{k=0}^{1 - a_{ij}} 
(-1)^k \bmatrix 1 - a_{ij} \\ k
\endbmatrix_{q^{d_i}}  e_i^k e_j e_i^{1 - a_{ij} - k} = 0,  
\end{equation}  where $\bmatrix m \\ p \endbmatrix_{q} = {{ [m]^!_q
}\over{ [p]^!_q  [m-p]^!_q }}$, $[k]^!_q = [1]_q\cdots[k]_q$, $[k]_q = 
{{q^k - q^{-k}}\over{q - q^{-1}}}$, and $q = e^\hbar$
(\cite{Drinf:def,Jimbo}). 
   
We will show: 

\begin{thm} \label{thm:first}
$U_\hbar\N_+$ is a free $\CC[[\hbar]]$-module, and the map 
$e_i\mapsto \bar e_i$ defines an algebra isomorphism of $U_\hbar\N_+
/\hbar U_\hbar\N_+$ with $U\N_+$. 
\end{thm}

This Theorem may be derived from the Poincar\'e-Birkhoff-Witt (PBW)
results of Lusztig's book \cite{Lusztig}; in the case $\G = \SL_n$, it can
also be derived from those of Rosso (\cite{Rosso:PBW}), and in the cases when 
$\G$ is semisimple or untwisted affine, from those of \cite{Kh:T}. 

The proof presented here is based on the comparison of $U_\hbar\N_+$ 
with a quantum shuffle algebra, Lie bialgebra duality and the
Deodhar-Gabber-Kac theorem. 

As a corollary of this  proof, we show 

\begin{cor} \label{cor:comparison}
The map $p_\hbar$ defined in Lemma \ref{pourim} is an algebra 
isomorphism from $U_\hbar\N_+$ to the subalgebra $\langle V \rangle$
of the shuffle algebra $\Sh(V)$ defined in sect.\ \ref{sect:shuffles}. 
\end{cor}

This result was proved in \cite{Rosso:shuffle}; it can also be derived from the 
results of \cite{Schau}. 
Rosso's proof uses the nondegeneracy of the pairing between opposite 
Borel quantum algebras (\cite{Lusztig}, Cor.\ 33.1.5; see also Thm.\
\ref{thm:second}).  Schauenburg shows that $\langle V \rangle$ is
isomorphic to the quotient of  the free algebra generated by the $e_i$
by the radical of a braided Hopf pairing.  Together with
\cite{Lusztig}, Cor.\ 33.1.5, this implies Cor.\ \ref{cor:comparison}. 
 
Define $U_\hbar\N_-$ as the algebra with generators $f_i$, $i = 1,\ldots,n$,  
and the same defining relations as $U_\hbar\N_+$ (with $e_i$ replaced by $f_i$). 
Define a grading on $U_\hbar\N_\pm$ by $(\pm\NN)^n$ by $\deg(e_i) = \eps_i$, 
$\deg(f_i) = -\eps_i$, where $\eps_i$ is the $i$th basis vector of $\NN^n$, and
define the braided tensor products $U_\hbar\N_\pm \bar\otimes U_\hbar\N_\pm$
as the algebras isomorphic to $U_\hbar\N_\pm \otimes_{\CC[[\hbar]]} U_\hbar\N_\pm$
as $\CC[[\hbar]]$-modules, with mutiplication rule 
\begin{equation} \label{braided:tensor:pdt} 
(x\otimes y)(x'\otimes y')
= q^{-\langle \deg(x'), \deg(y)\rangle}(xx'\otimes yy');
\end{equation} 
we set $\langle \eps_i,\eps_j\rangle = d_i a_{ij}$. 

Then $U_\hbar\N_\pm$ are endowed with braided Hopf algebra structures,
defined by $\Delta_+(e_i) = e_i\otimes 1 + 1\otimes e_i$,  and 
$\Delta_-(f_i) = f_i\otimes 1 + 1\otimes f_i$.  

In \cite{Drinf:ICM}, Drinfeld showed that there exists a unique pairing 
$\langle , \rangle_{U_\hbar\N_\pm}$ of $U_\hbar\N_+$ and $U_\hbar\N_-$ 
with values in $\CC((\hbar)) = \CC[[\hbar]][\hbar^{-1}]$, defined by 
\begin{equation} \label{pair:1}
\langle e_i, f_{i'} \rangle_{U_\hbar\N_\pm} = {1\over\hbar}
d_i^{-1}\delta_{ii'}, 
\end{equation}
\begin{equation} \label{pair:2}
\langle x,yy' \rangle_{U_\hbar\N_\pm} = \sum \langle x^{(1)},
y\rangle_{U_\hbar\N_\pm}   \langle x^{(2)}, y'\rangle_{U_\hbar\N_\pm},
\end{equation}
and 
\begin{equation} \label{pair:3}
\langle xx',y \rangle_{U_\hbar\N_\pm} =  \sum \langle x,
y^{(1)}\rangle_{U_\hbar\N_\pm}   \langle x',
y^{(2)}\rangle_{U_\hbar\N_\pm}
\end{equation} 
for $x,x'$ in $U_\hbar\N_+$ and $y,y'$ in $U_\hbar\N_-$, and
$\Delta_\pm(z) = \sum z^{(1)} \otimes z^{(2)}$ (braided Hopf pairing
axioms). 

As a direct consequence of Cor.\ \ref{cor:comparison}, we show: 

\begin{thm}  \label{thm:second}
The pairing $\langle , \rangle_{U_\hbar\N_\pm}$ between $U_\hbar\N_+$
and $U_\hbar\N_-$ is nondegenerate.  
\end{thm}

This result can be found in Lusztig's book  (\cite{Lusztig},  Cor.\
33.1.5, Def.\ 3.1.1 and Prop.\ 3.2.4); it relies  on the construction of
dual PBW bases. Another argument using  Lusztig' results on integrable modules
is in \cite{Tanisaki}, and an argument using irreducible Verma modules is in
\cite{Rosso:CRAS}. 

We also show: 
\begin{prop} \label{R:mat}
For any $\al$ in $\NN^n$, let $U_\hbar\N_\pm[\pm\al]$ be the part of 
$U_\hbar\N_\pm$ of degree $\pm\al$, and let $P[\al]$ be the element of
$U_\hbar\N_+[\al] \otimes U_\hbar\N_-[-\al]$ induced by $\langle ,
\rangle_{U_\hbar\N_\pm}$.  Let $\Delta_+$ be the set of positive roots
of $\G$ (the $\eps_i$ are the simple roots).  Let $(\bar
e_{\al,i})_{\al\in\Delta_+}$  and $(\bar f_{\al,i})_{\al\in\Delta_+}$ be dual
Cartan-Weyl bases of $\N_+$ and $\N_-$, and let $e_{\al,i},f_{\al,i}$ be lifts
of the $\bar e_{\al,i},\bar f_{\al,i}$ to $U_\hbar\N_\pm$.  Then, if $k$
is the integer such that $\al$ belongs to $k\Delta_+ \setminus
(k-1)\Delta_+$, we have 
$$
P[\al] = {{\hbar^k}\over{k!}} \sum_{\al_1,\ldots,\al_k\in\Delta_+ | \sum_i\al_i = \al; i_j}
e_{\al_1,i_1} \cdots e_{\al_k,i_k} \otimes f_{\al_1,i_1}\cdots f_{\al_k,i_k} + o(\hbar^k). 
$$
\end{prop}
The fact that $P[\al]$ has $\hbar$-adic valuation equal to $k$ was
stated by   Drinfeld in \cite{Drinf:ICM}. 

\subsubsection{The case of a generic deformation parameter}
\label{generic}

It is easy to derive from the above results, PBW and nondegeneracy  
results in the case where the parameter $q = e^\hbar$ is 
generic. 

\begin{cor} \label{generic:1}
Let $q'$ be an indeterminate, and let $U_{q'}\N_+$ be the algebra over
$\CC(q')$ with generators $e'_i$, $i = 1,\ldots,r$, and relations
(\ref{quantum:serre}),  with $e_i$ and $q = e^\hbar$ replaced by $e'_i$
and $q'$.  We have for any $\al$ in $\NN^n$, $\dim_{\CC(q')}
U_{q'}\N_+[\al]   = \dim_{\CC(q')} U_{q'}\N_+[\al]$. 

Let $(\bar e_\nu)_{\nu\in I}$ be a basis of homogeneous elements of $\N_+$. 
Let 
$$
\bar e_\nu = \sum_{i_j\in \{1,\ldots,n\}} C_{\nu; i_1,\ldots,i_k}
\bar e_{i_1} \cdots \bar e_{i_k}
$$ 
be expressions of the $\bar e_\nu$ in terms
of the generators $\bar e_1,\ldots,\bar e_n$. 

Let $C_{\nu; i_1,\ldots,i_k}(q')$ be rational functions of $q'$, such that 
$C_{\nu; i_1,\ldots,i_k}(1) = C_{\nu; i_1,\ldots,i_k}$, and set 
$e'_\nu = \sum_{i_j\in \{1,\ldots,n\}} C_{\nu; i_1,\ldots,i_k}(q')
e_{i_1} \cdots e_{i_k}$. Then the family 
$(\prod_{\nu} e_\nu^{n_\nu})$, where the $n_\nu$ are in $\NN$ and 
vanish except for a finite number of them, forms a basis of 
$U_{q'}\N_+$ (over $\CC(q')$).  
\end{cor}

\begin{cor} \label{generic:2}
Let $U_{q'}\N_-$ be the $\CC(q')$-algebra with generators $f'_i,1\leq i \leq n$, 
and relations (\ref{quantum:serre}), with $e_i$ replaced  by $f'_i$. Define
the braided tensor squares $U_{q'}\N_\pm\bar\otimes U_{q'}\N_\pm$ using 
(\ref{braided:tensor:pdt}). We have a braided Hopf pairing 
$\langle , \rangle_{\CC(q')}$ between $U_{q'}\N_+$ and $U_{q'}\N_-$, 
defined by (\ref{pair:1}), (\ref{pair:2}) and (\ref{pair:3}).  
The pairing $\langle , \rangle_{\CC(q')}$ is nondegenerate. 
\end{cor}

\subsection{Quantum current and Feigin-Odesskii algebras} \label{QC}

Our next results deal with quantum current algebras. Assume that 
the Cartan matrix $A$ is of finite type. Let $L\N_+$  be
the current Lie algebra $\N_+\otimes \CC[t,t^{-1}]$, endowed with the
bracket $[x\otimes t^n,y\otimes t^m] = [x,y]\otimes t^{n+m}$. 

\subsubsection{Quantum affine algebras}

Let $\cA$ be the quotient of the free algebra
$\CC[[\hbar]]\langle e_i[k],  i = 1, \ldots,n, k\in\ZZ\rangle$ by the
two-sided ideal generated by the coefficients  of monomials in the
formal series identities  
\begin{equation} \label{crossed:vertex}
(q^{d_i a_{ij}}z - w)e_i(z) e_j(w)  = (z - q^{d_i a_{ij}}w) e_j(w)
e_i(z) ,  
\end{equation} 
\begin{equation} \label{q:serre:nr}
\Sym_{z_1,\ldots,z_{1-a_{ij}}} \sum_{k=0}^{1-a_{ij}} (-1)^k \bmatrix 1 -
a_{ij} \\ k \endbmatrix_{q^{d_i}} e_i(z_1) \ldots e_i(z_{k}) e_j(w)
e_i(z_{k+1}) \ldots e_i(z_{1-a_{ij}})  =0,  
\end{equation} 
where $e_i(z)$ is the generating series $e_i(z) = \sum_{k\in z}
e_i[k]z^{-k}$, and $q = e^\hbar$.

Let $U_\hbar L\N_+$ be the quotient $\cA / (\cap_{N>0}\hbar^N \cA)$. 

Define $\wt A$ as the quotient of the free algebra $\CC[[\hbar]]\langle 
e_i[k]^{\wt \cA}, i = 1, \ldots, n k\in\ZZ\rangle$ by the two-sided ideal 
geberated by the coefficients of monimials of (\ref{crossed:vertex})
and 
\begin{equation} \label{variant:serre}
\sum_{k = 0}^{1 - a_{ij}} (-1)^k \bmatrix 1 - a_{ij} \\ k
\endbmatrix_{q^{d_{i}}} (e_i[0]^{\wt \cA})^k e_j[l]^{\wt \cA}
(e_i[0]^{\wt \cA})^{1 - a_{ij} - k} = 0 ,
\end{equation}
for any $i,j = 1,\ldots,n$ and $l$ integer. Define 
$\wt U_\hbar L\N_+$ as the quotient 
$\wt\cA / (\cap_{N>0}\hbar^N \wt\cA)$.

\begin{thm} \label{thm:third} 
1) $U_\hbar L\N_+$ is a free $\CC[[\hbar]]$-module, and the map
$e_i[n]\mapsto  \bar e_i\otimes t^n$ defines an algebra isomorphism from
$U_\hbar L\N_+ / \hbar U_\hbar L\N_+$ to  $U L\N_+$.  

2) Let $U_\hbar L\N_+^{top}$ be the quotient of $\CC\langle  e_i[k], i =
1,\ldots,n,k\in\ZZ\rangle[[\hbar]]$ by the $\hbar$-adically closed
two-sided ideal generated by the coefficients of monomials in  relations
(\ref{crossed:vertex}) and (\ref{q:serre:nr}). Then  $U_\hbar
L\N_+^{top}$ is a topologically free $\CC[[\hbar]]$-module;  it is
naturally the $\hbar$-adic completion of  $U_\hbar L\N_+$, and the map 
$e_i[n]\mapsto  \bar e_i\otimes t^n$ defines an algebra isomorphism 
from $U_\hbar L\N_+^{top} / \hbar U_\hbar L\N_+^{top}$ to  $U L\N_+$.  

3) There is a unique algebra map from $\wt U_\hbar L\N_+$ to $U_\hbar L\N_+$, 
sending each $e_i[k]^{\wt \cA}$ to $e_i[k]$; it is an algebra isomorphism. 

4) Let $\wt U_\hbar L\N_+^{top}$ be the quotient of $\CC\langle 
e_i[k]^{\wt \cA}, i = 1,\ldots,n,k\in\ZZ\rangle[[\hbar]]$ by the $\hbar$-adically closed
two-sided ideal generated by the coefficients of monomials in 
relations (\ref{crossed:vertex}) and (\ref{variant:serre}). Then 
$e_i[k]^{\wt \cA} \mapsto e_i[k]$ defines an algebra isomorphism between 
$\wt U_\hbar L\N_+^{top}$ and  $U_\hbar L\N_+^{top}$. 
\end{thm} 

The statements 1) and 2) of this Theorem can be derived from the results
of \cite{Beck}. 

In \cite{FO1,FO2}, Feigin and Odesskii defined the algebra $FO$, which may be
viewed  as a functional version of the shuffle algebra. $FO$ is defined
as  
\begin{equation} \label{grading} 
FO = \oplus_{\kk\in\NN^n} FO_\kk, 
\end{equation} 
where if $\kk = (k_i)_{1\leq i\leq n}$, we set 
$$
FO_\kk =
{1\over{\prod_{i<j} \prod_{1\leq\al\leq k_i, 1\leq\beta\leq k_j}
( t^{(i)}_\al  - t^{(j)}_\beta ) }}\CC[[\hbar]][(t^{(i)}_j)^{\pm 1}
, i = 1,\ldots,n, j  = 1,\ldots,k_i]^{\gotS_{k_1}\times\cdots\times
\gotS_{k_n}} ,  
$$  
where the product of symmetric groups acts by permutation of variables
of each  group of variables $(t^{(i)}_j)_{1\leq j \leq k_i}$.  $FO_\kk$ therefore
consists of rational functions in the $t^{(i)}_j$, symmetric in each
group  $(t^{(i)}_j)_{1\leq j \leq k_i}$, regular except for poles when the variables go to 
$0$ or infinity, and simple poles when variables of different ``colors'' collide.   
(\ref{grading}) defines a grading of $FO$ by $\NN^n$.
The product on $FO$ is also graded, and we have, for $f$ in $FO_\kk$ and $g$ in
$FO_\bl$ ($\bl = (l_i)_{1\leq i\leq n}$) ,  
\begin{align} \label{pdt:FO}
& (f*g)(t^{(i)}_j)_{1\leq i\leq n, 1\leq j \leq k_i + l_i}  \\ & \nonumber 
=  \Sym_{\{ t^{(1)}_j \}} 
 \ldots \Sym_{\{ t^{(n)}_j \}}  [\prod_{1\leq i\leq N} \prod_{N+1\leq
j\leq N+M} {{q^{ \langle \eps(i), \eps(j)\rangle} t_i - t_j}\over{t_i - t_j}}
f(t_1,\ldots,t_N) g(t_{N+1},\ldots,t_{N+M})] ,  
\end{align}
where $N = \sum_{i=1}^n k_i$ and $M = \sum_{i=1}^n l_i$; we set  
for any $s$, 
$$
t_{k_1 + \cdots + k_{s-1} + 1} = t^{(s)}_1 , 
\ldots , 
t_{k_1 + \cdots + k_{s}} = t^{(s)}_{k_s} , 
\quad 
t_{N+l_1 + \cdots + l_{s-1} + 1} = t^{(s)}_{k_s + 1} , \ldots , 
t_{N+l_1 + \cdots + l_s} = t^{(s)}_{k_s + l_s} ; 
$$
we also define
$\eps(i) = \eps_k$ if $t_i = t^{(k)}_l$ for some $l$;  
as before, $\langle \eps_i, \eps_j\rangle = d_i a_{ij}$ for $i,j =
1,\ldots,n$. 

In the right side of (\ref{pdt:FO}), each symmetrization can be
replaced by a sum over shuffles, since the argument in symmetric in 
each group of variables $(t^{(s)}_1,\ldots, t^{(s)}_{k_s})$ and
$(t^{(s)}_{k_s + 1},\ldots, t^{(s)}_{k_s + l_s})$. 

In Prop.\ \ref{Hopf:FO}, we define a topological braided Hopf structure on 
$FO$.  

In \cite{Enr}, we showed: 

\begin{prop}
There is a unique algebra morphism $i_\hbar$ from $U_\hbar L\N_+$ to $FO$, 
such that $i_\hbar(e_i[n])$ is the element $(t^{(i)}_1)^n$ of 
$FO_{\eps_i}$. 
\end{prop}

Let us denote by $LV$ the direct sum $\oplus_{i=1}^n FO_{\eps_i}$ and let 
$\langle LV\rangle$ be the sub-$\CC[[\hbar]]$-algebra of $FO$ generated by 
$LV$. As a corollary of the proof of Thm.\ \ref{thm:third}, we prove: 

\begin{cor} \label{cor:second}
$i_\hbar$ is an algebra isomorphism between $U_\hbar L\N_+$ and 
$\langle LV\rangle$. 
\end{cor}

$U_\hbar L\N_+$ is also endowed with a topological braided Hopf structure
(the Drinfeld comultiplication); it is then easy to see that $i_\hbar$ is 
compatible with both Hopf structures. 

Define $T(LV_\pm)$ as the free algebras
$\CC[[\hbar]]\langle e_i[k]^{(T)}, 1\leq i\leq n, k\in\ZZ \rangle$   and
$\CC[[\hbar]]\langle f_i[k]^{(T)}, 1\leq i\leq n, k\in\ZZ \rangle$. 
We
have a pairing $\langle , \rangle_{ T(LV_\pm) }$  between
$T(LV_+)$ and $T(LV_-)$ defined by
\begin{align} \label{pairing:introd}
& \langle e_{i_1}[k_1]^{(T)} \cdots e_{i_p}[k_p]^{(T)} , 
f_{j_1}[l_1]^{(T)} \cdots f_{j_{p'}}[l_{p'}]^{(T)}
\rangle_{T(LV_\pm)} = \delta_{pp'} {1\over{\hbar^p}} \cdot \\ & \nonumber 
\cdot \sum_{\sigma\in \gotS_p} \res_{z_1 = 0} \cdots \res_{z_p = 0}
(\prod_{s>t,\sigma^{-1}(s) < \sigma^{-1}(t)} { { q^{(\eps_{i_s}, \eps_{i_t})} 
z_s - z_t} \over{ z_s - q^{(\eps_{i_s}, \eps_{i_t})} z_t}}
\prod_{s=1}^p {1\over {d_{j_s}}}\delta_{i_{\sigma(s),j_s}} \prod_{s=1}^p
z_i^{k_s + l_s} \prod_{s=1}^p {{dz_s}\over{z_s}})
\end{align}
where the ratios ${ { q^{(\eps_{i_s}, \eps_{i_t})} 
z_s - z_t}\over{ z_s - q^{(\eps_{i_s}, \eps_{i_t})} z_t }}$ are expanded 
for $z_t << z_s$. 

Let $U_\hbar L\N_-$ be the quotient of $T(LV_-)$ by the homomorphic image of the 
ideal defining $U_{-\hbar} L\N_+$ by the map $e_i[k]^{(T)} \mapsto f_i[k]^{(T)}$. 

\begin{prop} \label{prop:morphism} (see \cite{Enr})  
This pairing induces a pairing $\langle , \rangle_{U_\hbar L\N_\pm}$ between
$U_\hbar L\N_+$ and $U_\hbar L\N_-$. 
\end{prop}

We then prove: 
\begin{thm} \label{nondeg:L}
The pairing  $\langle, \rangle_{U_\hbar L\N_\pm}$ is nondegenerate. 
\end{thm}
 
\subsubsection{The form of the $R$-matrix}

Let us set $A_\pm = U_\hbar L\N_\pm$.  Let $a$ and $b$ be two integers. 
Define $I^+_{\geq a}$ and $I^+_{\leq a}$, as the right, resp.\ left
ideals of  $A_+$ generated by the  $e_i[k],k\geq a$ , resp.\
the $e_i[k],k\leq a$.  The ideals $I^+_{\geq a}$  and  $I^+_{\leq a}$
are graded; for $\al$ in $(\pm\NN)^n$, denote by $I^+_{\geq a}[\al]$ 
and  $I^+_{\leq a}[\al]$  their component of degree $\al$. 

\begin{prop} \label{R:mat:QC} 
For any $\al$ in $\NN^n$, for any integers $a$ and $b$, $(I^+_{\leq a} +
 I^+_{\geq b})^\perp[-\al]$ and $[A_+ / (I^+_{\leq a} +  I^+_{\geq
b})^{\perp\perp}][\al]$ are free finite-dimensional
$\CC[[\hbar]]$-modules.  The pairing between $A_+$ and $A_-$ induces a
nondegenerate pairing between them. Moreover, the intersection
$\cap_{a,b} (I^+_{\leq a} +  I^+_{\geq b})^{\perp\perp}$ is zero. 
 
Denote by $P_{a,b}[\al]$ the corresponding element of  $[A_+ /
(I^+_{\leq a} +  I^+_{\geq b})^{\perp\perp}][\al] \otimes  (I^+_{\leq a}
+  I^+_{\geq b})^\perp[-\al][\hbar^{-1}]$.   $P_{a,b}[\al]$ defines an
element of $\limm_{\leftarrow a,b} A_+ / (I^+_{\leq a} +  I^+_{\geq
b})^{\perp\perp} \otimes_{\CC[[\hbar]]}  A_-[\hbar^{-1}]$. 

Let $(\bar e_{\beta,i})_{\beta\in\Delta_+}$ and $(\bar
f_{\beta,i})_{\beta\in\Delta_+}$ be dual Cartan-Weyl bases of $\N_+$ and
$\N_-$. Let $e_{\beta,i}[p]$ and $f_{\beta,i}[p]$ be lifts to $U_\hbar L\N_\pm$
of $\bar e_{\beta,i}\otimes t^p$ and $\bar f_{\beta,i}\otimes t^p$. 

Then if $\al$ belongs to 
$k\Delta_+ - (k-1)\Delta_+$, $P[\al]$ has the form
\begin{align*}
& P[\al] = {{\hbar^k}\over{k!}}
\sum_{\al_1,\ldots,\al_k\in\Delta_+,\sum_i\al_i = \al; i_j}
\sum_{p_1,\ldots,p_k\in \ZZ}  e_{\al_1,i_1}[p_1] \ldots e_{\al_k,i_k}[p_k]
\otimes  f_{\al_1,i_1}[-p_1] \ldots f_{\al_k,i_1}[-p_k]  \\ & + o(\hbar^k) 
\end{align*}
(all but a finite number of elements of this sum belong to  $(I^+_{\leq
a} +  I^+_{\geq b})^{\perp\perp} \otimes_{\CC[[\hbar]]} 
A_-[\hbar^{-1}]$.) 
\end{prop}

Let $(h'_i)_{i = 1,\ldots,n}$ be the basis of $\HH$, dual to 
$(h_i)_{i = 1,\ldots,n}$. Set $\cK = \exp(\sum_{i = 1}^n 
h_i[0] \otimes h'_i[0] + \sum_{p>0}
h_i[p] \otimes h'_i[-p]) $. Then the elements $\cR[\al] = \cK P[\al]$ 
of $\limm_{\leftarrow N} 
(U_\hbar L\B_+\otimes U_\hbar L\B_-) / I_N^{\B_\pm,(2)}$, 
where $I_N^{\B_\pm,(2)}$ is the ideal generated by the $h_i[p]\otimes 1, 
e_i[p]\otimes 1, p \geq N$, and the $1\otimes f_i[p], p\geq N$, 
satisfy the $R$-matrix identity
$$
\sum_{\gamma \in\NN^n, \beta \in (\pm \NN)^n, \beta + \gamma = \la }
\cR[\gamma] \Delta(x)_{(\beta,\al - \beta)} 
= 
\sum_{\gamma \in\NN^n, \beta \in (\pm \NN)^n, \beta + \gamma = \la }
\Delta'(x)_{(\beta,\al - \beta)} \cR[\gamma] , 
$$
for any $\la\in \ZZ^r$ and $x$ in the double $U_\hbar L\G$ of  $U_\hbar
L\B_+$ of degree $\al$ (the sums over the root lattice are obviously
finite, and each product makes sense in $\limm_{\leftarrow N}
(U_\hbar L\G \otimes U_\hbar L\G) / I_N^{\G,(2)}
$, where $I_N^{\G,(2)}$ is the left ideal generated by the $x[p]\otimes 1$ and 
$1\otimes x[p], p\geq N$, $x = e_i,h_i,f_i$).





\subsubsection{Yangians}

Let us describe how the above results are modified in the case of Yangians. 
Let $\cA^{rat}$ be the quotient of the free algebra 
$\CC[[\hbar]]
\langle e_i[k]^{rat}
, i = 1,\ldots,n,k\in\ZZ \rangle$ by the 
two-sided ideal generated by the coefficients of the relations
\begin{equation} \label{crossed:vx:rat}
(z-w + \hbar a_{ij}) e_i(z)^{rat} e_j(w)^{rat} =    
(z-w - \hbar a_{ij}) e_j(w)^{rat} e_i(z)^{rat},  
\end{equation}
\begin{equation} \label{yg:serre}
\Sym_{z_1,\ldots,z_{1  - a_{ij}}} \ad(e_i(z_1)^{rat}) \cdots 
\ad(e_i(z_{1 - a_{ij}})^{rat}) (e_j(w)^{rat}) = 0,  
\end{equation}
where we set $e_i(z)^{rat} = \sum_{k\in\ZZ} e_i[k]^{rat} z^{-k-1}$, and let us set 
$U_\hbar^{rat}L\N_+ = \cA^{rat} / \cap_{N>0} \hbar^N \cA^{rat}$.  
 
Define also $\wt \cA^{rat}$ as the quotient of the free algebra 
$\CC[[\hbar]]
\langle e_i[k]^{\wt\cA^{rat}}
, i = 1,\ldots,n,k\in\ZZ \rangle$ by the 
two-sided ideal generated by the coefficients of the relations
(\ref{crossed:vx:rat}) and 
\begin{equation} \label{yg:variant:serre}
\ad (e_i[0]^{\wt\cA^{rat}})^{1 - a_{ij}} e_j[l]^{\wt\cA^{rat}} = 0, 
\end{equation}
for any $i,j = 1,\ldots, n$ and integer $l$. 

\begin{thm}
1) $U_\hbar L\N_+$ is a free $\CC[[\hbar]]$-module. There is a unique 
algebra isomorphism from $U_\hbar^{rat} L\N_+ /\hbar U_\hbar^{rat} L\N_+$ to 
$U L\N_+$, sending the class of $e_i[k]^{rat}$ to $\bar e_i\otimes t^k$.  

2) There is a unique algebra morphism from $\wt U_\hbar L\N_+$ to 
$U_\hbar L\N_+$, sending $e_i[k]^{rat}$ to $e_i[k]^{\wt\cA^{rat}}$; 
it is an isomorphism between these algebras. 

3) Let $U_\hbar L\N_+^{rat,top}$ and  $\wt U_\hbar L\N_+^{rat,top}$ be
the  quotients of $\CC\langle  e_i[k]^{rat}, i =
1,\ldots,n,k\in\ZZ\rangle[[\hbar]]$  and  $\CC\langle 
e_i[k]^{\wt\cA^{rat}}, i = 1,\ldots,n,k\in\ZZ\rangle[[\hbar]]$  by the
$\hbar$-adically closed two-sided ideals generated by the coefficients
of monomials in  relations (\ref{crossed:vx:rat}) and
(\ref{yg:serre}), resp.\  (\ref{crossed:vx:rat}) and
(\ref{yg:variant:serre}). Then  $e_i[k]^{\wt \cA^{rat}} \mapsto
e_i[k]^{rat}$ defines an algebra isomorphism between  $U_\hbar
L\N_+^{rat,top}$ and  $\wt U_\hbar L\N_+^{rat,top}$.  $U_\hbar
L\N_+^{rat,top}$ is the $\hbar$-adic completion of $U_\hbar
L\N_+^{rat}$;  it is a topologically free $\CC[[\hbar]]$-module. 
\end{thm}

Define $FO^{rat}$ as the graded space $FO$, endowed with the product
obtained from (\ref{pdt:FO}) by the replacement of each $q^{\la}z -
q^{\mu}w$ by $z - w+\hbar(\la - \mu)$. $FO^{rat}$ is an associative algebra 
and we have 

\begin{thm}
There is a unique algebra map $i_\hbar^{rat}$ from $U_\hbar L\N_+^{rat}$ to 
$FO^{rat}$, sending each $e_i[k]^{rat}$ to $e_i\otimes t^k$. $i_\hbar$
is an isomorphism between $U_\hbar L\N_+^{rat}$ and its subalgebra 
$\langle LV\rangle^{rat}$ generated by the degree one elements. 
\end{thm}

Define a pairing $\langle , \rangle_{T(LV_\pm), rat}$ between $T(LV_+)$ and 
$T(LV_-)$ by the formula (\ref{pairing:introd}), where each
$q^{\la} z - q^{\mu}w$ is replaced by $z - w + \hbar(\la - \mu)$. 
Let $U_\hbar L\N_-^{rat}$ be the quotient $T(LV_-)$ by the homomorphic 
image of the ideal defining $U_\hbar L\N_+^{rat}$ by the map $e_i[k]^{(T)} 
\mapsto f_i[k]^{(T)}$. 

\begin{thm}
$\langle , \rangle_{T(LV_\pm), rat}$ induces a nondegenerate pairing 
$\langle , \rangle_{U_\hbar L\N_\pm, rat}$ between $U_\hbar L\N_\pm^{rat}$.

Define $I_a^{\pm, rat}$ as the right, resp.\ left ideals of $U_\hbar L\N_\pm^{rat}$
generated by the $e_i[k]^{rat}$, $k\geq a$, resp.\ by the $f_i[k]^{rat}$, $k\geq a$. 
Then for any $\al$ in $\NN^n$, $I^{+,rat}_a[\al] \cap (I^{-,rat}_b[-\al])^\perp$ is a 
space of finite codimension in $(I^{-,rat}_b[-\al])^\perp$. 
Let $P[\al]^{rat}$ be the corresponding element of 
$\limm_{\leftarrow a,b} (U_\hbar L\N_+^{rat} / I_a^{+,rat})[\al] 
\otimes (U_\hbar L\N_-^{rat} / I_b^{-,rat})[-\al][\hbar^{-1}]$. 
Let $e_{\beta,i}[p]^{rat}$ be lifts to $U_\hbar L\N_\pm$ of the 
$\bar e_{\beta,i}[p]$. 

Then if $\al$ belongs to 
$k\Delta_+ - (k-1)\Delta_+$, $P[\al]^{rat}$ has the form
\begin{align*}
& P[\al]^{rat} = {{\hbar^k}\over{k!}}
\sum_{\al_1,\ldots,\al_k\in\Delta_+,\sum_i\al_i = \al; i_j}
\sum_{p_1,\ldots,p_k\in \ZZ}  \\ & e_{\al_1,i_1}[p_1]^{rat} \ldots
e_{\al_k,i_k}[p_k]^{rat} \otimes  f_{\al_1,i_1}[-p_1-1]^{rat} \ldots
f_{\al_k,i_1}[-p_k - 1]^{rat}  + o(\hbar^k). 
\end{align*}
\end{thm}

The proofs of the statements of this section are analogous to those of 
the quantum affine case and will be omitted.

\subsection{Quantum current algebras of affine type (toroidal 
algebras) } \label{QC:affine}

Assume that $A$ is an arbitrary symmetrizable Cartan matrix.   Define
$U_\hbar L\N_+$  and $\wt U_\hbar L\N_+$ as in sect.\ \ref{QC}.  

\begin{prop} \label{oz}
1) Let $\wt F_+$ be the Lie algebra with generators $\wt x_i^+[k], 1\leq i \leq n, 
k\in\ZZ$, and relations given by the coefficients of monomials in 
$$
(z - w) [\wt x_i^+(z), \wt x_j^+(w)] = 0, \quad 
\ad (\wt x_i^+(z_1))  \cdots \ad (\wt x_i^+(z_{1 - a_{ij}})) 
(\wt x_j^+(w)) = 0,  
$$
where for $x$ any $\wt x_i^+$, $x(z)$ is the generating series
$\sum_{k\in\ZZ} x[k]z^{-k}$.  If we give degree $\eps_i$ to $\wt x_i^+[k]$, 
$\wt F_+$ is graded by set $\Delta_+$ of the roots of $\N_+$. 
 
Then $U_\hbar L\N_+ / \hbar U_\hbar L\N_+$  and $\wt U_\hbar L\N_+ / \hbar \wt
U_\hbar L\N_+$  are both isomorphic to the  enveloping algebra $U(\wt
F_+)$. 

2) There is a unique Lie algebra morphism  $j_+: \wt F_+ \to L\N_+$,  such
that $j_+(\wt e_i^+[k]) = \bar e_i\otimes t^k$. $j_+$ is graded and
surjective.  The kernel of $j_+$ is contained in $\oplus_{\al\in \Delta_+,
\al\ \on{imaginary}} \wt F_+[\al]$. 
\end{prop}

Let us assume now that $A$ is untwisted affine. $\N_+$ is the isomorphic
to a subalgebra of the loop algebra $\bar\G[\la,\la^{-1}]$, with $\bar\G$ semisimple.
Define $\T_+$ as the direct sum $L\N_+ \oplus (\oplus_{k>0,l\in\ZZ}
\CC K_{k\delta}[l])$, and endow it with the bracket such that the 
$K_{k\delta}[l]$ are central, and
$$
[(x \otimes t^l,0), (y \otimes t^m,0)] = ([x \otimes t^l, y\otimes t^m], 
\langle \bar x, \bar y\rangle_{\bar\G} (lk'' - m k') K_{(k' + k'')\delta}[l+m] ) 
$$ 
for $x \mapsto \bar x \otimes \la^{k'}, y \mapsto \bar y \otimes \la^{k''}$ 
by the 
inclusion $\N_+ \subset \bar\G[\la,\la^{-1}]$, where
$\langle , \rangle_{\bar\G}$ is an invariant scalar product on 
$\bar\G$. 

Then $\T_+$ is a Lie subalgebra of the 
toroidal algebra $\T$, which is the universal central
extension of $L\G$ (\cite{Kassel,Moody}). In what follows, we 
will set  $x[k]^{\T} = (x\otimes t^k,0)$.

\begin{prop} \label{nagila}
1) When $A$ is of affine Kac-Moody type, the kernel of $j_+$ is equal to 
the center of $\wt F_+$, so that $\wt F_+$ is a central extension of
$L\N_+$. 

2) We have a unique Lie algebra map $j'$ from $\T_+$ to $\wt F_+$ such
that $j'(\bar e_i \otimes t^n) = \wt e_i[n]$. This map is an isomorphism
iff $A$ is not of type $A_1^{(1)}$.  If $A$ is of type $A_1^{(1)}$, $j'$
is surjective, and its kernel is $\oplus_{n\in\ZZ} \CC K_{\delta}[n]$. 
\end{prop}

In Remark \ref{rem:generalizations}, we discuss possible generalizations
of Thm.\ \ref{thm:third} to the case of affine quantum current 
algebras, and the connection Prop.\ \ref{oz} with the results of
\cite{FO:new}. 

\subsection{}

The basic idea of the constructions of the two first parts of this work 
is to  compare the quantized algebras defined by generators and
relations   with quantum shuffles algebras. The idea to use shuffle
algebras to provide examples of Hopf algebras dates back to  Nichols
(\cite{Nichols}). Later, Schauenburg (\cite{Schau}) and Rosso 
(\cite{Rosso:shuffle}) showed that the positive part $U_\hbar\N_+$ of
the Drinfeld-Jimbo  quantized enveloping algebras are isomorphic to the
subalgebra $\Sh(V)$ of quantum shuffle Hopf algebras generated in degree
$1$.  Their results rely on Lusztig's PBW or duality (nondegeneracy of
Drinfeld's pairing) results. A nonabelian generalization of Schauenburg's result can be found 
in \cite{AG}.

In sect.\  \ref{proof:KM}, we show that applying Drinfeld's theory of Lie 
bialgebras  to $\Sh(V)$ yields  at the same time
proofs of these results (PBW for $U_\hbar\N_+$ and isomorphism of  
$U_\hbar\N_+$ with $\Sh(V)$, and nondegeneracy of the pairing as a simple
consequence), when the deformation parameter is formal or generic. 

In sect.\ \ref{proof:QC}, we apply the same idea to quantum current 
algebras. These algebras, also know as ``new realizations'' 
algebras, depend on the datum of a Cartan matrix.
In that situation, the proper replacement of shuffle  algebras
are the functional shuffle algebras introduced by  Feigin and Odesskii
(\cite{FO1,FO2}). We show that when the Cartan matrix is of 
finite type, the ideas of sect.\ \ref{proof:KM}  allow
to complete the results of \cite{Enr} on comparison of the quantum
current algebras  and the Feigin-Odesskii algebras. However, there 
are still some open problems in this direction, see Rem.\ \ref{open}. 
We hope that the ideas of this section will help generalize
the  results of \cite{Enr:Rub} from $\SL_2$ to arbitrary semisimple Lie
algebras. For this one should, in particular, find analogues of the
quantum Serre relations for the algebras in  genus $\geq 1$. 

In sect.\ \ref{toroidal}, we consider the classical limit of the 
quantum current algebras in the case of an affine Kac-Moody  Cartan
matrix. We show that this  classical limit is the enveloping algebra of
a Lie algebra  $\wt F_+$, which is a central extension of the Lie
algebra $L\N_+$ of  loops with values in the positive subalgebra $\N_+$
of the  affine Kac-Moody algebra.  $\wt F_+$ is graded by the roots of
$\N_+$, and its center is contained in the part of imaginary degrees. 
We show that in all affine untwisted cases, except the $A_1^{(1)}$ case,
 $\wt F_+$ is isomorphic to a subalgebra $\T_+$ of the toroidal
algebra $\T$ introduced and studied  in 
\cite{Moody,GKV,Vass}. In the 
$A_1^{(1)}$ case, we identify $\wt F+$ with a quotient of $\T_+$. 

In the quantum case, the center of  $U_\hbar L\N_+$ seems closely
connected  with the central part of the affine  elliptic algebras
constructed in the recent work of Feigin and Odesskii  (\cite{FO:new}).
We hope that a better understanding of this center will  enable to
extend to toroidal algebras the results of  sect.\ \ref{proof:QC}. 

This paper grew from the notes of the DEA course I taught at  univ.\ Paris
6 in february-april 1999. I would like to thank P.\ Schapira for giving
me the opportunity to give this course and its participants, notably C.\
Grunspan, O.\ Schiffmann and V.\ Toledano Laredo, for their patience and
attention. 

I also would like to thank N.\ Andruskiewitsch and B.\ Feigin for
valuable discussions. In particular, the idea that quantum shuffle algebras 
could be a tool to construct quantizations of Lie bialgebras is due to
N.\ Andruskiewitsch; it is clear that this idea plays an important role
in the present work. I also would like to thank N.\ Andruskiewitsch 
for his kind invitation to the univ.\ of C\'ordoba in August 1998, where these
discussions took place.

\section{Quantum Kac-Moody algebras (proofs of the results of sect.\
\ref{KM})} \label{proof:KM}

\subsection{PBW theorem and comparison with shuffle algebra
(proofs of Thm.\ \ref{thm:first} and Cor.\ \ref{cor:comparison})}

\subsubsection{Definition of $\Sh(V)$ and $\langle V \rangle$} \label{sect:shuffles}

Let us set $V = \oplus_{i=1}^n \CC v_i$. Let $\eps_i$ be the $i$th basis
vector of $\NN^n$.  Define the grading of $V$ by $\NN^n$ by $\deg(v_i) =
\eps_i$. Let $\Sh(V)$ be the quantum shuffle  algebra constructed from
$V$ and the braiding  $V\otimes V \to V\otimes V[[\hbar]]$, $v_i \otimes
v_j \mapsto q^{ - d_i a_{ij}} v_j \otimes v_i$. That is,  $\Sh(V)$ is
isomorphic, as $\CC[[\hbar]]$-module, to $\oplus_{i\geq 0} V^{\otimes
i}[[\hbar]]$. Denote the element $z_1\otimes\cdots\otimes z_k$ as
$[z_1|\cdots|z_k]$. The product is defined on $\Sh(V)$ as follows:  
\begin{align*}
[z_1|\cdots|z_k] 
& [z_{n+1}|\cdots|z_{k+l}] \\ & = \sum_{\si\in \Sigma_{k,l}}
q^{ - \sum_{1\leq i < j\leq k+l, \si(i) > \si(j)} \langle \deg(z_i)
, \deg(z_j) \rangle} [z_{\si(1)}|\cdots|z_{\si(k+l)}] ,    
\end{align*}
if the $z_i$ are homogeneous elements of $V$, and  where $\Sigma_{k,l}$ is
the subset of the symmetric group $\gotS_{k+l}$ consisting of  shuffle
permutations $\sigma$ such that $\si(i)< \si(j)$ if $1\leq i < j \leq k$
 or $k+1\leq i < j \leq k+l$; the bilinear form on $\NN^n$ is defined by 
form  $\langle \eps_i,\eps_j \rangle = d_i a_{ij}$.  

\begin{lemma} \label{<V>:free}
$\langle V \rangle$  is the direct sum of its  graded components, which
are free $\CC[[\hbar]]$-modules. It follows that $\langle V \rangle$  is a free
$\CC[[\hbar]]$-module. 
\end{lemma}

{\em Proof.} That $\langle V \rangle$  is the direct sum of its  graded
components follows from its definition. These graded components are 
$\CC[[\hbar]]$-submodules of finite-dimensional free
$\CC[[\hbar]]$-modules (the graded components of  $\Sh(V)$). Each graded
component is therefore a   finite-dimensional free module over 
$\CC[[\hbar]]$. The Lemma follows. \hfill \qed \medskip

\subsubsection{Crossed product algebras $\cV$ and $\cS$}

Define linear endomorphisms $\wt h_i,i= 1,\ldots,n$ and 
$\wt D_j, j = 1,\ldots, n-r$ of $V$ by the formulas 
$$
\wt h_i(v_j) = a_{ij}v_j, \quad \wt D_j(v_i) = \delta_{ij}v_i.  
$$
Extend the $\wt x$, $x\in \{h_i,D_j\}$ to linear endomorphisms of $\Sh(V)$ by the formulas  
$$
\wt x([x_1|\cdots|x_n]) = \sum_{k=1}^n [x_1|\cdots|\wt x(x_k)|\cdots |x_n].   
$$ 
It is clear that the $\wt x$ define derivations of $\Sh(V)$. These 
derivations preserve $\langle V \rangle$. 

Define $\cV$ and $\cS$ as the crossed product algebras of $\langle V
\rangle$ and  $\Sh(V)$ with the derivations $\wt h_i,\wt D_j$. More precisely,
$\cV$ and $\cS$ are  isomorphic, as $\CC[[\hbar]]$-modules, to their
tensor products  $\langle V\rangle
\otimes_{\CC[[\hbar]]}\CC[h^{\cV}_i,D_j^{\cV}][[\hbar]]$ and  $\Sh(V)
\otimes_{\CC[[\hbar]]}\CC[h^{\cV}_i,D_j^{\cV}][[\hbar]]$ with $\hbar$-adically
completed polynomial
algebras in $n+r$ variables.  The products on $\cV$ and $\cS$ are then
defined by the rules
$$
(x\otimes \prod_{s = 1}^{2n-r} (X_s^{\cV})^{\al_s})  
(y\otimes \prod_{s = 1}^{2n-r} (X_s^{\cV})^{\beta_s})  
= 
\sum_{(i_s)} \prod_{s = 1}^{2n-r} 
\pmatrix \al_s \\ i_s \endpmatrix
\left( x \prod_{s = 1}^{2n-r} \wt X_s^{i_s} (y) \right) \otimes \left( 
\prod_{s = 1}^{2n-r}  
(X_s^{\cV})^{\al_s + \beta_s - i_s} \right) ,   
$$
where we set $X_s = h_s$ for $s = 1,\ldots , n$, and 
$X_s = D_{s-n}$ for $s = n+1,\ldots, 2n-r$. 
In what follows, we will denote $x\otimes 1$ and $1\otimes X_s$ simply by 
$x$ and $X_s$, so that $x\otimes \prod_s (X_s^\cV)^{\al_s}$ will be 
$x\prod_s (X_s^\cV)^{\al_s}$. 

$\cS$ is then endowed with a Hopf $\CC[[\hbar]]$-algebra structure (that is,
all maps of Hopf algebra axioms are $\CC[[\hbar]]$-module maps, and the tensor
products are completed in the $\hbar$-adic topology),
defined by
$$
\Delta_\cV(h^{\cV} ) =  h^{\cV} \otimes 1 + 1 \otimes h^{\cV} \on{\ for\ } 
h\in \{h_i,D_j\},   
$$
$$
\Delta_\cV([v_{i_1} | \cdots | v_{i_m} ])  = \sum_{k=0}^m [v_{i_1} |
\cdots | v_{i_k} ] \otimes  \exp( \hbar \sum_{j=1}^{k} d_{i_j}
h_{i_j}^{\cV})[v_{i_{k+1}} | \cdots | v_{i_m} ]) . 
$$
$\cV$ is then a Hopf subalgebra of $\cS$. 

Assign degrees $0$ to the elements  $h_i^{\cV},D_j^{\cV}$, and $\eps_i$
to $v_i$. $\cV$ is then the  direct sum of its homogeneous components,
which are free  finite-dimensional modules over
$\CC[h_i^{\cV},D_j^{\cV}][[\hbar]]$; the grading of $\cV$ is compatible
with its algebra structure. 

\subsubsection{Hopf co-Poisson and Lie bialgebra structures} \label{adama}

Define $\cV_0$ as $\cV / \hbar\cV$. 

\begin{lemma}
$\cV_0$ is a cocommutative Hopf
 algebra.
\end{lemma}

{\em Proof.} Define $\Delta'_{\cV}$ as $\Delta_{\cV}$ composed 
with the permutation of factors. 
 We have to show that for $x$ in $\cV$, we have 
\begin{equation} \label{balus}
(\Delta_{\cV} - \Delta'_{\cV})(x)
\subset \hbar(\cV\otimes_{\CC[[\hbar]]}\cV).
\end{equation} 
For $x$ one of the $\bar h_i,\bar D_j$, (\ref{balus}) is clearly 
satisfied. On the other hand, if (\ref{balus}) is satisfied for $x$ and $y$
in $\cV$, then $(\Delta_{\cV} - \Delta'_{\cV})(xy)$ is equal to 
$(\Delta_{\cV} - \Delta'_{\cV})(x)\Delta_{\cV}(y) + 
\Delta'_{\cV}(x) (\Delta_{\cV} - \Delta'_{\cV})(y)$ and therefore belongs to 
$\hbar(\cV\otimes_{\CC[[\hbar]]}\cV)$. It follows that (\ref{balus}) 
holds for any $x$ in $\cV$. \hfill\qed\medskip

\begin{lemma} \label{pourim} 

1) There exists a unique surjective Hopf algebra morphism $p_\hbar$ 
from   $U_\hbar\B_+$ to $\cV$, such that $p_\hbar(h_i) = h_i^{\cV}$ and
$p_\hbar(x_i^+) =  [v_i]$. 

2) The map $D_j\mapsto \bar D_j,h_i\mapsto \bar h_i,x_i^+\mapsto \bar x_i^+$
 extends to an isomorphism from $U_\hbar\B_+ / \hbar U_\hbar\B_+$ to
$U\B_+$.  

3) $p_\hbar$ induces a surjective cocommutative Hopf algebra morphism 
$p$ from  $U_\hbar\B_+ = U\B_+$ to $\cV / \hbar\cV = \cV_0$. 
\end{lemma}

{\em Proof.} That the quantum Serre relations are satisfied in $\Sh(V)$
by the $[v_i]$ follows from \cite{Rosso:shuffle}, Lemma 14 (the proof
relies on $q$-binomial coefficients identities, which are proved by
induction); this proves the first part of the Lemma. 

Let us show that $U_\hbar\B_+ / \hbar U_\hbar\B_+$ is isomorphic
$U\B_+$. $U_\hbar\B_+ / \hbar U_\hbar\B_+$ is equal to the quotient of 
$\CC\langle h_i,D_j,x_i^+, i = 1,\ldots,n,  j = 1,\ldots,2n-r
\rangle[[\hbar]]$  by the sum of $\hbar \CC\langle
h_i,D_j,x_i^+\rangle[[\hbar]]$ and the closed ideal generated by the relations
$[h_i,e_j] = a_{ij}e_j$,  $[D_i,e_j] = \delta_{ij}e_j$ 
 and the quantum Serre relations (\ref{quantum:serre}).  This sum is the
same as that of $\hbar \CC\langle h_i,D_j,x_i^+\rangle[[\hbar]]$ and the
closed ideal generated by $[h_i,e_j] = a_{ij}e_j$,  $[D_i,e_j] =
\delta_{ij}e_j$  and the classical Serre relations. The quotient of 
$\CC\langle h_i,D_j,x_i^+ \rangle[[\hbar]]$ by this last space  is equal
to $U\B_+$.  This proves the second part  of the Lemma. 

The third part is immediate. 
\hfill \qed \medskip 

\begin{prop} \label{ratiu}
Let $\LL$ be a Lie algebra and let $J$ be a two-sided  ideal of $U\L$
such that $\Delta_{U\LL}(J)\subset U\LL \otimes J  +  J \otimes U\LL$.
Then $\J = J\cap \LL$ is an ideal of the Lie algebra $\LL$ and we have
$J = (U\LL)\J = \J(U\LL)$. 
\end{prop}

{\em Proof.}
We first show: 

\begin{lemma} \label{lausanne}
Let $\LL$ be a Lie algebra and let $J$ be a left ideal of $U\LL$
such that $\Delta_{U\LL}(J) \subset J \otimes U\LL + U\LL \otimes J$. 
Let $\J$ be the intersection $\LL\cap J$. Then $J$ is equal to 
($U\LL)\J$. 
\end{lemma}

{\em Proof of Lemma.} Denote by $(U\LL)_n$ the subspace of $U\LL$ spanned
by the monomials in elements of $\LL$ of degree $\leq n$. Let us set
$\bar\Delta_{U\LL}(x) = \Delta_{U\LL}(x) - x\otimes 1 - 1\otimes x$.  We
have $\bar\Delta_{U\LL}((U\L)_n) \subset \sum_{p,q>0, p+ q = n} (U\LL)_p
\otimes (U\LL)_q$. Denote by $J_n$ the intersection   $J\cap (U\LL)_n$.
Then we have $\bar\Delta_{U\LL}(J_n) \subset  \sum_{p,q>0, p+q = n}
J_p\otimes (U\LL)_q + (U\LL)_p \otimes J_q$.   

Let us show by induction that $J_n$ is contained in $(U\LL)_{n-1}\J$. 
This is clear if $n= 1$; assume it is true at order $n-1$. 
Let $x$ be an element of $J_n$. Then 
$\bar\Delta_{U\LL}(x)$ is contained in $\sum_{p,q>0, p+q = n}
(U\LL)_{p-1}\J \otimes (U\LL)_q + (U\LL)_p \otimes (U\LL)_{q-1}\J$.   

Let $\bar x$ be the image of $x$ in $(U\LL)_n / (U\LL)_{n-1}$. 
$(U\LL)_n / (U\LL)_{n-1}$ is isomorphic to the $n$th symmetric 
power $S^n\LL$. Let $\Delta_{S\LL}$ be the coproduct of the 
symmetric algebra $S\LL$, for which elements of degree $1$ are 
primitive, and set $\bar\Delta_{S\LL}(\bar x) = \Delta_{S\LL}(\bar x)
- \bar x \otimes 1 - 1\otimes \bar x$. Then $\bar \Delta(\bar x)$
is contained in $\sum_{p,q>0,p+q = n} (S^{p-1}\LL) \J \otimes S^q \LL 
+ S^p \LL \otimes (S^{q-1}\LL) \J$. It follows that $\bar x$ belongs to 
$(S^{n-1}\LL) \J$.  The difference of $x$ with some element of 
$(U\LL)_{n-1}\J$ therefore belongs to $(U\LL)_{n-1}$, so that it belongs to 
$J_{n-1}$ and by hypothesis to $(U\LL)_{n-2}\J$. Therefore, $x$ belongs to 
$(U\LL)_{n-1}\J$. This proves the Lemma. 
\hfill \qed\medskip 

Let us prove Prop.\ \ref{ratiu}. Lemma \ref{lausanne} implies that $J =
(U\LL) \J$ and its analogue for right Hopf ideals implies that $J = \J
(U\LL)$. Therefore,  we have $(U\LL)\J = \J (U\LL)$. Let us fix $x$ in
$\LL$ and $j$ in $\J$, then  $[x,j]$ belongs $J$; since it also belongs
to $\LL$, $[x,j]$ belongs to $\J$.  Therefore $\J$ is an ideal of $\LL$.
 \hfill \qed\medskip 

Let $J$ be the kernel of the cocommutative Hopf algebras morphism $p$
defined in Lemma \ref{pourim}, 3). Let us set $\J = J\cap \B_+$ and 
$\A = \B_+ / \J$. prop.\ \ref{ratiu} then implies 

\begin{lemma} The Lie algebra structure on $\B_+$ induces a Lie algebra
structure on $\A$. Moreover, $\cV_0$ is isomorphic with $U\A$, and $p$
can be identified with the quotient map $U\B_+ \to  U\A$.  
\end{lemma}

Define $\delta_{\cV_0}$ as 
${{\Delta_\cV -  \Delta_\cV'}\over\hbar}$ mod $\hbar$.  $\delta_{\cV_0}$
is a linear map from $\cV_0$ to the antisymmetric part of  its tensor
square $\wedge^2\cV_0$. It obeys the rules 
\begin{equation} \label{co-leib} 
 (\Delta_{\cV_0}\otimes id)\circ \delta_{\cV_0} = 
 (\delta_{\cV_0}^{2\to 23}  + \delta_{\cV_0}^{2\to 13}) \circ
 \Delta_{\cV_0},   
\end{equation} 
$$ 
\Alt ( \delta_{\cV_0}\otimes
id)\circ\delta_{\cV_0} ) = 0,  
$$ 
$$ \delta_{\cV_0}(xy)  = \delta_{\cV_0}(x) \Delta_{\cV_0}(y)   + \Delta_{\cV_0}(x)
\delta_{\cV_0}(y) \on{\ for\ }  x,y \on{\ in\ } \cV_0,  
$$ if
$\delta_{\cV_0}(y) = \sum_i y'_i \otimes y''_i$, we set  $\delta^{2\to
23}_{\cV_0}(x\otimes y) = \sum_i x\otimes y'_i \otimes y''_i$,  and
$\delta^{2\to 13}_{\cV_0}(x\otimes y) = \sum_i y'_i \otimes x\otimes
y''_i$.  These rules are the co-Leibnitz, co-Jacobi and Hopf
compatibility  conditions; they mean that $(\cV_0,\delta_{\cV_0})$ is a
Hopf co-Poisson algebra (see \cite{Drinf:ICM}).  

\begin{lemma} \label{A:bialg} 
$\delta_{\cV_0}$ maps $\A$ to $\wedge^2 \A$. 
\end{lemma}

{\em Proof.} Let $a$ be an element of $\A$ and set $\delta_{\cV_0}(a) =  \sum_i
x_i \otimes y_i$, where $(y_i)$ is a free family. Then (\ref{co-leib})
implies that $\Delta_{\cV_0}(x_i) \otimes y_i = \sum_i x_i \otimes
1\otimes y_i  + 1 \otimes x_i\otimes y_i$, so that each $x_i$ is
primitive and therefore  belongs to  $\A$.  So $\delta_{\cV_0}(a)$
belongs to $\A\otimes \cV_0$. Since  $\delta_{\cV_0}(a)$ is also antisymmetric, 
it belongs to $\wedge^2\A$.   
\hfill \qed\medskip 

Call $\delta_\A$ the map from $\A$ to $\wedge^2\A$ defined as the
restriction  of $\delta_{\cV_0}$ to $\A$. $(\A,\delta_\A)$ is then a Lie
bialgebra, which means that   $\delta_\A$ is a $1$-cocyle of $\A$ with
values in the antisymmetric part of the  tensor square of its
adjoint representation,  satisfying the co-Jacobi  identity
$\Alt(\delta_\A\otimes id)\delta_\A = 0$. 
 
\begin{remark} 
The Hopf co-Poisson algebra and Lie bialgebra axioms were introduced  by
Drinfeld in \cite{Drinf:ICM}.  Drinfeld showed that the quantization of
a cocommutative Hopf algebra lead to such structures. He also  stated
that there is an equivalence of categories between the category of  Hopf
co-Poisson algebras and that of Lie bialgebras. Lemma \ref{A:bialg} can
therefore be viewed as the proof of one part of this statement (from
Hopf co-Poisson to Lie bialgebras).  It is also not difficult to prove the 
other part (from Lie bialgebras to Hopf co-Poisson).  
\end{remark}

\subsubsection{Kac-Moody Lie algebras} \label{KM:rappels}

Let $\G$ be the Kac-Moody Lie algebra associated with $A$.
$\G$ has generators $\bar x_i^\pm,\bar h_i$, $i = 1,\ldots, n$ and 
$\bar D_j, j  = 1, \ldots, n-r$,  and relations
\begin{equation} \label{comm:cartan:class}
[h,h'] = 0 \on{\ if\ } h,h'\in \{\bar h_i, \bar D_j\}, 
\end{equation}
\begin{equation} \label{cartan:n+:class}
[\bar h_i, \bar x_{i'}^\pm] = \pm a_{ii'} \bar x_{i'}^\pm,  
\quad [\bar D_j, \bar x_{i}^\pm] = \pm \delta_{ij} \bar x_{i}^\pm, 
\end{equation}
\begin{equation} \label{serre:class}
\ad(\bar x_i^\pm)^{1 -a_{ij}} (\bar x_j^\pm) = 0,  
\end{equation}
$$
[\bar x_i^+, \bar x_{i'}^-] = \delta_{ii'} \bar h_i, \on{\ for \ all\ }
i,i' = 1, \ldots,n, j = 1,\ldots,n-r. 
$$

Let $\HH$ and $\N_\pm$ be the subalgebras of $\G$ generated by  $\{\bar
h_i,\bar D_j\}$ and $\{\bar x_i^\pm\}$.  Then $\G$ has the Cartan decomposition
$\G = \N_+ \oplus \HH \oplus \N_-$.  Let us set $\B_\pm = 
\HH\oplus\N_\pm$. $\N_\pm$ and $\B_\pm$  are the Lie algebras with 
generators $\{\bar x_i^\pm\}$ and $\{\bar h_i, \bar D_j,\bar x_i^\pm\}$  and
relations (\ref{serre:class}) for $\N_\pm$ and 
(\ref{comm:cartan:class}), (\ref{cartan:n+:class}) and 
(\ref{serre:class}) for $\B_\pm$.  $\G$ is endowed with a nondegenerate
bilinear  form $\langle , \rangle_\G$, which is determined by 
$\langle \bar h_i, \bar h_{i'} \rangle_\G = d_{i'}^{-1} a_{ii'}$,
$\langle \bar x_i^+,  \bar x_{i'}^- \rangle_\G = d_i^{-1}\delta_{ii'}$, 
$\langle \bar h_i, \bar D_j \rangle_\G = d_i^{-1}\delta_{ij} $, 
and that its values for all other pairs of generators 
is zero (see \cite{DGK,Kac}).  

\subsubsection{Hopf algebras $U_\hbar\B_\pm$}

Define $U_\hbar\B_\pm$ as the algebras with generators $h_i^\pm,D_j^\pm$ and 
$x_i^\pm$, and relations
$$
[h_i^\pm,x^\pm_{i'}] = \pm a_{ii'} x^\pm_{i'},   \quad
[D_j^\pm,x^\pm_{i}] = \pm \delta_{ij} x^\pm_{i},  
$$
and relations (\ref{quantum:serre}), with $e_i$ replaced by $x^\pm_i$. 
It is easy to see that the maps $e_i \mapsto x^+_i$ and $f_i\mapsto x^-_i$ define 
algebra inclusions of $U_\hbar\N_\pm$ in $U_\hbar\B_\pm$.   
We have  Hopf algebra structures on  $U_\hbar\B_+$ and
$U_\hbar\B_-$,  defined by 
$$
\Delta_\pm(h^\pm) = h^\pm \otimes 1 + 1\otimes h^\pm \on{\ for\ }
h\in \{\bar h_i,\bar D_j\}, 
\quad  
\Delta_\pm(x_i^\pm) = x_i^\pm \otimes e^{\pm\hbar d_i h_i^\pm} 
+ 1\otimes x_i^\pm. 
$$

\subsubsection{Comparison lemmas}

Recall that $\cV$ is the direct sum of its graded components. 
Its component of degree zero is $\CC[h_i^{\cV}, D_j^{\cV}][[\hbar]]$. 
Let $h_i^{\A}, D_j^{\A}$ be the images of  $h_i^{\cV}, D_j^{\cV}$ 
in $\cV_0$. We have an inclusion of $\CC[h_i^{\A},D_j^{\A}]$
in $\cV_0$. It follows that the  $h_i^{\A}$ and $D_j^{\A}$ are 
linearly independent and commute to each other. 

On the other hand, as the elements $h_i^{\cV}$ and $D_j^{\cV}$ are 
primitive in $\cV$, the $h_i^{\A}$ and $D_j^{\A}$ are also primitive; it
follows that they belong to $\A$, and we have 
\begin{equation} \label{shoshan}
p(\bar h_i) = h_i^\A,  p(\bar D_j) = D_j^\A,  \quad \delta_\A(h_i^\A)  =
\delta_\A(D_j^\A)  = 0. 
\end{equation}

The degree $\eps_i$ component of $\cV$ is  
$\CC[h_i,D_j][[\hbar]] \cdot[v_i]$, and the map $x\mapsto x[v_i]$ is a 
$\CC[[\hbar]]$-module isomorphism from $\CC[h_i,D_j][[\hbar]]$ to this 
component. Therefore, $[v_i]$ has a nonzero
image in $\cV / \hbar\cV = \cV_0$. Since we have 
$$
\Delta_\cV([v_i]) = [v_i]\otimes e^{\hbar d_i h_i^{\cV}} + 
1\otimes [v_i], 
$$
$[v_i]\ mod\ \hbar$ is primitive in $\cV_0$. Therefore, $[v_i]\ 
mod\ \hbar$ belongs to $\A$; call $v_i^\A$ this element of $\A$. 
It is clear that 
\begin{equation} \label{habira}
\delta_\A(v_i^\A) = d_i v_i^\A \wedge h_i^{\A}.
\end{equation}

On the other hand, $p_\hbar(\bar x_i^+) = [v_i]$ implies that 
\begin{equation} \label{megila}
p(\bar x_i^+)=v_i^\A. 
\end{equation}

Recall that we have a Lie bialgebra structure on $\B_+$; it consists in
a map $\delta_{\B_+}$ from $\B_+$ to $\wedge^2\B_+$, which is uniquely 
determined by the conditions $\delta_{\B_+}(\bar x_i^+) = d_i \bar x_i^+
\wedge \bar h_i$,  $\delta_{\B_+}(\bar h_i)  = \delta_{\B_+}(\bar D_j) 
= 0$ and that it satisfies the $1$-cocycle identity. 

\begin{lemma} \label{amo}
$p$ is also a Lie bialgebra morphism from $(\B_+,\delta_{\B_+})$ to 
$(\A,\delta_\A)$. 
\end{lemma}

{\em Proof.} This means that 
\begin{equation} \label{esther}
(\wedge^2 p) \circ \delta_{\B_+}(x) =  \delta_{\A}
\circ p(x), \on{\ for\ }x\on{\ in \ }\B_+.
\end{equation} 

For $x$ equal to $\bar h_i$ and $\bar D_j$, (\ref{esther}) follows from 
(\ref{shoshan}). It follows from (\ref{habira}) and (\ref{megila}) that 
$$
\delta_\A\circ p(\bar x_i^+) = \delta_\A(v_i^\A) = d_i v_i^\A \wedge h_i^\A
= (\wedge^2 p)(d_i \bar x_i^+\wedge \bar h_i) =\wedge^2 p(\delta_\A(\bar x_i^+)),     
$$
so that (\ref{esther}) is also true for $x = \bar x_i^+$. 

Since $p$ is a Lie algebra morphism, both sides of (\ref{esther})  are
$1$-cocycles of $\B_+$ with values in the antisymmetric part of the
tensor square of $(\A,\ad_\A\circ p)$. Therefore, (\ref{esther}) holds
on the subalgebra of $\B_+$ generated by the $\bar h_i,\bar D_j$ and
$\bar x_i^+$, which is $\B_+$ itself. 
\hfill \qed \medskip

Denote by $\HH_\A$ the subspace of $\A$ spanned by the $h_i^\A$ and
$D_j^\A$; it forms an abelian Lie  subalgebra of $\A$. 

Since $p_\hbar(h_i) = h_i^{\cV}$ and $p_\hbar(D_j) = D_j^{\cV}$, 
the restriction of $p$ to the Cartan subalgebra $\HH$ of $\B_+$ 
is a Lie algebra isomorphism from $\HH$ to $\HH_\A$. 

Define, for $\al$ in $\HH^*$, the root subspace $\A[\al]$ associated
with $\al$ by 
$$
\A[\al] =\{x\in \A |   [p(h),x] = \al(h) x, \on{\ for\  all\ } 
h \on{\ in\ } \HH\}, 
$$
and as usual 
$$
\B_+[\al] =\{x\in \B_+ |   [h,x] = \al(h) x, \on{\ for\  all\ } 
h \on{\ in\ } \HH\}.  
$$

\begin{lemma}
$\A$ is the direct sum of its root subspaces $\A[\al]$, where $\al$ belongs to 
the set $\Delta_+\cup\{0\}$ of roots of $\B_+$. Each $\A[\al]$
is finite dimensional. $\delta_{\A}$ is a graded map from $\A$ to $\wedge^2\A$, 
therefore the graded dual $\A^*$ of $\A$, defined as $\oplus_{\al}\A[\al]^*$, 
has a Lie bialgebra structure. 
\end{lemma}

{\em Proof.} For any $\al$ in $\HH^*$, $p$ maps $\B_+[\al]$ to $\A[\al]$. It
follows that $\A$ is the sum of the root subspaces $\A[\al]$, where $\al$ belongs to 
the set of roots of $\B_+$. That this
sum is direct is  proved as in the case of $\B_+$:  let $\al_i$ be the
root such that $\bar x_i^+$ belongs to $\B_+[\al_i]$. Then $\al_1,\ldots,\al_n$
form a basis of $\HH^*$ (they are the simple roots of $\G$).  Let $(\bar
H_1,\ldots,\bar H_n)$ be the basis of $\HH$ dual to 
$(\al_1,\ldots,\al_n)$.  Then for any family $x_{\al}$ of $\A[\al]$ such
that $\sum_\al x_\al = 0$, we have, by applying $\ad(p(\bar H_1))^k$ to
this equality, $\sum_{(n_i)\in \NN^n} n_1^k x_{\sum_i n_i \al_i} = 0$,
which gives for any integer $a_1$,  $\sum_{(n_i)\in \NN^n | n_1 = a_1}
x_{\sum_i n_i\al_i} = 0$; applying  $\ad(p(\bar H_2))^k$,  we get
$\sum_{(n_i)\in \NN^n | n_1 = a_1, n_2 = a_2 } x_{\sum_i n_i\al_i} = 0$;
finally  each $x_\al$ vanishes. 

That $\delta_\A$ is graded follows from the fact that its restriction to 
$\HH_\A$ vanishes and from the cocycle identity. 
\hfill \qed\medskip 

\begin{lemma} \label{isom:Lie}
 $p$ is a Lie bialgebra isomorphism.  
\end{lemma}

{\em Proof.} It follows from Lemma \ref{amo} that $p^*$ is an injective
Lie bialgebra  morphism from $\A^*$ in the graded dual $\B_+^*$ of $\B_+$. 
Recall that $\B_+^*$ is
isomorphic, as a Lie algebra, to $\B_- = \HH \oplus \N_-$;  this relies
on the nondegeneracy  of the invariant pairing between $\B_+$ and
$\B_-$, itself a consequence of  \cite{DGK} (in what follows, we will
denote by $\HH_-$ the Cartan subalgebra $\HH$ of $\B_-$).  We will show that
the image of $p^*$ contains a generating family of $\B_-$. 

Let us denote by $\HH_{\A^*}$ the space of forms on $\A$, which  vanish
on all the $\A[\al],\al\neq 0$. The duality beween  $\B_+$ and $\B_-$
identifies $\HH_-$ with the space of the forms on  $\B_+$ which vanish on
$\N_+ = \oplus_{\al\neq 0}\B_+[\al]$.  We have $p(\B_+[\al])\subset
\A[\al]$ for any $\al$, therefore  
$$
p^*[(\oplus_{\al\neq 0}\A[\al])^{\perp}] \subset 
(\oplus_{\al\neq 0}\B_+[\al])^{\perp}, 
$$
which means that $p^*(\HH_{\A^*})\subset \HH_-$. Since $p^*$ is
injective and $\HH_{\A^*}$ and $\HH_-$ have the same dimension, $p^*$
induces an isomorphism  between $\HH_{\A^*}$ and $\HH_-$. It follows
that the image of  $p^*$ contains $\HH_-$. 

Since $\bar x_i^+$ belongs to $\B_+[\al_i]$,  the element $v_i^\A$ of $\A$
defined before Lemma \ref{amo} belongs to  $\A[\al_i]$. We have seen that
$v_i^\A$ is nonzero. 

Let $\xi_i$ be the element of $\A^*$ which is $1$ on $v_i^\A$ and zero
on each $\A[\al]$, $\al\neq\al_i$.  For $x$ in $\B_+[\al]$,
$\al\neq\al_i$, $\langle p^*(\xi_i), x \rangle_{\B_+\times\B_-} =
\langle \xi_i, p(x) \rangle_{\A^*\times\A} = 0$ because $\xi_i$ vanishes
on $\A[\al]$.  It follows that $p^*(\xi_i)$ has weight  $-\al_i$ in
$\B_-$. On the other hand,  $p^*(\xi_i)$ is nonzero, because $p^*$ is
injective, so it is a nonzero constant times $\bar x_i^-$. 

Since the image of $p^*$ contains $\HH_-$ and the
$\bar x_i^-$, $p^*$ is an isomorphism.  \hfill \qed \medskip

\begin{lemma} \label{isom:n+}
$p_\hbar$ mod $\hbar$ restricts to an isomorphism of $\NN^n$-graded 
algebras from $U\N_+$ to  $\langle V \rangle / \hbar \langle V \rangle$.
\end{lemma}

{\em Proof.} It follows from Lemma \ref{isom:Lie} that $p_\hbar$ mod
$\hbar$ induces an isomorphism from $U\B_+ = U_\hbar\B_+ / \hbar
U_\hbar\B_+$ to  $\cV_0 = \cV / \hbar\cV$. Therefore it induces an
isomorphism from $U\N_+$ to its image in $\cV_0$. Since $U\N_+$
coincides with the image of $U_\hbar\N_+$ by the  projection 
$U_\hbar\B_+ \to U_\hbar\B_+ / \hbar U_\hbar\B_+  = U\B_+$, this  
image coincides with that of the composed map 
\begin{equation} \label{composed}
U_\hbar\N_+\to U_\hbar\B_+\to \cV\to \cV_0.
\end{equation}  The  image
of the composed map $U_\hbar\N_+\to U_\hbar\B_+\to \cV$ is equal to 
$\langle V \rangle$. We have a $\CC[[\hbar]]$-module isomorphism  of
$\cV$ with $\langle V \rangle \otimes_{\CC[[\hbar]]}
\CC[h_i^{\cV},D_j^{\cV}][[\hbar]]$, so that $\hbar\cV\cap \langle V
\rangle = \hbar\langle V \rangle$. It follows that  the image of $\langle
V \rangle$ by $\cV\to \cV_0$ is $\langle V \rangle / \hbar\langle V
\rangle$. Therefore the image of (\ref{composed}) is 
$\langle V \rangle / \hbar\langle V \rangle$. 
\hfill \qed \medskip

{\em Proof of Thm.\ \ref{thm:first}.} Assign degree $\eps_i$ to the
generator  $e_i$ of $U_\hbar\N_+$. Then $U_\hbar\N_+$ is the direct sum
of its  homogeneous components $(U_\hbar\N_+)[\al]$, $\al\in\NN^n$,
which are finitely generated  $\CC[[\hbar]]$-modules. As a
$\CC[[\hbar]]$-module, $(U_\hbar\N_+)[\al]$ is therefore isomorphic  to
the direct sum $\oplus_{i=1}^{q_\al} \CC[[\hbar]] /
(\hbar^{n_i^{(\al)}}) \oplus \CC[[\hbar]]^{p_\al} $ of its torsion part 
with a free module (see Lemma \ref{str:modules}). 

Lemma \ref{pourim}, 2) implies that 
$U_\hbar\N_+[\al] / \hbar U_\hbar\N_+[\al] = U\N_+[\al]$ so 
\begin{equation} \label{pal}
p_\al+q_\al = \dimm U\N_+[\al].
\end{equation} 

On the other hand, $\langle V \rangle[\al]$ is a free finite dimensional
module over $\CC[[\hbar]]$, by Lemma \ref{<V>:free} above, so it is
isomorphic, as a $\CC[[\hbar]]$-module, to $\CC[[\hbar]]^{p'_\al}$.
$p_\hbar$ restricts to a surjective  $\CC[[\hbar]]$-module morphism from
 $U_\hbar\N_+[\al]$ to $\langle V \rangle[\al]$, therefore $p_\hbar$ maps 
the torsion part of $U_\hbar\N_+[\al]$ to zero and 
\begin{equation} \label{ma}
p_\al\geq p'_\al.
\end{equation} 
Moreover, 
\begin{equation} \label{kh}
p'_\al = \dimm U\N_+[\al]
\end{equation} 
by Lemma \ref{isom:n+}. 
It follows from (\ref{pal}), (\ref{ma}) and (\ref{kh}) that 
$p_\al = p'_\al$ and $q_\al=0$. 

This means that $U_\hbar\N_+$ has no torsion, and is isomorphic to 
$\langle V\rangle$. In view of Lemma \ref{pourim}, 2), this proves 
Thm.\ \ref{thm:first}. This also proves Cor.\ \ref{cor:comparison}. 
 \hfill \qed \medskip

\subsection{Nondegeneracy of Hopf pairing (proof of 
Thm.\ \ref{thm:second})}

Let $(v_i^*)$ be the basis of $V^*$ such that $\langle v_i^*,v_j\rangle 
= d_i^{-1}\delta_{ij}$.  Assign to $v_i^*$ the degree $-\eps_i$.  Let
$T(V^*)$ be the tensor algebra $\oplus_i (V^*)^{\otimes i}[[\hbar]]$. 
Define the braided tensor product strucutre on the tensor square of 
$T(V^*)$ according to (\ref{braided:tensor:pdt}).   $T(V^*)$ is endowed
with the braided Hopf structure defined by $\Delta_{T(V^*)}(v_i^*)  =
v_i^*\otimes 1 + 1\otimes  v_i^*$, for any $i = 1,\ldots,n$. We have a
surjective braided Hopf algebra morphism from $T(V^*)$ to 
$U_\hbar\N_-$, defined by $v_i^*\mapsto f_i$, for $i = 1,\ldots,n$.

Then we have a braided Hopf pairing 
$$
\langle , \rangle_{\Sh(V) \times T(V^*)}
:  \Sh(V) \times T(V^*)  
\to \CC((\hbar)), 
$$
defined by the rules
\begin{equation} \label{formula:crossed:pdt}
\langle [v_{i_1} | \cdots | v_{i_k}],\xi_{i'_1}\cdots 
\xi_{i'_{k'}}
 \rangle_{\cS \times U_\hbar\wt \B_- } =
{1\over\hbar} \delta_{kk'} \prod_{j=1}^k\langle v_{i_j}, 
\xi_{i'_j}\rangle_{V \times V^*} . 
\end{equation}

The ideal of $T(V^*)$ generated by the quantum Serre
relations is in the radical of this pairing (see e.g.\ \cite{Lusztig},
chap.\ 1; this is a consequence of $q$-binomial identities). 
 
It follows that 
$\langle , \rangle_{\Sh(V) \times T(V^*)}$ induces a braided Hopf 
pairing
$$
\langle , \rangle_{\Sh(V) \times U_\hbar\N_- }:  
\Sh(V) \times U_\hbar\N_- \to \CC((\hbar)). 
$$

By Thm.\ \ref{thm:first}, $U_\hbar\N_+$ is a braided Hopf subalgebra of 
$\Sh(V)$. The restriction of
$\langle , \rangle_{\Sh(V) \times U_\hbar\N_- }$ 
to $U_\hbar\N_+\times
U_\hbar\N_-$ therefore induces a braided Hopf pairing  between $U_\hbar\N_+$
and $U_\hbar\N_-$; since it coincides on generators with $\langle
, \rangle_{U_\hbar\N_+\times U_\hbar\N_-}$,  it is equal to 
$\langle , \rangle_{U_\hbar\N_+\times U_\hbar\N_-}$. 

View $V^{\otimes k}$ as a subspace of $\Sh(V)$. Assign degree $1$ to each 
element of $V^*$ in $T(V^*)$; then 
$T(V^*)$ is a graded algebra; we denote by 
$T(V^*)^{(k)}
$ is homogeneous component of degree $k$.  The
restriction of $\langle , \rangle_{\Sh(V) \times T(V^*) }$
to
$V^{\otimes k} \times T(V^*)^{(k)}$ can be identified with  the
natural pairing of $V^{\otimes k}$ with $(V^*)^{\otimes k}$, which is
nondegenerate. Therefore the annihilator of  $T(V^*)$ in
$\Sh(V)$ for $\langle , \rangle_{\Sh(V) \times T(V^*) }$ 
is zero. By Thm.\ \ref{thm:first}, it follows that the annihilator 
of $U_\hbar\N_-$ in $U_\hbar\N_+$ for $\langle , \rangle_{U_\hbar\N_\pm}$
is zero.

Since the pairing $\langle , \rangle_{U_\hbar\N_+\times U_\hbar\N_-}$ is
graded and the graded components of $U_\hbar\N_+$  and  $U_\hbar\N_-$
have the same dimensions (as $\CC[[\hbar]]$-modules), the pairing 
$\langle , \rangle_{U_\hbar\N_+\times U_\hbar\N_-}$ is nondegenerate.  
 \hfill\qed \medskip 

\subsection{The form of the $R$-matrix (proof of Prop.\ \ref{R:mat})}

Let us endow $U_\hbar\G = U_\hbar \B_+ \otimes U_\hbar \N_-$ with the double 
algebra structure such that $U_\hbar \B_+ \to U_\hbar\G$, $x_+ \mapsto 
x_+ \otimes 1$, and $U_\hbar \N_-\to U_\hbar\G$, $x_- \mapsto 1\otimes x_-$ 
are algebra morphisms and, if $e_i^{\G} = e_i\otimes 1$, 
$f_i^{\G} = 1\otimes f_i$, $h_i^\G = h_i \otimes 1$ and $D_j^\G = D_j \otimes 1$, 
$$
[e^\G_i ,  f^\G_j] = \delta_{ij}{{q^{d_ih^\G_i}  
- q^{ - d_ih_i^\G}}\over{q^{d_i} - q^{-d_i}}} ,  
$$
and 
$$
[h^\G_i,f^\G_j] = - a_{ij} f^\G_j, \quad [D^\G_j, f^\G_j] = -\delta_{ij} f^\G_j. 
$$

$U_\hbar\G$ is endowed with a topological Hopf algebra structure 
$\Delta: U_\hbar \G\mapsto U_\hbar \G \hat\otimes U_\hbar\G =
\limm_{\leftarrow N} ( U_\hbar \G \otimes U_\hbar \G )  / \hbar^N (
U_\hbar \G \otimes  U_\hbar \G )$,  extending $\Delta_+$ and $\Delta_-$
(\cite{Drinf:ICM}).  

Let $t_0$ be the element of $\HH \otimes \HH$ corresponding to the 
restriction of the invariant pairing of $\G$ to $\HH$ and let 
$\cR[\al]$ be the element of $[\limm_{\leftarrow N} ( U_\hbar \G
\otimes U_\hbar\G) /
\hbar^N (U_\hbar\G\otimes U_\hbar\G)][\hbar^{-1}]$ 
$$
\cR[\al] = \exp(\hbar t_0)P[\al]. 
$$
Then we have the equalities
\begin{equation} \label{kolia}
\cR[\al - \al_i] (e^\G_i \otimes q^{d_i h^\G_i }) +   \cR[\al] (1
\otimes e^\G_i) =  (q^{d_i h^\G_i } \otimes  e^\G_i )  \cR[\al]  +  
(e^\G_i \otimes 1)  \cR[\al - \al_i] , 
\end{equation}
for any $i = 1,\ldots, n$. 

\begin{lemma}
For any nonzero $\al$ in $\NN^n$, $P[\al]$ belongs to $\hbar U_\hbar \N_+ \otimes
U_\hbar\N_-$. 
\end{lemma}

{\em Proof.} Let us show this by induction on the height of $\al$ (we
say that the height of $\al = (\al_i)_{1\leq i \leq n}$ is $\sum_{i=
1}^n \al_i$). If $\al$ is a simple root $\al_i$,  $P[\al] = \hbar
e_i^{\G} \otimes f_i^{\G}$, so that the statement  holds when $\deg(\al)
= 1$. 

Assume that we know that $P[\al]$ belongs to $\hbar U_\hbar\N_+ \otimes
U_\hbar \N_-$ for any $\al$ of height $<\nu$. Let $\al$ 
be of height $\nu$. Let $v$ be the $\hbar$-adic valuation of 
$P[\al]$, and assume that $v\leq 0$. $\cR[\al]$ belongs to $\hbar^{v}
(U_\hbar\G \hat\otimes U_\hbar\G)$, and since $v\leq 0$  the equality
(\ref{kolia}) takes place in $\hbar^{v} (U_\hbar\G \hat\otimes
U_\hbar\G)$.  Let us set $R_\al = \hbar^{-v} \cR[\al]$ mod $\hbar$;
$R_\al$ is an element of $U\B_+ \otimes U\B_-$. Since $\hbar^{-v}
\cR[\al - \al_i]$ is zero mod $\hbar$, (\ref{kolia}) implies that 
$R_\al$ commutes with each $1\otimes \bar e_i$. 

Lemma 1.5 of \cite{Kac} says that if $a$ belongs to $\N_-$ and  commutes
with each $e_i$, then $a$ is zero. It follows that if $x$ belongs to 
$U\N_-$ and commutes with each $\bar e_i$, $x$ is scalar; and if in 
addition $x$ has nonzero degree, $x$ is zero. Therefore $R_\al$ is zero.

It follows that $v\geq 1$, which proves the induction. 
\hfill \qed\medskip 

It follows from \cite{Drinf:ICM} that the $\cR[\al]$ satisfy 
the quasi-triangular identities
\begin{equation} \label{QT1}
(\Delta\otimes id)\cR[\al] = \sum_{\beta,\gamma\in\NN^n, 
\beta + \gamma = \al} \cR[\beta]^{(13)} \cR[\gamma]^{(23)} , 
\end{equation}
\begin{equation} \label{QT2} 
(id \otimes \Delta)\cR[\al] = \sum_{\beta,\gamma\in\NN^n, 
\beta + \gamma = \al} \cR[\beta]^{(13)} \cR[\gamma]^{(12)}  . 
\end{equation}
Let us set, for $\al\neq 0$,  
$r[\al] = R[\al] / \hbar$ mod $\hbar$; $r[\al]$ belongs
to $U\B_+ \otimes U\B_-$. Dividing the equalities (\ref{QT1}) and 
(\ref{QT2}) by 
$\hbar$, we get $(\Delta_{U\B_+} \otimes id) r[\al] = 
r[\al]^{(13)} + r[\al]^{(23)}$ and   
$(id \otimes \Delta_{U\B_+}) r[\al] = 
r[\al]^{(12)} + r[\al]^{(13)}$. Therefore $r[\al]$ belongs to 
$\B_+ \otimes \B_-$. 

Moreover, (\ref{kolia}) implies the identity 
$$
\delta(x)_{(\beta,\al - \beta)} = [r[\beta - \al], x\otimes 1] + 
+ [r[\beta], 1\otimes x]
$$
for $x$ in $\G[\al]$, where we set 
$(\sum_i x_i\otimes y_i)_{(\al,\beta)} =  \sum_i  (x_i)_{(\al)} \otimes
(y_i)_{(\beta)}$ and we denote by $x_{(\al)}$ the degree $\al$ component
of an element $x$ of $U\G$. 
  
It follows that $r[\al]$ is the element of $\N_+[\al] \otimes 
\N_-[-\al]$ corresponding to the invariant pairing of $\G$. 

Let us now prove by induction on $k$ that it $\al$ belongs to 
$k\Delta_+ - (k-1)\Delta_+$, $P[\al]$ belongs to $\hbar^k U_\hbar 
\N_+ \otimes U_\hbar \N_-$ and 
$$
P[\al] = {{\hbar^k}\over{k!}}  
\sum_{\al_1,\ldots, \al_k \in \Delta_+,
\sum_{i = 1}^k \al_i = \al } r[\al_1] \cdots r[\al_k] + o(\hbar^k).  
$$
Assume that the statement is proved up to order $k-1$ and let $\la$
belong to $k\Delta_+ - (k-1)\Delta_+$. Then (\ref{QT1}) and the 
induction hypothesis implies that 
\begin{align} \label{hannibal}
& (\wt\Delta \otimes id) (P[\la]) = 
\sum_{\al_1, \ldots,\al_k\in\Delta_+, \sum_{i = 1}^k \al_i =  \la; i_j} 
\sum_{l,l' >0, l + l' = k} {{\hbar^k}\over{l! l'! }}
\\ & e_{\al_1,i_1} \cdots e_{\al_l,i_l} \otimes
e_{\al_{l+1},i_{l+1}} \cdots e_{\al_k,i_k} \otimes
f_{\al_1,i_1} \cdots f_{\al_k,i_k} , 
\end{align}
where $\wt\Delta(x) = \Delta (x) - x\otimes 1 - 1 \otimes x$. Let $\sigma$ be any 
permutation of $\{1,\ldots,k\}$. 

For any $\al_1, \ldots,\al_k$ in $\Delta_+$, such that $\sum_{i = 1}^k \al_i =  \la$, 
we have 
$$
f_{\al_1,i_1} \cdots f_{\al_k,i_k}  = f_{\al_{\sigma(1)},i_{\sigma(1)}} 
\cdots f_{\al_{\sigma(k)},i_{\sigma(k)}} + o(\hbar). 
$$ 
Indeed, the difference  of both sides is a sum of products of the
$[f_{\al_s,i_s}, f_{\al_t,i_t}]$ with elements of $\N_+$; but $\al_s +
\al_t$ does not belong to  $\Delta_+$ by hypothesis on $\la$, so
$[f_{\al_s,i_s}, f_{\al_t,i_t}] = o(\hbar)$. 

The right side of (\ref{hannibal}) can the be rewritten as
\begin{align*}
& \sum_{\al_1, \ldots,\al_k\in\Delta_+, \sum_{i = 1}^k \al_i =  \la; i_j} 
\sum_{l,l' >0, l + l' = k} 
{{\hbar^k}\over{l! l'! }}
{1\over{\on{card} \Sigma_{l,l'}}}
\\ & 
\sum_{\sigma\in\Sigma_{l,l'}}
e_{\al_{\sigma(1)},i_{\sigma(1)}} \cdots e_{\al_{\sigma(l)},i_{\sigma(l)}} \otimes
e_{\al_{\sigma(l+1)},i_{\sigma(l+1)}} \cdots e_{\al_{\sigma(k)},i_{\sigma(k)}} \otimes
f_{\al_1,i_1} \cdots f_{\al_k,i_k} , 
\end{align*}
where $\Sigma_{l,l'}$ is the set of shuffle transformations of $((1, \ldots,l), 
(l+1,\ldots, l+l'))$. 
Therefore the right side of (\ref{hannibal})  is equal to 
$$
{{\hbar^k}\over{k!}}\wt\Delta 
\left( \sum_{\al_1, \ldots,\al_k\in\Delta_+, \sum_{i = 1}^k \al_i =  \la; i_j} 
e_{\al_1;i_1} \cdots e_{\al_k;i_k} \otimes f_{\al_1;i_1} \cdots f_{\al_k;i_k}
\right) + o(\hbar^k). 
$$

Let $v$ be the $\hbar$-adic valuation of $P[\la]$. Assume that $v<k$. 
Set $\bar P[\al] = \hbar^{-v} P[\al]$ mod $\hbar$. Then if we call
$\Delta_0$ the coproduct of $U\N_+$, and we set $\wt \Delta_0(x) =
\Delta_0(x) - x \otimes 1 - 1 \otimes x$,  (\ref{hannibal}) gives $(\wt
\Delta_0 \otimes id) (\bar P[\al]) = 0$, so that  $\bar P[\al]$ belongs
to $\N_+ \otimes U \N_-$; since $\bar P[\al]$ also belongs to $U\N_+[\al]\otimes
U\N_-[-\al]$ and $\al$ does not belong to $\Delta_+$, $\bar P[\al]$ is zero, 
contradiction. Therefore $v\geq k$. Let us set $P'[\al] = \hbar^{-k}P[\al]$
mod $\hbar$; we find that 
$$
(\wt\Delta_0\otimes id) \left( P'[\al] - {1\over {k!}}
\sum_{\al_1, \ldots,\al_k\in\Delta_+, \sum_{i = 1}^k \al_i =  \la; i_j} 
\bar e_{\al_1;i_1} \cdots \bar e_{\al_k;i_k} \otimes 
\bar f_{\al_1;i_1} \cdots \bar f_{\al_k;i_k}
\right)  = 0, 
$$
so that $P'[\al]$ belongs to  ${1\over {k!}}
\sum_{\al_1, \ldots,\al_k\in\Delta_+, \sum_{i = 1}^k \al_i =  \la; i_j} 
\bar e_{\al_1;i_1} \cdots \bar e_{\al_k;i_k} \otimes 
\bar f_{\al_1;i_1} \cdots \bar f_{\al_k;i_k}
+ \N_+ \otimes U\N_-$; as $\N_+[\al]$ is zero, $P'[\al]$ is equal to 
${1\over {k!}}
\sum_{\al_1, \ldots,\al_k\in\Delta_+, \sum_{i = 1}^k \al_i =  \la; i_j} 
\bar e_{\al_1;i_1} \cdots \bar e_{\al_k;i_k} \otimes 
\bar f_{\al_1;i_1} \cdots \bar f_{\al_k;i_k}$, 
which proves the induction. 
\hfill \qed\medskip

\subsection{The generic case (proof of Cors.\ \ref{generic:1} and 
\ref{generic:2})}

We have the equality  $U_\hbar \N_+ \otimes_{\CC[[\hbar]]} \CC((\hbar))
= U_{q'}\N_+ \otimes_{\CC(q')} \CC((\hbar))$, therefore  the graded
components of $U_{q'}\N_+$ have the same dimension as  those of $U_\hbar
\N_+$. Cors.\ \ref{generic:1} and  \ref{generic:2}) follow.

\section{Quantum current algebras of finite type 
(proofs of results of sect.\ \ref{QC})} \label{proof:QC}

\subsection{PBW theorem and comparison with Feigin-Odesskii 
algebra (proofs of Thm.\ \ref{thm:third} and Cor.\ \ref{cor:second})}

\subsubsection{Identification of algebras generated by the classical limits 
of quantum currents relations}

Recall that $A$ is now assumed of finite type. Define $L\B_+$ as the Lie
subalgebra  $(\HH\otimes\CC[t^{-1}]) \oplus (\N_+\otimes\CC[t,t^{-1}])$
of $\G\otimes\CC[t,t^{-1}]$. 

\begin{prop} \label{fargo}
Define $U_\hbar L\B_+$  and $\wt U_\hbar L\B_+$ 
as the algebra with generators $h_i[k], 
i = 1,\ldots, n, k\leq 0$ and $x^+_i[k], i = 1,\ldots, n, k\in\ZZ$, 
and relations 
$$
[h_i[k], h_j[l]] = 0, \quad [h_i[k], x^+_j[l]] = {{ q^{k d_i a_{ij}}  
- q^{ - k d_i a_{ij}} }\over{2\hbar k d_i}} x^+_j[k+l], 
$$
and relations (\ref{crossed:vertex}) and (\ref{q:serre:nr})  
among the $x^+_j[k]$ (with $e_i$ replaced by $x_i^+$), 
resp.\ (\ref{crossed:vertex}) and (\ref{variant:serre}). 
There are algebra
isomorphisms  from $U_\hbar L\B_+ / \hbar U_\hbar L\B_+$  and 
$\wt U_\hbar L\B_+ / \hbar \wt U_\hbar L\B_+$  to $U L\B_+$, 
sending $h_i[k]$ to $\bar h_i\otimes t^k$ and $x^+_i[k]$ to 
$\bar x_i^+\otimes t^k$. 
\end{prop}

{\em Proof.} $U_\hbar L\B_+ / \hbar U_\hbar L\B_+$ is the algebra 
with generators $\bar h_i[k],\bar e_i[l], 
1\leq i,j\leq n, k\leq 0,l\in\ZZ$  and relations 
$$
[\bar h_i[k], \bar e_i[l]] = a_{ij} \bar e_i[k+l],  
$$
and
\begin{equation} \label{ceb:1}
(z-w)[\bar e_i(z),\bar e_j(w)] = 0, 
\end{equation}
$$ 
\Sym_{z_1,\ldots,z_{1 - a_{ij}}}
\left( \ad(\bar e_i(z_1))  \cdots \ad(\bar e_i(z_{1 - a_{ij}}))
(\bar e_j(w)) \right)  = 0 , 
$$ 
where $\bar e_i(z) = \sum_{k\in \ZZ} \bar e_i[k] z^{-k}$. 
It follows from (\ref{ceb:1}) with $i = j$ that we have $[\bar
e_i[n],\bar e_j[m]] = 0$ for all $n,m$. Therefore,   $\ad(\bar e_i(z_1))
 \cdots \ad(\bar e_i(z_{1 - a_{ij}})) (\bar e_j(w))$ is symmetric in the
$z_i$, so that the last equation is equivalent to 
\begin{equation} \label{ceb:2}
\ad(\bar e_i(z_1))  \cdots \ad(\bar e_i(z_{1 - a_{ij}}))
(\bar e_j(w)) = 0. 
\end{equation}

On the other hand, $\wt U_\hbar L\B_+ / \hbar \wt U_\hbar L\B_+$ is the algebra 
with generators $\bar h_i[k]',\bar e_i[l]', 
1\leq i,j\leq n, k\leq 0,l\in\ZZ$  and relations 
$$
[\bar h_i[k]', \bar e_i[l]'] = a_{ij} \bar e_i[k+l]',  
$$
and
\begin{equation} \label{cef:1}
(z-w)[\bar e_i(z)',\bar e_j(w)'] = 0, 
\end{equation}
and
\begin{equation} \label{cef:2}
(\ad e_i[0]')^{1 - a_{ij}} e_j[k]' = 0. 
\end{equation}
The algebras presentad by the pairs of relations (\ref{ceb:1}) and (\ref{ceb:2}) on one hand, 
and (\ref{cef:1}) and (\ref{cef:2}) on the other, are isomorphic.  
Indeed, (\ref{ceb:1}) and (\ref{cef:1}) are equivalent, 
and (\ref{ceb:2}) implies (\ref{cef:2}); on the other hand,
(\ref{ceb:1}) implies that 
$[e_i[0]', e_j[k+l]'] = [e_i[k]', e_j[l]']$, 
so that  
$[e_i[0]',[e_i[0]', e_j[k+k'+l]']] = [e_i[0]',[e_i[k]', e_j[k'+l]']]
= [e_i[k]',[e_i[0]', e_j[k'+l]']]$, because the (\ref{ceb:1}) 
implies that $[e_i[0]',e_i[k]'] = 0$, therefore  
$[e_i[0]',[e_i[0]', e_j[k+k'+l]']] 
= [e_i[k]',[e_i[k']', e_j[l]']]$; one then proves by 
induction that 
$(\ad e_i[0]')^p ( e_j[k + k_1 + \cdots + k_p]) 
= \ad e_i[k_1]' \cdots \ad e_i[k_p]'( e_j[k])$. 
With $p = 1 - a_{ij}$, this relation shows that the $e_i[k]'$ 
satisfy (\ref{ceb:2}).

If follows that if $\wt F_+$ is the Lie algebra 
defined in Prop.\ \ref{oz}, both quotient algebras $U_\hbar L\B_+ / \hbar U_\hbar L\B_+$  and 
$\wt U_\hbar L\B_+ / \hbar \wt U_\hbar L\B_+$  are isomorphic to 
the crossed product of $U\wt F_+$ with the derivations 
$\tilde h_i[k]'$, defined by $\wt h_i[k]'(\bar e_i[l]) = a_{ij}
\bar e_i[k+l]$. 

It is clear that there is a unique Lie algebra morphism $j_+$ from the
Lie algebra  $\wt F_+$ defined in Prop.\ \ref{oz} to  $\N_+\otimes
\CC[t,t^{-1}]$, sending $\bar e_i[k]$ to $\bar x_i^+ \otimes t^k$. Let
us prove that it is an isomorphism. 

For this, let us define $\wt F$ as the Lie algebra with generators
$\wt x_i^\pm[k],\wt h_i^\pm[k]$, $1\leq i\leq n, k\in\ZZ$,
and relations given by the coefficients of the monomials in 
$$
(z-w)[\wt x_i^\pm(z), \wt x_j^\pm(w)] = 0,  \quad\on{if}\quad 
x,y\in\{\wt x_i^\pm\}, 
$$
$$
[\wt h_i(z), \wt h_j(w) ] = 0, 
$$
$$
[\wt h_i(z), \wt x_j^\pm(w)] = \pm a_{ij} \delta(z/w) \wt x_j^\pm(w), 
$$
$$
[\wt x_i^+(z), \wt x_j^-(w) ] = \delta_{ij}\delta(z/w)
\wt h_i(z),  
$$
$$
\ad(\wt x_i^\pm(z_1))  \cdots \ad(\wt x_i^\pm ( z_{1- a_{ij}}) ) 
(\wt x_j^\pm(w)) = 0,  
$$
where we set $\wt x(z) = \sum_{k\in\ZZ} x[k]z^{-k}$ for $x$ in 
$\{\wt x_i^\pm,\wt h_i\}$.  

\begin{lemma} \label{cite:U}
(In this Lemma, $\G$ may be an arbitrary Kac-Moody Lie algebra.) 
Let $W$ be the Weyl group of $\G$, and $s_i$ be its elementary 
reflection associated to the root $\al_i$. Then there is a unique 
action of $W$ on $\wt F$ such that 
$$
s_i(\wt x_i^\pm[k]) = \wt x_i^\mp[k], 
$$
$$
s_i(\wt x_j^\pm[k]) = \ad(\wt x_i^\pm[0])^{-a_{ij}} 
(\wt x_j^\mp[k]), \quad \on{if}\quad j\neq i, 
$$
$$
s_i(\wt h_j^\pm[k]) = \wt h_j^\pm[k] - a_{ij}\wt h_i^\pm[k]. 
$$
\end{lemma}

{\em Proof of Lemma.} The proof follows the usual proof for Kac-Moody Lie 
algebras. For example, if $j,k$ are different from $i$, we have 
\begin{align*}
(z-w) [s_i(\wt x_j^\pm(z) 
) , s_i( \wt x_k^\pm(w)  ) ] & = 
(z-w) [ \ad(\wt x_i^\pm[0])^{-a_{ij}}
(\wt x_j^\pm(z)) 
,  \ad(\wt x_i^\pm[0])^{-a_{ik}}(\wt x_k^\pm(w)) ]
\\ & = 
\ad (\wt x_i^\pm[0])^{-a_{ij} -a_{ik} } \left(
(z-w) [\wt x_j^\pm(z) ,  \wt x_k^\pm(w)  ]  \right )  
\end{align*}
because the Serre relations imply that $\ad(\wt x_i^\pm[0])^{1 
- a_{il}} (\wt x_l(u)) = 0$ for $j = k,l$ and $u = z,w$; therefore 
$(z-w) [s_i(\wt x_j^\pm(z) 
) , s_i( \wt x_k^\pm(w)  ) ]$ is zero. 
\hfill \qed\medskip 

\begin{lemma} \label{india} 
There is a unique Lie algebra isomorphism $j$ from $\wt F$ to 
$\G\otimes \CC[t,t^{-1}]$, such that $j(\wt x[k])  = \bar x\otimes t^k$, 
for any $x$ in $\{x_i^\pm, h_i\}$ and $k$ in $\ZZ$. 
\end{lemma}

{\em Proof of Lemma.} Let $\wt F_-$ be the Lie algebra with generators 
$\bar x_i^-[k], 1\leq i\leq n, k\in\ZZ$ and relations (\ref{ceb:1}) and (\ref{ceb:2}),
with $\bar x_i^+[k]$ replaced by $\bar x_i^-[k]$, and let $\wt H$ be the 
abelian Lie algebra with generators $\bar h_i[k], 1\leq i \leq n, k\in\ZZ$. 
There are unique Lie algebra morphisms from $\wt F_\pm$ and $\wt H$ to $\wt F$, 
sending the $\bar x_i^\pm[k]$ to $\wt x_i^\pm[k]$ and the $\bar h_i[k]$ to 
$\wt h_i[k]$. These morphisms are injections, so that we will indentify 
$\wt F_\pm$ and $\wt H$ with their images in $\wt F$. 

Moreover, let $F_\pm$ be the free Lie algebras with generators 
$x_i^\pm[k]^{F}$, $i = 1,\ldots,n$, $k$ integer. Endow $F_\pm \oplus \wt H$
with the Lie algebra structure such that $\wt H$ is abelian, $F_\pm$ is a 
Lie subalgebra of $F_\pm \oplus \wt H$, and  $[\wt h_i[k], x_i^\pm[l]^F]
= \pm a_{ij} x_i^\pm[k+l]^F$

There are unique derivations  $\Phi^{\mp}_{i,k}$ from $F_\pm$ to 
$F_\pm \oplus \wt H$ such that 
$$
\Phi^\mp_{i,k}(x_{i'}^\pm[l]) = \delta_{ii'} \wt h_i[k+l].  
$$
Let $I^F_\pm$ be the ideals of $F_\pm$ generated by relations (\ref{ceb:1})
and (\ref{ceb:2}); then computation shows that $I^F_\pm$ are 
preserved by the $\Phi^\mp_{i,k}$. It follows that $\wt F$ is the
direct sum of its subspaces $\wt F_\pm$ and $\wt H$. 

The rules $\deg(\wt x_i^\pm[k]) = (\pm\eps_i,k)$ and $\deg(\wt h_i[k]) =
(0,k)$ define a Lie algebra grading of $\wt F$ by $\ZZ^n\times \ZZ$,
because the relations of $\wt F$ are homogeneous for this grading.
Clearly,  $\dimm \wt F_\pm[(\pm\eps_i,k)] = 1$ for any $i$ and $k$, so
that  $\dimm \wt F[(\pm\eps_i,k)] = 1$. 

Let $\al$ be any root on $\G$. Then there is some simple root
$\pm\eps_i$ an element $w$  of $W$ such that $\al = w(\pm\eps_i)$. 
Then $\wt F[(\pm\eps_i,k)] = \wt F[(\al,k)]$ so that $\dimm\wt F[(\al,k)] = 1$. 

It is clear that the map $j$ defined in the statement of the Lemma
defines a Lie algebra morphism. Define a grading by $\ZZ^r\times\ZZ$ on 
$\G\otimes\CC[t,t^{-1}]$,  by the rules $\deg(x\otimes t^k) =
(\deg(x),k)$, for $x$ a homogeneous  (for the root grading) element of
$\G$.  Then $j$ is a graded map. Moreover, if $\al$ is in $\pm\Delta_+$
and $x$ is a nonzero element of $\G$ od degree $\al$, then $x$ can be
written  as a $\sum \la_{i_1,\cdots,i_p} [x^\pm_{i_1}, [\ldots,
x^\pm_{i_p}]]$;   then the image by $j$ of  $\sum \la_{i_1,\cdots,i_p}
[x^\pm_{i_1}[0], [\ldots, x^\pm_{i_p}[k]]]$ is equal to $x\otimes t^k$;
therefore the map induced by $j$ from $\wt F[(\al,k)]$ to
$\G\otimes\CC[t,t^{-1}][(\al,k)]$ is nonzero and therefore an
isomorphism.   

It follows that $\Ker j$ is equal to $\sum_{\al\in\ZZ^n \setminus
(\Delta_+\cup \{0\}\cup (-\Delta_+)), k\in\ZZ} \wt F[(\al,k)]$. Any element of 
$\wt F[(\al,k)]$ is a linear combination of brackets 
$[x^\pm_{i_l}[k_l], [\ldots, x^\pm_{i_1}[k_1]]]$, with $\sum_{s=1}^l \pm
\eps_{i_s} =  \al$. Assume $\al$ is not a root of $\G$ and let $l'$ be
the smallest integer such that $\sum_{s=1}^{l' + 1} \pm \eps_{i_s}$ is
not a root of $\G$. Let us show that each $[x^\pm_{i_{l' + 1}}[k_{l' +
1}], [\ldots, x^\pm_{i_1}[k_1]]]$ vanishes. It follows from the
fact that $j$ is an isomorphism when restricted to the parts of 
degree in $\Delta_+ \cup (-\Delta_+)$ 
that we may write each $[x^\pm_{i_{l'}}[k_{l'}], [\ldots,
x^\pm_{i_1}[k_1]]]$  as a linear combination $\sum_s \la_s
[x^\pm_{i_{l'}}[k^{(s)}(i_{l'})] \cdots x^\pm_{i_{1}}[k^{(s)}(i_{1})]]$,
where for each $s$, $i\mapsto k^{(s)}(i)$  is a map from
$\{1,\ldots,n\}$ to $\ZZ$, such that $k^{(s)}(i_{l' + 1}) =  k_{l' +
1}$. The defining relations for $\N_\pm$ hold among the 
$x^\pm_1[k^{(s)}_1], \ldots, x^\pm_n[k^{(s)}_n]$, therefore we have Lie
algebra maps from $\N_\pm$ to $\wt F^\pm$ sending each $x_i^\pm$ to
$x^\pm_i[k^{(s)}_i]$.  $[x^\pm_{i_{l' + 1}}[k^{(s)}(i_{l' + 1})]  \cdots
x^\pm_{i_{1}}[k^{(s)}(i_{1})]]$ is the image of zero by one of these
maps, and is therefore zero. It follows that $[x^\pm_{i_{l'}}[k_{l'}], [\ldots,
x^\pm_{i_1}[k_1]]]$ vanishes, so that $\Ker j$ is zero.
\hfill \qed \medskip

{\em End of proof of Proposition.} The restriction of $j$ to 
$\wt F_+$ coincides with the map $j_+$ define before Lemma 
\ref{cite:U}, therefore $j_+$ induces an isomorphism between 
$\wt F_+$ and $\N_+\otimes\CC[t,t^{-1}]$. 
\hfill \qed \medskip 

\subsubsection{Crossed product algebras $\cV^L$ and $\cS^L$}

For $k$ an integer $\leq 0$ and  $1\leq i \leq n$,
define  endomorphisms $\wt{h_i[k]}$  of $FO$ by  
\begin{equation} \label{def:hik}
(\wt{h_i[k]} f) (t^{(i)}_l) = [\sum_{j=1}^n \sum_{l=1}^{k_j} 
{{q^{kd_i a_{ij}} - q^{- kd_i a_{ij}}}\over{2\hbar kd_i}}
(t^{(j)}_l)^k ] f(t^{(i)}_l) 
\end{equation}
if $f\in FO_\kk$. The $\wt{h_i[k]}$ are 
derivations of $FO$. These derivations preserve $LV$, therefore 
they preserve $\langle LV\rangle$. 

Define $\cV^L$ and $\cS^L$ as the crossed product algebras of  $\langle
LV \rangle $ and $FO$ with the derivations $\wt{x[k]}$:  $\cV^L$, resp.\
$\cS^L$ is equal to $\langle LV\rangle \otimes \CC[h_i[k]^{\cV^L}, k\leq
0]$,  resp.\ $FO \otimes \CC[h_i[k]^{\cV^L}, k\leq 0]$; both spaces are
endowed with the products given by formula  (\ref{formula:crossed:pdt}),
where $x$ now belongs to  $\langle LV\rangle$, resp.\ $FO$ and the $X_s$
are replaced by $h_i[k]$,  $k\leq 0$.

Define $\hat\cV^L$ and $\cS^L$ as the partial $\hbar$-adic completions 
$\CC[h_i[k]^{\cV^L}, k<0] \otimes \langle LV\rangle 
\otimes_{\CC[[\hbar]]} \CC[h_i[0]^{\cV^L}][[\hbar]]$ and 
$\CC[h_i[k]^{\cV^L}, k<0] \otimes FO \otimes_{\CC[[\hbar]]}
\CC[h_i[0]^{\cV^L}][[\hbar]]$.  .

\begin{lemma} \label{LV:free}

1) The rules $\deg(h_i[k]) = 0,\deg(t_i^k) = \eps_i$ define gradings of 
$\langle LV\rangle$, $FO$, $\cV^L$, $\cS^L$, $\hat\cV^L$ and $\hat\cS^L$ by
$\NN^n$, which are compatible with the inclusions. For $X$ any of
these algebras, we denote by $X_\kk$ its homogeneous component of degree
$\kk$. $X$ is therefore the direct sum of the $X_\kk$. 

2) For any $\kk$, $\langle LV\rangle_\kk$ is a free
$\CC[[\hbar]]$-modules with a countable basis. 

3) For any $\kk$, $\cV^L_\kk$ and $\cS^L$ are free 
$\CC[[\hbar]][h_i[k]^{\cV^L}]$-modules; and   $\hat\cV^L_\kk$ and
$\hat\cS^L$ are free  $\CC[[\hbar]][h_i[k]^{\cV^L},k<0]
\otimes_{\CC[[\hbar]]} \CC[h_i[0]^{\cV^L}][[\hbar]]$-modules. 

\end{lemma}

{\em Proof.} 1) is clear.  $\langle LV\rangle_\kk$ is a 
$\CC[[\hbar]]$-submodule of $FO_\kk$, and by  Lemma \ref{free:inf:dim},
it is a  free $\CC[[\hbar]]$-module with a countable basis.  This shows
2). 3) is a direct consequence of 2). \hfill \qed \medskip

\subsubsection{Ideals and completions}

Define for $N$ positive integer, $I_N$ as the left ideal of $\langle
LV\rangle$ generated by the elements $(t_i^k)$ of $FO_{\eps_i}$, $k\geq
N$,  $i = 1,\ldots,n$. Define $\cI_N$ and $\hat \cI_N$ as the left
ideals of $\cV^L$ and $\hat\cV^L$ generated by the same family.  For
$s\geq 0$, set $I_N^{(s)} = \hbar^{-s} (I_N \cap \hbar^s \langle
LV\rangle)$,  and $I_N^{(\infty)} = \cup_{s\geq 0} I_N^{(s)}$; define
$\cI_N^{(s)}$,  $\cI_N^{(\infty)}$ and $\hat\cI_N^{(s)}$,
$\hat\cI_N^{(\infty)}$ in the same way.

For any integer $a$, define $LV^{\geq a}$ as the subspace of $FO$ equal to 
the direct sum $\oplus_{i= 1}^n t_i^a \CC[[\hbar]][t_i]$ and let 
$\langle LV^{\geq a}\rangle$ be the subalgebra of $FO$ generated by 
$LV^{\geq a}$. Define $I_N^{(0),\geq a}$ as the left ideal of $\langle LV^{\geq a}\rangle$
generated by the $t_i^{k},k\geq N$ and 
$I_N^{\geq a}$ as the ideal of $\langle LV^{\geq a}\rangle$ formed of the elements 
$x$ such that for some $k\geq 0$, $\hbar^k x$ belongs to $I_N^{(0),\geq a}$.

For any integer $a$ and $\kk$ in $\NN^r$, define $FO^{\geq a}$ as the subspace of 
$FO_\kk$
consisting of the rational functions 
$$
g(t^{(i)}_\al) = {1\over{\prod_{i = 1}^n} \prod_{1\leq \al \leq k_i, 1\leq \beta\leq k_j (t^{(i)}_\al
- t^{(j)}_\beta)}} f(t^{(i)}_\al),  
$$
where the $f(t^{(i)}_\al)$ have degree $\geq a$ in each variable $t^{(i)}_\al$
and the total degree of $g$ is $\geq (\sum_i k_i)a$. 
Set $FO^{\geq a} = \oplus_{\kk\in\NN^n} FO^{\geq a}_\kk$. Then 
$FO^{\geq a}$ is a subalgebra of $FO$. 

Define $\cI_{(N)}^{\geq a}$ as the set of elements of $FO^{(\geq a)}$, where 
$f(t^{(i)}_\al)$ has total degree $N$ in the variables $t^{(i)}_\al$, and let 
$\cI_N^{\geq a}$ be the direct sum $\oplus_{k\geq N} \cI_N^{\geq a}$.

\begin{lemma} \label{ladino}
For $(J_N)_{N>0}$ a family of left ideals of some algebra $A$, say that 
$(J_N)_{N>0}$ has property $(*)$ if for any integer $N>0$ and element
$a$ in $A$,  there is an integer  $k(a,N)>0$ such that $J_N a\subset
J_{k(a,N)}$ for any $N$ large enough, and $k(N,a)$ tends to infinity
with $N$, $a$ being fixed. Then the inverse limit $\limm_{\leftarrow N}
A/J_N$ has an algebra structure. 

Say that $J_N$ has property $(**)$ if for any integer $N>0$ and element
$a$ in $A$,  there is are integer  $k'(a,N)$ and $k''(a,N)>0$ such that $J_N a\subset
J_{k'(a,N)}$ and $a J_N\subset
J_{k''(a,N)}$ for any $N$ large enough, and $k(N,a)$ tends to infinity
with $N$, $a$ being fixed. In that case also, the inverse limit $\limm_{\leftarrow N}
A/J_N$ has an algebra structure.

1) The family $(I^{(\infty)}_N)_{N>0}$ of ideals of $\langle LV\rangle$
has property $(*)$; 

2) the family $(\cI^{(\infty)}_N)_{N>0}$ of ideals of $\cV^L$ has property $(*)$; 

3) the family $(\hat\cI^{(\infty)}_N)_{N>0}$ of ideals of $\hat\cV^L$ has 
property $(*)$; 

4) the family $(I_N^{\geq a})_{N>0}$ of ideals of $\langle LV^{\geq a} \rangle$
has property $(*)$; 

5) the family $(\cI^{\geq a}_N)_{N>0}$ of ideals of $FO^{\geq a}$ has property $(**)$. 
\end{lemma}

{\em Proof.} Set for any $a$ in $A$ and $N>0$, $k'(a,N) =$ inf 
$\{k| J_N a \subset J_k\}$; then $k'(N,a)$ tends to infinity with $N$, 
$a$ being fixed and we have $k'(N,a')\geq$ inf $(k(N,a), p)$ if 
$a'$ belongs to $a + J_p$. 
 
An element of $\limm_{\leftarrow N} A / J_N$ is a family $(a_N)_{N>0}$, 
$a_N\in A / J_N$, such that $a_{N+1} + J_N = a_N$. For $a = (a_N)_{N>0}$
and $b = (b_N)_{N>0}$ in $\limm_{\leftarrow N} A / J_N$,  choose
$\beta_N$ in $b_N$ and let $N'(N,\beta_N)$ be the smallest integer $N'$ 
such that $k'(N',\beta_N)\geq N$; $N'(N,\beta_N)$ is independent of the
choice of  $\beta_N$, we denote it $N'(N,b)$.  Choose then $\al_N$ in
$a_{N'(N,b)}$; then $\al_N\beta_N + J_N$ is independent of the choice of
$\al_N$ and $\beta_N$; one checks that $\al_{N+1}\beta_{N+1} + J_N =
\al_N\beta_N + J_N$, so that $(\al_N\beta_N + J_N)_{N>0}$ defines an
element of $\limm_{\leftarrow N} A / J_N$. The product $ab$ is defined
to be this element.
The construction is similar in the case of property $(**)$.

1) The equality 
\begin{align*}
& t_i^k * t_j^l = q^{-n d_i a_{ij}} t_i^{k-n} * t_j^{l+n}
+ q^{- d_i a_{ij}} t_j^l * t_i^k  \\ & + \sum_{n' = 1}^{n-1}
(q^{ - (n' + 1) d_i a_{ij}}  - q^{ - (n' - 1) d_i a_{ij}} ) 
t_j^{l+n'} * t_i^{k-n'} - q^{- (n-1)d_i a_{ij}} t_j^{l+n} * 
t_i^{k-n} ,  
\end{align*}
where $n = N - l$ and $k\geq 2N - l$,  implies that if $k\geq 2N - l$,
$t_i^k * t_j^l$ belongs to $I_N$.  It follows that $I_{2N - l} * t_j^l
\subset I_N$. Set then  $k(t_j^l,N) = [{1\over 2}(N + l/2)] + 2$; for
$a$ in $\langle LV \rangle$,  and any decomposition $dec$ of $a$ as a
sum $\sum_{(j_i),(l_i)} \la_{(j_i),(l_i)} t_{j_1}^{l_1} * \cdots *
t_{j_p}^{l_p}$, define $k(a,N,dec)$ as the smallest of integers
$k(t_{j_p}^{l_p}, \cdots, k(t_{j_1}^{l_1},N))$; finally,  define
$k(a,N)$ as the largest of all $k(a,N,dec)$. 
The family $(I_N)_{N>0}$ has  property $(*)$, with this function
$k(a,N)$. Then for any $a$ in $\langle LV\rangle$,   $(I_N \cap \hbar^s
\cV^L) a \subset \hbar^s \cV^L \cap I_{k(a,N)}$, so that  $I_N^{(s)} a
\subset I^{(s)}_{k(a,N)}$. It follows that we have also    
$I_N^{(\infty)} a \subset I^{(\infty)}_{k(a,N)}$, so that the families 
$(I_N^{(s)})_{N>0}$ and $(I_N^{(\infty)})_{N>0}$ have property $(*)$. 

2) follows from the fact that $I_{N} h_j[l]^{\cV^L}
\subset I_{N + l}$. 

3) follows from the fact that for $a$ any element of $\CC[h_i[0]^{\cV^L}][[\hbar]]$, 
we have $\hat\cI_N^{(\infty)}a \subset \hat\cI_N^{(\infty)}$.

4) is proved in the same way as 1) ans 2). 

5) For $\kk$ in $\NN^n$, set $|\kk| = \sum_{i = 1}^n k_i$. Then if $f$ belongs
to $FO^{\geq a}_{\kk}$, we have $f * \cI_N^{\geq a} \subset \cI^{\geq a}_{N + a|\kk|}$ 
and $\cI_N^{\geq a} * f\subset \cI^{\geq a}_{N + a|\kk|}$. Therefore the family 
$(\cI_N^{\geq a})_N$ has property $(**)$. 
\hfill \qed \medskip 

It follows that the inverse limits
$\limm_{\leftarrow N} \langle LV\rangle / I_N^{(\infty)}$, 
$\limm_{\leftarrow N} \cV^L / \cI_N^{(\infty)}$, 
$\limm_{\leftarrow N} \hat\cV^L / \hat\cI_N^{(\infty)}$, 
$\langle LV^{\geq a}\rangle / I_N^{\geq a}$ and 
$\limm_{\leftarrow N}FO^{\geq a} / \cI_N^{\geq a}$ 
have algebra structures.
Moreover, as we have $\cI_N^{(\infty)}\cap \langle LV\rangle = I_N^{(\infty)}$ and
$\hat\cI_N^{(\infty)}\cap \cV^L = \cI_N^{(\infty)}$, we have natural algebra inclusions 
$$
\limm_{\leftarrow N} \langle LV\rangle / I_N^{(\infty)} \subset 
\limm_{\leftarrow N} \cV^L / \cI_N^{(\infty)} \subset 
\limm_{\leftarrow N} \hat\cV^L / \hat\cI_N^{(\infty)}.
$$
Moreover, there exists a function $\phi(\kk,N)$, tending to infinity with $N$, 
such that $\cI_{N,\kk}^{\geq a} \cap \langle LV^{\geq a}
\rangle \subset I_{\phi(\kk,n)}$. Indeed, if the $k_i$ are $\geq a$ and 
$t_{i_1}^{k_1} * \cdots * t_{i_{l}}^{k_l}$ belongs to $\cI^{\geq a}_{N,\kk}$
($l = |\kk|$), then $k_1 + \cdots + k_l \geq N$ so that one of the $k_i$ is 
$\geq N/l$. The statement then follows from the proof of the above Lemma, 1).  
It follows that we have an algebra inclusion 
$$\limm_{\leftarrow N}
\langle LV^{\geq a} \rangle / I^{\geq a}_N \subset \limm_{\leftarrow N}
FO^{\geq a} / \cI_N^{\geq a}. 
$$
 
If a family $(J_N)_{N>0}$ of left ideals of the algebra $A$ has
property $(*)$, and $B$ is any algebra, the family $(J_N \otimes B)_{N>0}$
also satisfies $(*)$; therefore the inverse limit $\limm_{\leftarrow N}
(A\otimes B) / (J_N \otimes B)$ has an algebra structure. It follows
that we have algebra structures on $\limm_{\leftarrow N} \langle
LV\rangle \otimes A/ I_N^{(\infty)}\otimes A$,   $\limm_{\leftarrow N} \cV^L
\otimes A/ \cI_N^{(\infty)}\otimes A$ and  $\limm_{\leftarrow N} \hat\cV^L \otimes
A/ \hat\cI_N^{(\infty)}\otimes A$ for any  algebra $A$.   

\subsubsection{Topological Hopf structures on $\cV^L$ and $\cS^L$}

For $i = 1,\ldots,n$ and $l\geq 0$, define $K_i[-l]$ as the element of $\cV^L$
(or $\cS^L$)   
$$
K_i[-l] = e^{ - \hbar d_i h_i[0]^{\cV^L}} S_l(-2\hbar d_ih_i[k]^{\cV^L},k<0),
$$ 
where 
$S_l(z_1,z_2,\ldots)$ are the Schur polynomials in variables $(z_i)_{i<0}$, 
which are determined by the relation $\exp(\sum_{i<0} z_it^{-i}) 
= \sum_{l\leq 0} S_l(z_k) t^{-l}$. 

\begin{prop} \label{Hopf:FO}
There is unique graded algebra morphism $\Delta_{\cS^L}$ from $\cS^L$ 
to $\limm_{\leftarrow N}(\cS^L \otimes_{\CC[[\hbar]]}\hat\cS^L) / 
(\cI_N^{\cS} \otimes_{\CC[[\hbar]]} \hat\cS^L )$, such that 
$$
\Delta_{\cS^L}(h_i[k]) = h_i[k] \otimes 1 + 1\otimes h_i[k]   
$$
for $1\leq i \leq n$ and $k\leq 0$, and 
its restriction to $\cS^L_\kk$ is the direct sum of the maps
$\Delta_{\cS^L}^{\kk',\kk''}: FO^{\geq a}_\kk \to 
\limm_{\leftarrow N}(\cI^{\geq a}_{\kk'} \otimes_{\CC[[\hbar]]} \hat\cS^L_{\kk''}) 
/ [FO_{\kk',N}^{\geq a}\otimes_{\CC[[\hbar]]} \hat \cS^L_{\kk''} ] $, 
where $\kk = \kk' + \kk''$, defined by 
$$
\Delta_{\cS^L}^{\kk',\kk''}(P) = \sum_{p_1,\ldots,p_{N'}\geq 0}
\left( \prod_{i = 1}^{N'} u_i^{p_i} P'_{\al}(u_1,\ldots,u_{N'}) 
\right) \otimes \left( 
P''_\al(u_{N' + 1}, \ldots, u_N) 
\prod_{i = 1}^{N'} K_{\eps(i)}[-p_i]
\right),  
$$
where $N' = \sum_{i=1}^n k'_i$, $N'' = \sum_{i=1}^n k''_i$, 
$N = N' + N''$, the arguments of the functions in $FO_{\kk'}$ 
and $FO_{\kk''}$ are respectively $(t^{(i)}_j)_{1\leq i\leq n, 
1\leq j\leq k'_i}$ and $(t^{\prime (i)}_j)_{1\leq i\leq n, 
1\leq j\leq k''_i}$; we set 
$$
(u_1,\ldots, u_{k'_1}) = (t^{(1)}_1, \ldots, t^{(1)}_{k'_1}), \ldots, 
(u_{k'_1 + \ldots + k'_{n-1} + 1}, \ldots, u_{N'}) = 
(t^{(n)}_1, \ldots, t^{(n)}_{k'_n}), 
$$
$$
(u_{N' + 1},\ldots, u_{N' + k''_1}) = (t^{\prime (1)}_1, \ldots, t^{\prime 
(1)}_{k''_1}), \ldots, 
(u_{N' + k''_1 + \ldots + k''_{n-1} + 1}, \ldots, u_{N}) = 
(t^{\prime (n)}_1, \ldots, t^{\prime (n)}_{k''_n}), 
$$
$$
t^{(i)}_{k'_i + j} = t^{\prime (i)}_j, \quad j = 1,\ldots, k''_i, 
$$
$$
(t_1,\ldots, t_{k_1}) = (t^{(1)}_1, \ldots, t^{(1)}_{k_1}), \ldots, 
(t_{k_1 + \ldots + k_{n-1} + 1}, \ldots, t_{N}) = 
(t^{(n)}_1, \ldots, t^{(n)}_{k_n}), 
$$
and 
\begin{equation} \label{weil}
\sum_{\al} P'_\al(u_1,\ldots,u_{N'})  P''_\al(u_{N'+1},\ldots,u_{N})
= P(t_1,\ldots,t_N) \prod_{1\leq l\leq N', N' + 1\leq l'\leq N}
{{u_{l'} - u_l}\over{q^{\langle \eps(l),\eps(l') \rangle} u_{l'} - 
u_l}},   
\end{equation}
for $l$ in $\{1,\ldots,N'\}$, resp.\ $\{N' + 1,\ldots,N\}$, 
$\eps(l)$ is the element of $\{1,\ldots,n\}$ 
such that $u_l = t^{(\eps(l))}_j$ , resp.\ $u_l = t^{\prime(\eps(l))}_j$ 
for some $j$; in (\ref{weil}), the ratios are expanded for $u_l << u_{l'}$. 
\end{prop}

{\em Proof.} Define $FO^{(2)}$ as the $(\NN^r)^2$-graded $\CC[[\hbar]]$-module
$FO^{(2)} = \oplus_{\kk,\kk'\in\NN^r} FO^{(2)}_{\kk,\kk'}$, where 
\begin{align*}
&  FO^{(2)}_{\kk,\kk'} = 
\\ & {1\over{
\prod_{i<j, 1\leq \al \leq k_i,
 1\leq \beta \leq k_j} ( t^{(i)}_\al  - t^{(j)}_{\beta} ) 
\prod_{i<j, 1\leq \al \leq k'_i, 1\leq \beta \leq k'_j
} ( u^{(i)}_\al  - u^{(j)}_{\beta} ) 
\prod_{i,j, 1\leq \al \leq k_i, 1\leq \beta \leq k'_j} (t^{(i)}_\al - u^{(j)}_\beta) }}
\cdot \\ & \cdot \CC[[\hbar]][(t^{(i)}_j)^{\pm1}, (u^{(i)}_{j'})^{\pm1}, j = 1,\ldots, k_i , 
j' = 1,\ldots, k'_i ]^{\prod_{i=1}^n\gotS_{k_i} \times \gotS_{k'_i}} , 
\end{align*}
where the groups $\gotS_{k_i}$ and $\gotS_{k'_i}$ act by permutation of the
variables $t^{(i)}_j$ and $u^{(i)}_j$. Define on $FO^{(2)}$, the graded composition 
map $*$ as follows: for $f$ in $FO^{(2)}_{\kk,\kk'}$ and $g$ in $FO^{(2)}_{\bl,\bl'}$, 
\begin{align*}
& (f*g)(t_1,\ldots, t_P, u_1, \ldots, u_{P'}) = 
\Sym_{t^{(1)}_j}  \cdots \Sym_{t^{(n)}_j}   
\Sym_{u^{(1)}_j}  \cdots \Sym_{u^{(n)}_j}
\\ & ( 
\prod_{1\leq i \leq N, N+1\leq j \leq P} { {q^{ \langle \eps_t(i), \eps_t(j)  \rangle } 
t_i - t_j }\over{t_i - t_j } } 
\prod_{1\leq i \leq N', N+1\leq j \leq P} { {q^{ \langle \eps_u(i), \eps_t(j)  \rangle } 
u_i - t_j }\over{u_i - t_j }} 
\\ & \prod_{1\leq i \leq N, N'+1\leq j \leq P'} { {q^{ \langle \eps_t(i), \eps_u(j)  \rangle } 
t_i - u_j }\over{t_i - u_j }} 
\prod_{1\leq i \leq N', N'+1\leq j \leq P'} {{q^{ \langle \eps_u(i), \eps_u(j)  \rangle } 
u_i - u_j }\over{u_i - u_j }}
\\ &  
f(t_1,\ldots, t_N, u_1,\ldots, u_M)  
g(t_{N+1},\ldots, t_P, u_{M+1},\ldots, u_P)  ) ,    
\end{align*}
where $N = \sum_i k_i, N' = \sum_i k'_i, M = \sum_i l_i, M' = \sum_i l'_i, 
P = N + M, P' = N' + M'$, and $(t_1,\ldots, t_{k_1 + k'_1}) = (t^{(1)}_1, \ldots, 
t^{(1)}_{k_1 + k'_1})$, ..., 
$(t_{k_1 + \cdots + k'_{n-1} + 1 },\ldots, t_{P}) = (t^{(n)}_1, \ldots, 
t^{(n)}_{k_n + k'_n})$, and  $(u_1,\ldots, u_{l_1 + l'_1}) = (u^{(1)}_1, \ldots, 
u^{(1)}_{l_1 + l'_1})$, ..., 
$(u_{l_1 + \cdots + l'_{n-1} + 1 },\ldots, u_{P'}) = (u^{(n)}_1, \ldots, 
u^{(n)}_{l_n + l'_n})$. We set $\eps_x(\al) = i$ if $x_\al = x^{(i)}_j$ for some $j$
($x$ is $t$ or $u$). It is easy to check that $µ$ defines an algebra structure on 
$FO^{(2)}$. 

Let $\Delta_{FO}$ be the linear map from $FO$ to $FO^{(2)}$, which maps $FO_{\kk}$
to $\oplus_{\kk' + \kk'' = \kk} FO^{(2)}_{\kk',\kk''}$ as follows
$$
\Delta_{FO}(P)(t^{(1)}_1,\ldots, t^{(n)}_{k'_n},   
u^{(1)}_1,\ldots, u^{(n)}_{k''_n} ) = 
P(t^{(1)}_1,\ldots, t^{(1)}_{k'_1},  u^{(1)}_1,\ldots, u^{(1)}_{k''_1}, 
t^{(2)}_1, \ldots, u^{(n)}_{k''_n}) .    
$$
Then it is immediate that $\Delta_{FO}$ defines an algebra morphism. 

Consider the $(\NN^n)^2$-graded map $\mu : FO^{(2)} \to 
\cup_a \limm_{\leftarrow N}
(FO^{\geq a} \otimes_{\CC[[\hbar]]} \hat\cS^L)
 / [FO_N^{\geq a}\otimes_{\CC[[\hbar]]} \hat
\cS^L]$
defined by 
$$
\mu(P) = \sum_{p_1,\ldots,p_{N'}\geq 0} (t_1^{p_1}\cdots t_{N'}^{p_{N'}}
\bar P_\al) \otimes (\bar P'_\al K_{\eps_t(1)}[-p_1] \cdots K_{\eps_t(N')}[-p_{N'}])
$$
if $P(t_1,\ldots, t_N, u_1,\ldots, u_{N'})$ belongs to $FO^{(2)}_{\kk,\kk'}$, and we set 
$$
P(t_1,\ldots, u_{N'}) \prod_{1\leq l \leq N, 1\leq l'\leq N'} 
{{t_l - u_{l'}}\over{q^{\langle \eps_t(l), \eps_u(l')\rangle} t_l - u_{l'}}} 
= \sum_{\al} \bar P_\al(t_1,\ldots,t_{N})  \bar P'_\al(u_1,\ldots,u_{N'}) 
$$
(the expansion is for $u_{l'} << t_l$). 

Let $P,Q$ belong to $FO_{\kk,\bl}$ and $FO_{\kk',\bl'}$. $\mu(P)\mu(Q)$ is 
equal to 
\begin{align} \label{pioche}
& \sum_{\al,\beta} \left( t_1^{p_1} \cdots t_N^{p_N}
P_\al(t_1,\cdots,t_N) \otimes  P'_\al(u_1,\ldots,u_P)
K_{\eps_t(1)}[-p_1]\cdots  K_{\eps_t(N)}[-p_N]\right)  \cdot \\ &
\nonumber \cdot  \left( t_1^{\prime p'_1} \cdots t_{N'}^{\prime p'_{N'}}
Q_\al(t'_1,\cdots,t'_{N'}) \otimes  Q'_\al(u'_1,\ldots,u'_{P'})
K_{\eps_{t'}(1)}[-p'_1]\cdots  K_{\eps_{t'}(1)}[-p'_{N'}]\right) ; 
\end{align}
since 
$$
(\sum_{i\geq 0} t^i K_\al[-i]) 
Q(u_1,\ldots,u_P) =  \prod_{i=1}^P {{u_i - 
q^{\langle \eps_u(i), \al \rangle }t}\over
{ q^{\langle \eps_u(i), \al \rangle } u_i - t}}
Q(u_1,\ldots,u_P) (\sum_{i\geq 0} t^i K_\al[-i]) , 
$$
(\ref{pioche}) is equal to 
\begin{align} \label{pioche'}
& \sum_{\gamma}
\sum_{\al,\beta} \left( f_\gamma(t_1,\ldots, t_N) t_1^{p_1} \cdots t_N^{p_N} 
P_\al(t_1,\cdots,t_N) \otimes 
P'_\al(u_1,\ldots,u_P) 
\right) \cdot 
\\ & \nonumber \cdot 
\left( t_1^{\prime p'_1} \cdots t_{N'}^{\prime p'_{N'}} Q_\al(t'_1,\cdots,t'_{N'}) \otimes 
g_\gamma(u'_1,\ldots,u'_{P'}) 
Q'_\al(u'_1,\ldots,u'_{P'}) K_{\eps_t(1)}[-p_1]\cdots 
\right. \\ &  \nonumber
K_{\eps_t(N)}[-p_N] K_{\eps_{t'}(1)}[-p'_1]\cdots  K_{\eps_{t'}(1)}[-p'_{N'}] ) ,  
\end{align}
with 
$$
\sum_\gamma f_\gamma(t_1,\ldots, t_N) g_\gamma(u'_1,\ldots, u'_{P'}) = 
\prod_{1\leq i \leq N, 1\leq j\leq P'} 
{{u_j - q^{\langle \eps_u(j), \al \rangle }t_i}\over
{q^{\langle \eps_u(j), \al \rangle } u_j - t_i}}.  
$$
After some computation, one finds that (\ref{pioche'}) coincides with 
$\mu(PQ)$. Therefore $\mu$ is an algebra morphism. 

It follows that the composition $\mu\circ \Delta_{FO}$ is an algebra morphism. 
This composition coincides with $\Delta_{\cS^L}$, which is therefore an algebra 
morphism. 
\hfill\qed\medskip 

\begin{remark}
$\Delta_{\cS^L}^{\kk',\kk''}$ may also be expressed by the formula
$$
\Delta_{\cS^L}^{\kk',\kk''}(P) = \sum_{p_1,\ldots,p_{N'}\geq 0}
\left( \prod_{i = 1}^{N'} u_i^{p_i} P'''_{\al}(u_1,\ldots,u_{N'}) 
\right) \otimes \left( \prod_{i = 1}^{N'} K_{\eps(i)}[-p_i]
P''''_\al(u_{N' + 1}, \ldots, u_N) 
\right),  
$$
where 
$$ 
\sum_{\al} P'''_\al(u_1,\ldots,u_{N'})  P''''_\al(u_{N'+1},\ldots,u_{N})
= P(t_1,\ldots,t_N) \prod_{1\leq l\leq N', N' + 1\leq l'\leq N}
{{u_{l'} - u_l}\over{u_{l'} - q^{\langle \eps(l),\eps(l') \rangle} 
u_l}} . 
$$
\hfill \qed\medskip 
\end{remark}

\begin{cor}
There is a unique algebra morphism $\Delta_{\cV^L}$ 
from $\cV^L$
to $\limm_{\leftarrow N} (\cV^L \otimes_{\CC[[\hbar]]} \hat \cV^L)  / 
(\cI_N^{(\infty)}\otimes_{\CC[[\hbar]]}  \hat\cV^L)$, such that 
$$
\Delta_{\cV^L}(h_i[k]^{\cV^L}) = h_i[k]^{\cV^L}\otimes 1 + 1 \otimes h_i[k]^{\cV^L}, 
\quad 1\leq i \leq n, \quad k\leq 0, 
$$
$$ 
\Delta_{\cV^L}(t_i^k) = \sum_{l\geq 0}t_i^{k+l}\otimes K_i[-l]
+ 1\otimes t_i^{k}
\quad 1\leq i \leq n, \quad k\in\ZZ. 
$$
\end{cor}

{\em Proof.}  For each $i$, $\Delta_{\cS^L}$  maps $FO_{\eps_i}$ to
$\limm_{\leftarrow N} (\cV^L \otimes_{\CC[[\hbar]]} \hat \cV^L)  / 
(\cI_N^{(\infty)}\otimes_{\CC[[\hbar]]}  \hat\cV^L)$, because  
$\cI_{N,\kk}^{\geq a} \cap \langle LV^{\geq a}
\rangle \subset I_{\phi(\kk,n)}$. It follows  that
$\Delta_{\cS^L}$  maps $\cV^L$ to the same space. Call $\Delta_{\cV^L}$
the restriction  of $\Delta_{\cS^L}$ to $\cV^L$. This restriction is
clearly  characterized by its values on $h_i[k]$ and $t_i^k$.  \hfill
\qed\medskip

\subsubsection{Construction of a Hopf algebra structure}

Recall that we showed in Lemma \ref{LV:free} that $\cV^L$ is a free 
$\CC[[\hbar]]$-module. Let us set $\cV^L_0 = \cV^L / \hbar\cV^L$.  

Let $\overline{\cI}_N$ be the image of the ideal $\cI_N^{(\infty)}$ of $\cV^L$ by the
projection from $\cV^L$ to $\cV^L_0$. By Lemma  \ref{ladino}, 2), and since
the map from $\cV^L$ to $\cV^L_0$ is surjective,  the ideals
$\ol{\cI}_N$ have property $(*)$, so that   $\cap_{N>0}\ol{\cI}_N$  is
a two-sided ideal of $\cV^L_0$. Define $\cW_0$ as the quotient  algebra
$\cW_0 = \cV^L_0 / \cap_{N>0}\ol{\cI}_N$.  We are going to define a
Hopf algebra structure on $\cW_0$.

Since the $\cV^L_0 \otimes\ol{\cI}_N$ have property $(*)$,  
$\limm_{\leftarrow N} (\cV^L_0 \otimes \cV^L_0)  
/ (\cV^L_0 \otimes \ol{\cI}_N)$ has an algebra structure.
Moreover, the projection 
\begin{align*}
& [\limm_{\leftarrow N} (\cV^L \otimes_{\CC[[\hbar]]} \cV^L)  
/ (\cV^L \otimes_{\CC[[\hbar]]} \cI_N^{(\infty)})] 
/ \hbar [\limm_{\leftarrow N} (\cV^L \otimes_{\CC[[\hbar]]} \cV^L)  
/ (\cV^L \otimes_{\CC[[\hbar]]} \cI_N^{(\infty)})] 
\\ & \to 
\limm_{\leftarrow N} (\cV^L_0 \otimes \cV^L_0)  
/ (\cV^L_0 \otimes \ol{\cI}_N) , 
\end{align*}
is an algebra isomorphism. 
$\Delta_{\cV^L}$ 
induces therefore an algebra morphism $\Delta_{\cV^L_0}$ from $\cV^L_0$
to $\limm_{\leftarrow N} (\cV^L_0 \otimes \cV^L_0)  
/ (\cV^L_0 \otimes \ol{\cI}_N ) $. 
 
On the other hand, the $\ol{\cI}_N \otimes \cV^L_0  + 
\cV^L_0 \otimes \ol{\cI}_N$ have property $(*)$, so that  
$\limm_{\leftarrow N} (\cV^L_0 \otimes \cV^L_0)  
/ (\ol{\cI}_N \otimes \cV^L_0  + \cV^L_0 \otimes \ol{\cI}_N)$
has an algebra structure. The composition of  $\Delta_{\cV^L_0}$ 
with the projection 
$$
\lim_{\leftarrow N} (\cV_0^L \otimes \cV_0^L) / 
(\cV_0^L \otimes \ol{\cI}_N) \to  
\lim_{\leftarrow N} (\cV_0^L \otimes \cV_0^L) / 
(\ol{\cI}_N \otimes \cV_0^L  + \cV_0^L  \otimes \ol{\cI}_N ) 
$$
then yields an algebra morphism $p'\circ \Delta_{\cV_L^0}$ from $\cV_L^{(0)}$ to 
$\limm_{\leftarrow N} (\cV^L_0 \otimes \cV^L_0)  
/ (\ol{\cI}_N \otimes \cV^L_0  + \cV^L_0 \otimes \ol{\cI}_N)$. 

We have for any $k\geq N$, $\Delta_{\cV}(t_i^k)\in \cI_N\otimes\hat\cV +
\cV\otimes \hat\cI_N$, therefore 
$\Delta_{\cV^L}(\cI_N) \subset \cI_N \otimes_{\CC[[\hbar]]} \hat\cV
+ \cV\otimes_{\CC[[\hbar]]} \hat\cI_N$, therefore 
$\Delta_{\cV^L_0}(\ol{\cI}_N ) \subset \ol{\cI}_N \otimes \cV_0 +
\cV_0\otimes \ol{\cI}_N$.  It follows that $p'\circ \Delta_{\cV_L^0}$
maps the intersection $\cap_{N> 0}\ol{\cI}_N$ to the kernel of
the projection  $\limm_{\leftarrow N} (\cV^L_0 \otimes \cV^L_0)   /
(\cV^L_0 \otimes \ol{\cI}_N) \to  \limm_{\leftarrow N} (\cV^L_0
\otimes \cV^L_0)   / (\ol{\cI}_N \otimes \cV^L_0  + \cV^L_0 \otimes
\ol{\cI}_N)$. 

We have then an  algebra morphism $\tilde\Delta_{\cW_0}$ 
from $\cW_0$ to  $\limm_{\leftarrow N}
(\cV^L_0 \otimes \cV^L_0)   / (\ol{\cI}_N \otimes \cV^L_0  + \cV^L_0
\otimes \ol{\cI}_N)$. 

In Prop.\ \ref{prop:morphism}, we defined a surjective algebra morphism
$i_\hbar$ from $U_\hbar L\B_+$ to $\cV^L$. It induces an algebra
morphism $i$ from $U L\B_+$ to  $\cV^L_0$, which is also surjective.

Let $T$ be the free algebra generated by the $h_i[k]^{(T)},
k\leq 0, i = 1,\ldots,n$ and $e_i[k]^{(T)}, i =
1,\ldots,n, k\in \ZZ$. We have a  natural projection of $T$ on $U
L\B_+$, sending each  $x[k]^{(T)}$ to $x\otimes t^k$; composing it with
$\iota$,  we get a surjective algebra morphism $\pi$ from $T$ to
$\cV^L_0$. 
 
We have a unique algebra morphism $\Delta_T:T\to T\otimes T$, 
such that $\Delta_T(x_i[k]^{(T)}) = x_i[k]^{(T)} \otimes 1 + 
1\otimes x_i[k]^{(T)}$. 

\begin{lemma}
Let $\pi_{\cV_0\to \cW_0}$ be the natural projection from $\cV_0$ to $\cW_0$. 
We have the identity 
\begin{equation} \label{basic} 
 \tilde\Delta_{\cW_0} \circ (\pi_{\cV_0\to\cW_0}\circ\pi) 
= (\nu\circ(\pi_{\cV_0\to \cW_0}\otimes\pi_{\cV_0\to \cW_0}))
\circ\Delta_T, 
\end{equation}
where $\nu$ is the natural projection from $(\cV^L_0 \otimes \cV^L_0)$ to
$\limm_{\leftarrow N} (\cV^L_0 \otimes \cV^L_0) / (\ol{\cI}_N \otimes
\cV^L_0  + \cV^L_0 \otimes \ol{\cI}_N)$. 
\end{lemma}

{\em Proof.} The two sides are algebra morphisms from $\cW_0$ to
$\limm_{\leftarrow N} (\cV^L_0 \otimes \cV^L_0) / (\ol{\cI}_N \otimes
\cV^L_0  + \cV^L_0 \otimes \ol{\cI}_N)$. The identity is satisfied on
generators of $\cW_0$, therefore it is true.  \hfill \qed\medskip 

Let $J$ be the kernel of the projection $\pi_{\cV_0\to\cW_0}\circ \pi$
from  $T$ to $\cW_0$. It follows from (\ref{basic}) that $\Delta_T(J)$
in contained in the  kernel of $\nu\circ(\pi\otimes\pi)$, which is the
preimage by $\pi\otimes\pi$ of $\Ker\nu$; $\Ker\nu$ is equal to $\cap_{N>0} 
(\ol{\cI}_N \otimes \cV^L_0  + \cV^L_0 \otimes \ol{\cI}_N)$, which is 
$(\cap_{N>0} \ol{\cI}_N) 
\otimes \cV^L_0  + \cV^L_0 \otimes (\cap_{N>0} \ol{\cI}_N)$. Therefore 
$(\pi\otimes\pi)^{-1}(\Ker\nu)$ is  $J\otimes T + T\otimes J$. 
We have shown that $\Delta_T(J)\subset J\otimes T + T \otimes J$. 
We have shown: 

\begin{prop}
 $\Delta_T$ induces a cocommutative Hopf algebra structure on 
$T / J = \cW_0$. 
\end{prop}

We will denote by $\Delta_{\cW_0}$ the coproduct induced
by $\Delta_T$ on $\cW_0$.

\subsubsection{Compatibility of $\Delta_{\cW_0}$ with 
$\Delta_{\cV^L_0}$}

Recall that $\Delta_{\cV^L_0}$ is an algebra 
morphism from $\cV^L_0$ to $\lim_{\leftarrow N} 
(\cV^L_0\otimes\cV^L_0)  / (\ol{\cI}_N \otimes \cV^L_0)$. 
Let us denote by $\Delta_{\cV_0^L;N}$
the induced map from $\cV_0^L$ to $(\cV_0^L\otimes\cV^L_0) 
/ (\ol{\cI}_N \otimes \cV^L_0)$. 
We have seen that for any 
integer $p>0$, $\Delta_{\cV^L_0;N}(\ol{\cI}_p)$ is contained in  
the image of $\ol{\cI}_p \otimes \cV^L_0 + \cV^L_0 \otimes\ol{\cI}_p$
by the projection map $\cV^L_0\otimes \cV^L_0 \to (\cV^L_0\otimes\cV^L_0) 
/ (\ol{\cI}_N \otimes \cV^L_0)$. This image is
$$
[(\ol{\cI}_p + \ol{\cI}_N) \otimes \cV^L_0 + \cV^L_0 \otimes \ol{\cI}_p] 
/ (\ol{\cI}_N \otimes \cV^L_0) . 
$$
Therefore, $\Delta_{\cV^L_0;N}(\cap_{p>0}\ol{\cI}_p)$ is contained in  
the intersection of these spaces, which is 
$$
[ \ol{\cI}_N \otimes \cV^L_0 + \cV^L_0 \otimes (\cap_{p>0}\ol{\cI}_p) ]
/ (\ol{\cI}_N \otimes \cV^L_0) . 
$$

It follows that $\Delta_{\cV_0^L;N}$ induces a linear map from 
$\cW_0$ to $(\cW_0 \otimes \cW_0)/(\ol{\cJ}_N \otimes \cW_0)$, 
where $\ol{\cJ}_N$ is the image of $\ol{\cI}_N$ by the projection map 
$\cV_0 \to \cW_0$. 

It also induces an algebra morphism from $\cW_0$ to $\limm_{\leftarrow N}
(\cW_0\otimes\cW_0) / (\ol{\cJ}_N \otimes \cW_0)$. 

Then this algebra morphism factors through the coproduct map $\Delta_{\cW_0}$
defined above. To check this, it is enough to check it on generators $x[k]$
of $\cW_0$.   
 
It follows that  

\begin{lemma} \label{kably}
1) $i_\hbar$ induces a map $i : UL\B_+ \to \cW_0$, which is a surjective Hopf algebra
morphism. 

2) Let $\A_L$ be the Lie algebra of primitive elements of  $\cW_0$. 
The restriction $\iota_{|L\B_+}$ of $\iota$ to $L\B_+$ to $\A_L$ induces a surjective
Lie algebra morphism. 
\end{lemma}

{\em Proof.} $\Delta_{\cW_0}\circ i$ and  $(i \otimes i)\circ
\Delta_{U\B_+}$ are both algebra morphisms from $U\B_+$ to  $\cW_0
\otimes\cW_0$. Their values on the $x\otimes t^k$ coincide, therefore
they are equal. This shows 1). 

2) follows directly from 1) and from Prop.\ \ref{ratiu}. \hfill
\qed\medskip

\subsubsection{Construction of $\delta_{\cW_0}$}

Define $\hat\cV^L \hat\otimes \hat\cV^L$ as the tensor product 
$$ 
\CC[X_s[0]^{\cV^L(1)} ,X_s[0]^{\cV^L(2)}][[\hbar]]  \otimes
\CC[X_s[k]^{\cV^L(1)} ,X_s[k]^{\cV^L(2)},k>0] \otimes   \langle LV
\rangle \otimes_{\CC[[\hbar]]} \langle LV \rangle;
$$  
endow $\hat\cV^L\hat\otimes \hat\cV^L$  with the unique $\hbar$-adically
continuous algebra structure  such that $\cV^L \otimes_{\CC[[\hbar]]}
\cV^L \to  \cV^L \hat\otimes \cV^L$, $X_s[k]^{\cV^L(1)}\mapsto 
X_s[k]^{\cV^L} \otimes 1$, $X_s[k]^{\cV^L(2)}\mapsto  1\otimes
X_s[k]^{\cV^L}$, $1\otimes x\otimes 1\mapsto x\otimes 1$,   $1\otimes
1\otimes x\mapsto 1\otimes x$ (where $x$ is in $\langle LV \rangle$) is
an algebra morphism. 

Define in the same way $\hat\cV^L \hat\otimes \hat\cI_N^{(\infty)}$ as
the tensor product $ \CC[X_s[k]^{\cV^L(1)},X_s[k]^{\cV^L(2)}][[\hbar]]
\otimes \langle LV \rangle \otimes I_N^{(\infty)}$, where the tensor products are
over $\CC[[\hbar]]$.  Each $\hat\cV^L \hat\otimes \hat\cI_N^{(\infty)}$ is then a
left ideal of  $\hat\cV^L \hat\otimes\cV^L$. 

Clearly, we have $(\hat\cV^L \hat\otimes\hat\cV^L) / 
\hbar(\hat\cV^L \hat\otimes\hat\cV^L) = \cV_0^L\otimes\cV_0^L$. 
Moreover, 
$$
[(\hat\cV^L \hat\otimes\hat\cV^L) / (\hat\cV^L \hat\otimes \hat\cI_N^{(\infty)})]
/  \hbar [(\hat\cV^L \hat\otimes\hat\cV^L) / (\hat\cV^L \hat\otimes \hat\cI_N^{(\infty)})]
= (\cV_0^L\otimes\cV_0^L)  / (\cV_0^L\otimes\ol{\cI}_N) .  
$$

Then $\Delta_{\cV^L}$ is an algebra morphism from  $\cV^L$ to
$\limm_{\leftarrow N} (\hat \cV^L\hat\otimes \hat\cV^L)  / (\hat
\cI_N^{(\infty)} \hat\otimes \hat\cV^L)$. We again denote by $\Delta_{\cV^L}$ the 
composition of this map with the projection on  $\lim_{\leftarrow N} (
\hat\cV^L \hat\otimes \hat\cV^L) / (\hat\cI_N^{(\infty)} \hat\otimes\hat\cV^L  +
\hat\cV^L \hat\otimes \hat\cI_N^{(\infty)})$. Define $\Delta'_{\cV^L}$ as   
$\Delta_{\cV}^L$ composed with the exchange of factors. 

We have then $(\Delta_{\cV^L} - \Delta'_{\cV^L})(\cV^L)
\subset  \lim_{\leftarrow N} ( \hbar \hat\cV^L \hat\otimes
 \hat\cV^L) / [(\hat\cI_N^{(\infty)} \hat\otimes\hat\cV^L 
+ \hat\cV^L \hat\otimes \hat\cI_N^{(\infty)})\cap \hbar 
\hat\cV^L \hat\otimes \hat\cV^L]$.  

Since $\cI_N^{(\infty)}$ is divisible in $\cV^L$, we have 
$$
(\hat\cI_N^{(\infty)} \hat\otimes \hat\cV^L  + \hat\cV^L \hat\otimes 
\hat\cI_N^{(\infty)})
\cap \hbar (\hat\cV^L \hat\otimes \hat\cV^L) = 
\hbar (\hat\cI_N^{(\infty)} \hat\otimes \hat\cV^L 
+ \hat\cV^L \hat\otimes \hat\cI_N^{(\infty)} ) ,  
$$  
so 
${{\Delta_{\cV^L} - \Delta'_{\cV^L}
}\over{\hbar}}$ is a linear map from $\cV^L$ to 
$$
\lim_{\leftarrow N} ( \hat\cV^L \hat\otimes
 \hat\cV^L) / (\hat\cI_N^{(\infty)} \hat \otimes\hat\cV^L 
+ \hat\cV^L \hat\otimes \hat\cI_N^{(\infty)}) .  
$$

Define $\cI_N^{(1)}$ as the image of $\hat \cI_N^{(1)}$ in $\cV_0^L$ by the 
projection  $\cV^L\to\cV_0^L$. 

Let us set $\delta_{\cV_0^L} = {{\Delta_{\cV^L} - \Delta'_{\cV^L}
}\over{\hbar}}$ mod $\hbar$. Then $\delta_{\cV_0^L}$ is a linear map from 
$\cV^L_0$ to 
$$
\lim_{\leftarrow N} (\cV_0^L \otimes
 \cV_0^L) / (\ol{\cI}_N \otimes \cV_0^L 
+ \cV_0^L \otimes \ol{\cI}_N) .   
$$

Moreover, ${{\Delta_{\cV^L} - \Delta'_{\cV^L} }\over{\hbar}}$ maps
$\cI_N^{(\infty)}$ to the inverse limit  
$$
\lim_{\leftarrow M}
(\hat\cI_N^{(\infty)} \hat\otimes \hat\cV^L + \hat\cV^L\hat\otimes
\hat\cI_N^{(\infty)} ) / (\hat\cI_M^{(\infty)} \hat\otimes \hat\cV^L +
\hat\cV^L\hat\otimes \hat\cI_M^{(\infty)} ).
$$ 
Therefore, $\delta_{\cV_0^L}$
maps $\ol{\cI}_N$ to  $\lim_{\leftarrow M}
(\ol{\cI}_N \otimes \cV^L_0 + \cV^L_0\otimes \ol{\cI}_N ) 
/ (\ol{\cI}_M \otimes \cV^L_0 + \cV^L_0\otimes \ol{\cI}_M )$. 
Therefore, $\delta_{\cV_0^L} (\cap_N \ol{\cI}_N)$ is zero. It follows that 
$\delta_{\cV_0^L}$ induces a map $\delta_{\cW_0}$ from $\cW_0$ to 
$\lim_{\leftarrow N} (\cW_0 \otimes\cW_0) / (\ol{\cJ}_N \otimes \cW_0 + 
\cW_0 \otimes \ol{\cJ}_N)$. 

\subsubsection{Identities satisfied  by  $\delta_{\cW_0}$}

\begin{lemma} \label{peisakh}
$\delta_{\cW_0}$ satisfies 
\begin{equation} \label{co-leib:top}
(\Delta_{\cW_0}\otimes id)\circ \delta_{\cW_0} = (\delta_{\cW_0}^{2\to 23} 
+ \delta_{\cW_0}^{2\to 13}) \circ \Delta_{\cW_0},  
\end{equation}
\begin{equation} \label{co-jacobi:QC}
\Alt ( \delta_{\cW_0}\otimes id)\circ\delta_{\cW_0} ) = 0, 
\end{equation}
\begin{equation} \label{comp:Hopf:QC}
\delta_{\cW_0}(xy)  = \delta_{\cW_0}(x) \Delta_{\cW_0}(y)  
+ \Delta_{\cW_0}(x) \delta_{\cW_0}(y) \on{\ for\ } 
x,y \on{\ in\ } \cW_0, 
\end{equation}
where we use the notation of sect.\ \ref{adama}. 
The two first equalities are identities of maps from $\cW_0$ to 
$\limm_{\leftarrow N} \cW_0^{\otimes 3} / (\ol{\cJ}_N \otimes \cW_0^{\otimes 2}
+ \cW_0 \otimes \ol{\cJ}_N \otimes \cW_0 + \cW_0^{\otimes 2} \otimes 
\ol{\cJ}_N)$. 
\end{lemma}

{\em Proof.}
$\Delta_{\cV^L}$ maps $\cI_N^{(\infty)}$ to 
$$
\lim_{\leftarrow M}
(\hat\cI_N^{(\infty)} \hat\otimes \hat\cV^L + \hat\cV^L\hat\otimes
\hat\cI_N^{(\infty)} ) / (\hat\cI_M^{(\infty)} \hat\otimes \hat\cV^L +
\hat\cV^L\hat\otimes \hat\cI_M^{(\infty)} ).
$$ 
Therefore, $(\Delta_{\cV^L}\otimes id ) \otimes \Delta_{\cV^L}$
and $(id \otimes \Delta_{\cV^L}) \otimes \Delta_{\cV^L}$ both 
define algebra morphisms from $\cV^L$ to 
$\lim_{\leftarrow N} (\hat\cV^L)^{\hat\otimes 3}
 / [\hat \cI_N^{(\infty)} \hat\otimes (\hat\cV^L)^{\hat\otimes 2}  
+ \hat\cV^L \hat\otimes \hat \cI_N^{(\infty)} \hat\otimes \hat\cV^L
+ (\hat\cV^L)^{\hat\otimes 2} \hat\otimes 
\hat \cI_N^{(\infty)}]$. These morphisms are the restrictions
to $\cV^L$ of $(\Delta_{\cS^L}\otimes id ) \otimes \Delta_{\cS^L}$
and $(id \otimes \Delta_{\cS^L}) \otimes \Delta_{\cS^L}$, which 
coincide, therefore they coincide.  

The intersection $\cap_{N>0}\cI_N^{(\infty)}$ is a two-sided
ideal of $\cV^L$. Define $\cW$ as the quotient $\cV^L / 
\cap_{N>0}\cI_N^{(\infty)}$. Let $\cJ_N$ be the image of 
$\cI_N^{(\infty)}$ by the projection of $\cV^L$ on $\cW$. Define
$\hat\cW$ and $\hat \cJ_N$ in the same way, replacing $\cV^L$
and $\cI_N^{(\infty)}$ by $\hat\cV^L$
and $\hat\cI_N^{(\infty)}$. Then $\Delta_{\cV^L}$ induces an 
algebra morphism $\Delta_{\cW}$ from $\cW$ to $\limm_{\leftarrow N}
(\hat\cW\hat\otimes\hat\cW) / (\hat\cJ_N\hat\otimes\hat\cW + 
\hat\cW \hat\otimes\cJ_N)$. Moreover,   
$(\Delta_{\cW}\otimes id ) \circ \Delta_{\cW}$
and $(id \otimes \Delta_{\cW}) \circ \Delta_{\cW}$
define coinciding algebra morphisms from $\cW$ to 
$\limm_{\leftarrow N}
(\hat\cW^{\hat\otimes 3}
) / [\hat\cJ_N \hat\otimes \hat\cW^{\hat\otimes 2}
+ \hat\cW \hat\otimes \cJ_N \hat\otimes \hat\cW 
+ \hat\cW^{\hat\otimes 2} \hat\otimes \hat\cJ_N]$.  

Moreover, $\cW$ is a free $\CC[[\hbar]]$-module, and we have a 
topological Hopf algebra isomorphism of $\cW / \hbar\cW$ with $\cW_0$.
The usual  manipulations then imply the statements of the Lemma.
\hfill \qed\medskip

The identities of Lemma \ref{peisakh} are the topological versions of 
the co-Leibnitz, co-Jacobi and Hopf compatibility rules.

\subsubsection{Topological Lie bialgebra structure on $\A_L$}

Define $\A_L^{(N)}$ as the intersection $\A_L \cap \ol{\cJ}_N$. 

\begin{lemma} \label{coen}
$\ol{\cJ}_N$ is the left ideal $(U\A_L) \A_L^{(N)}$ of 
$\cW_0 = U\A_L$. Moreover, $\A_L^{(N)}$ is a Lie subalgebra of 
$\A_L$.   
\end{lemma}

{\em Proof.} $\cJ_N$ is a left ideal of $\cW$, therefore $\ol{\cJ}_N$
is a left ideal of $\cW_0$.  Moreover, $\Delta_{\cW}(\cJ_N)$
is contained in the inverse limit 
$\limm_{\leftarrow M} (\hat\cJ_N \hat\otimes \hat\cW + \hat\cW 
\hat\otimes \hat\cJ_N) / (\hat\cJ_M \hat\otimes \hat\cW + \hat\cW 
\hat\otimes \hat\cJ_M)$. It follows that $\Delta_{\cW_0}(\ol{\cJ}_N)$
is contained in $\ol{\cJ}_N \otimes \cW_0 + \cW_0 \otimes \ol{\cJ}_N$. 

The first statement of Lemma \ref{coen} now follows from Lemma \ref{lausanne}. 
For $x,y$ in $\A_L^{(N)}$, $[x,y] = xy - yx$ belongs to 
$\A_L$ and also to $(U\A_L)\A_L^{(N)}$, so it belongs to 
$\A_L^{(N)}$.  Therefore $\A_L^{(N)}$ is a Lie subalgebra of $\A_L$. 
\hfill \qed\medskip 

\begin{lemma} \label{basic:delta}
1) The restriction of $\delta_{\cW_0}$ to  $\A_L$ defines a map
$\delta_{\A_L} : \A_L \to \limm_{\leftarrow N} (\A_L \otimes\A_L) / 
(\A_L^{(N)} \otimes\A_L + \A_L \otimes \A_L^{(N)})$. 

2) For any element $x$ of $\A_L$, $\ad(x)(\A_L^{(N)})$ in contained in
$\A_L^{(N - k(x))}$. The tensor square of the  adjoint action therefore induces
a $\A_L$-module structure on  $\limm_{\leftarrow N} (\A_L \otimes\A_L) /
 (\A_L^{(N)} \otimes\A_L + \A_L \otimes \A_L^{(N)})$. $\delta_{\A_L}$ is
a $1$-cocycle of $\A_L$ with values in this module. 

3) We have $\delta_{\A_L}(\A_L^{(N)}) \subset \limm_{\leftarrow M}
(\A_L^{(N)} \otimes\A_L + \A_L \otimes\A_L^{(N)})   
/ (\A_L^{(M)} \otimes\A_L + \A_L \otimes\A_L^{(M)})$. 
$(\delta_{\A_L}\circ id) \circ \delta_{\A_L}$ therefore defines a map from 
$\A_L$ to $\limm_{\leftarrow N} \A_L^{\otimes 3} / (\A_L^{(N)} \otimes 
\A_L^{\otimes 2} + \A_L\otimes\A_L^{(N)} \otimes\A_L + \A_L^{\otimes 2}
\otimes\A_L^{(N)})$. It satisfies the rule 
\begin{equation} \label{coass}
\Alt(\delta_{\A_L}\otimes id)\circ\delta_{\A_L} = 0.
\end{equation} 
\end{lemma}

{\em Proof.} Let us show 1). $\delta_{\cW_0}$
induces an map $\delta_{\cW_0;N}$ from $\cW_0$ to $\cW_0^{\otimes 2} / 
(\ol{\cJ}_N \otimes \cW_0 + \cW_0 \otimes\ol{\cJ}_N)$. Let $a$ belong to $\A_L$. 
Let us write 
$\delta_{\cW_0;N}(a) = \sum_i a_i\otimes b_i$ mod 
$\ol{\cJ}_N \otimes \cW_0 + \cW_0 \otimes\ol{\cJ}_N$, with 
$(a_i)_i$ and $(b_i)_i$ finite families of $\cW_0$ such that 
$(b_i \on{\ mod\ }\ol{\cJ}_N)_i$
is a free family of $\cW_0 / \ol{\cJ}_N$. It follows from 
(\ref{co-leib:top}) that 
$$
\sum_i (\Delta_{\cW_0}(a_i) - a_i\otimes 1 - 1 \otimes a_i) \otimes b_i 
$$
belongs to $\ol{\cJ}_N \otimes\cW_0^{\otimes 2} + \cW_0\otimes
\ol{\cJ}_N \otimes \cW_0 + \cW_0^{\otimes 2} \otimes \ol{\cJ}_N$.  
Its image by the projection $\cW_0^{\otimes 3} \to [\cW_0^{\otimes 2} / 
(\ol{\cJ}_N  \otimes \cW_0 + \cW_0 \otimes \ol{\cJ}_N ) ]\otimes 
[\cW_0 / \ol{\cJ}_N]$ its therefore zero. 
It follows that each $a_i$ is such that $\Delta_{\cW_0}(a_i) 
- a_i\otimes 1 - 1 \otimes a_i$ belongs to $\ol{\cJ}_N \otimes\cW_0 
+ \cW_0\otimes \ol{\cJ}_N
$. Reasoning by induction on the degree of $a_i$
(for the enveloping algebra filtration of $\cW_0$), we find that $a_i$ 
belongs to $\A_L + \ol{\cJ}_N$. Therefore, $\delta_{\cW_0;n}(a)$ belongs to 
the image of $\A_L \otimes\cW_0$ in $\cW_0^{\otimes 2} / 
(\ol{\cJ}_N  \otimes \cW_0 + \cW_0 \otimes \ol{\cJ}_N )$. 
Since $\delta_{\cW_0;n}(a)$ is also antisymmetric, it belongs to the 
image of $\A_L \otimes\A_L$ in this space. This shows 1). 

Let us show 2). For $x$ en element of $\A_L$ and $y$ an element of 
$\A_L^{(N)}$, $[x,y] = xy - yx$ belongs to $\ol{\cJ}_N 
+ \ol{\cJ}_{N - k(x)} = \ol{\cJ}_{N - k(x)}$; since it also belongs to 
$\A_L$, $[x,y]$ belongs to $\A_L^{(N)}$. That $\delta_{\A_L}$ is a $1$-cocycle
then follows from (\ref{comp:Hopf:QC}). 

Let us show 3). $\Delta_{\cW;M}(\cJ_N)$ is contained in  $(\cJ_N \otimes
\cW + \cW \otimes \cJ_N)  / (\cJ_M \otimes \cW +  \cW \otimes \cJ_M)$. 
It follows that $\delta_{\cW_0;M}(\ol{\cJ}_N)$ is  contained in
$(\ol{\cJ}_N \otimes \cW_0 + \cW_0 \otimes \o{\cJ}_N)  /  (\ol{\cJ}_M
\otimes \cW_0 + \cW_0 \otimes \o{\cJ}_M)$. Therefore, 
$\delta_{\A_L}(\A_L^{(N)})$ is contained in $\limm_{\leftarrow M}
(\A_L^{(N)} \otimes\A_L + \A_L \otimes\A_L^{(N)}) / 
(\A_L^{(M)} \otimes\A_L + \A_L \otimes\A_L^{(M)})$. (\ref{coass})
in then a consequence of (\ref{co-jacobi:QC}). 
\hfill \qed\medskip 

Define the restricted dual $\A_L^\star$ of $\A_L$ as the subspace of 
$\A_L^*$ composed of the forms $\phi$ on $\A_L$, such for some $N$,
$\phi$  vanishes on $\A_L^{(N)}$. 

\begin{lemma} \label{szenes} 
The dual map to $\delta_{\A_L}$ defines a Lie algebra structure on $\A_L^\star$. 
\end{lemma}

{\em Proof.} Let $\phi,\psi$ belong to $\A_L^\star$. Let $N$ be an 
integer such that $\phi,\psi$ vanish on $\A_L^{(N)}$. For any integer $M$
let $\bar\delta_{\A_L;M}$ be a lift to $\A_L^{\otimes 2}$ of the 
map $\delta_{\A_L;M}$ from $\A_L$ to $\A_L^{\otimes 2} / (\A_L^{(M)})
\otimes\A_L + \A_L \otimes\A_L^{(M)}$ induced by $\delta_{\A_L}$. 

Let $x$ belong to $\A_L$. Then $M\geq N$, the number  $\langle
\phi\otimes\psi, \bar\delta_{\A_L;M}(x) \rangle$ is independent of 
the lift $\bar\delta_{\A_L;M}$ and of $M$; it defines
a linear form $[\phi,\psi]$ on $\A_L$. The first statement of Lemma
\ref{basic:delta},  3), implies that $[\phi,\psi]$ actually belongs to
$\A_L^\star$.  It is clear that $(\phi,\psi)\mapsto [\phi,\psi]$ is linear
and antisymmetric in $\phi$ and $\psi$. (\ref{coass}) implies that it
satisfies the Jacobi identity.  \hfill \qed \medskip 

\subsubsection{Topological Lie bialgebra structure on $L\B_+$}

Define for any integer $N$, $(L\B_+)^{(N)}$ as the Lie subalgebra of 
$L\B_+$ generated by the $\bar x_i^+\otimes t^k$, $k\geq N, i = 1,\ldots,n$.

For any $x$ in $L\B_+$, there exists an integer $l(x)$ such that 
$\ad(x)((L\B_+)^{(N)})$ is contained in $(L\B_+)^{(N - l(x))}$. It follows that 
$\limm_{\leftarrow N} (L\B_+)^{\otimes 2} / [(L\B_+)^{(N)} \otimes L\B_+ 
+ L\B_+ \otimes (L\B_+)^{(N)}]$
and $\limm_{\leftarrow N} (L\B_+)^{\otimes 3} / [(L\B_+)^{(N)} \otimes L\B_+^{\otimes 2} 
+ L\B_+ \otimes (L\B_+)^{(N)} \otimes L \B_+ + L\B_+^{\otimes 2} \otimes (L\B_+)^{(N)}]$
have $L\B_+$-module structures. 

\begin{lemma} There is a unique map $\delta_{L\B_+}$ from $L\B_+$ to 
$\limm_{\leftarrow N} (L\B_+)^{\otimes 2} / [(L\B_+)^{(N)} \otimes L\B_+ 
+ L\B_+ \otimes (L\B_+)^{(N)}]$, such that 
$\delta_{L\B_+} (\bar h_i\otimes t^k) = 0$ and 
$$
\delta_{L\B_+}(\bar x_i^+\otimes
t^k) = d_i \Alt [ {1\over 2 }(\bar h_i\otimes 1) \otimes (\bar x_i^+ \otimes 1)   
+ \sum_{l>0} (\bar h_i\otimes t^{-l}) \otimes (\bar x_i^+ \otimes t^{k+l}) ]
$$
and $\delta_{L\B_+}$ is a $1$-cocycle. Moreover, $\delta_{L\B_+}$ maps
$L\B_+^{(N)}$ to $\limm_{\leftarrow M} (L\B_+^{(N)} \otimes L\B_+ 
+ L\B_+ \otimes L\B_+^{(N)})  / (L\B_+^{(M)} \otimes L\B_+ 
+ L\B_+ \otimes L\B_+^{(M)})$, and it satisfies the co-Jacobi 
identity $\Alt (\delta_{L\B_+}\otimes id) \circ \delta_{L\B_+} = 0$. 
\end{lemma}

Define the restricted dual $(L\B_+)^\star$ to $L\B_+$ as the subspace of
 $(L\B_+)^*$ consisting of the forms on $L\B_+$, which vanish on some
$(L\B_+)^{(N)}$. The argument of Lemma \ref{szenes} implies that
$\delta_{L\B_+}$ induces a Lie algebra structure on $(L\B_+)^\star$. 

\begin{lemma}
Define on $\G\otimes\CC((t))$, the pairing $\langle , \rangle_{\G\otimes \CC((t))}$ 
as the tensor product of the invariant pairing on $\G$ and $\langle f,g \rangle 
= \res_0 (fg {{dt}\over t})$. $\langle , \rangle_{\G\otimes \CC((t))}$  
an isomorphism of $(L\B_+)^\star$ with the subalgebra $L\B_-$ of $L\G$
defined as $\HH \otimes \CC[[t]] \oplus \N_- \otimes \CC((t))$. This isomorphism 
is a Lie algebra antiisomorphism (that is, it is an isomorphism after we
change the bracket of $L\B_-$ into its opposite). 
\end{lemma}

The map $\iota_{|L\B_+}$ defined in Lemma \ref{kably}, 2), maps the
generators of  $L\B_+^{(N)} $ to $\A_L^{(N)}$; since  $\A_L^{(N)}$ is a
Lie subalgebra of $\A_L$ (Lemma \ref{coen}),  we have $\iota_{|
L\B_+}(L\B_+^{(N)}) \subset \A_L^{(N)}$. 

It follows that $\iota_{|L\B_+}$ induces a linear map $\iota^\star$ 
from $\A_L^\star$ to $(L\B_+)^\star = L\B_-$. 

Moreover, we have 
\begin{equation} \label{compat}
\delta_{\A_L} \circ \iota_{|L\B_+} = 
(\iota_{L\B_+}^{\otimes 2}) \circ \delta_{L\B_+}, 
\end{equation}
because both maps are $1$-cocycles of $L\B_+$ with values in 
$\limm_{\leftarrow N} \A_L^{\otimes 2} / (\A_L^{(N)} \otimes\A_L
+ \A_L \otimes \A_L^{(N)})$, and coincide on the generators of 
$L\B_+$. 

(\ref{compat}) then implies that $\iota^\star : \A_L^\star \to 
L\B_-$ is a Lie algebra morphism.

Let us set $(L\B_-)_{pol} = \HH\otimes \CC[t^{-1}] \oplus
\N_- \otimes\CC[t,t^{-1}]$; $(L\B_-)_{pol}$ is the polynomial 
part of $L\B_-$. 

\begin{lemma} \label{surj}
The image of $\iota^\star$ contains $(L\B_-)_{pol}$. 
\end{lemma}

{\em Proof.} 
As we have seen, $\cV^L$ is graded by $\NN^n$. Each ideal 
$\cI_N^{(\infty)}$is a graded ideal, so that $\cW = \cV / 
\cap_N \cI_N^{(\infty)}$ is also graded by $\NN^n$. Moreover, 
the degree $0$ and $\eps_i$ components of $\cI_N^{(\infty)}$
are respectively $0$ and $\oplus_{k\geq N}
\CC[[\hbar]][h_i[k]^{\cV^L},k\leq 0] t_i^k$. Therefore, 
 the components of $\cap_N \cI_N^{(\infty)}$ of degree $0$ and 
$\eps_i$ are zero. The components
of $\cW$ of degrees $0$ and $\eps_i$ are therefore respectively 
$\CC[[\hbar]] [h_i[k], k\leq 0]$ and $\oplus_{l\in\ZZ}
\CC[[\hbar]][h_i[k]^{\cV^L},k\leq 0] t_i^l$. 

$\cW_0$ is also graded by $\NN^n$, and its components of degrees 
$0$ and $\eps_i$ are $\CC[h_i[k], k\leq 0]$ and $\oplus_{l\in\ZZ}
\CC[h_i[k]^{\cV^L},k\leq 0] t_i^l$. 

The primitive part $\A_L$ of $\cW_0$ is therefore also graded by 
$\NN^n$, and the computation of $\Delta_{\cW_0}$ on $\cW_0[0]$  and
$\cW_0[\eps_i]$ shows that $\A_L[0] = \oplus_{1\leq i\leq n, k\leq 0}
\CC h_i[k]^{\cV^L}$ and $\A_L[\eps_i] = \oplus_{k\in\ZZ} \CC t_i^k$. 

Define linear forms $h_{i,k}^*$ and $e_{i,k}^*$ on $\A_L$ by the 
rules that $h_{i,k}^*$ vanishes on $\oplus_{\al\neq 0}\A_L[\al]$, 
and the restriction of $h_{i,k}^*$ to $\A_L[0]$ maps 
$h_j[l]^{\cV^L}$ to $\delta_{ij}\delta_{kl}$; and 
$e_{i,k}^*$ vanishes on $\oplus_{\al\neq \eps_i}\A_L[\al]$, 
and the restriction of $e_{i,k}^*$ to $\A_L[\eps_i]$ maps 
$t_i^l$ to $\delta_{kl}$.  

It follows from the computation of $\cI_N^{(\infty)}[0]$  and 
$\cI_N^{(\infty)}[\eps_i]$ that 
the $h_{i,k}^*$ and $e_{i,k}^*$ vanish on all the $\ol{\cJ}_N$, 
resp.\ on the $\ol{\cJ}_N, N\geq k$, adn therefore on all the 
$\A_L^{(N)}$, resp.\ on the $\A_L^{(N)}, N\geq k$. It follows that 
the $h_{i,k}^*$ and $e_{i,k}^*$ actully belong to $\A_L^\star$. 

Since the images of $\bar h_i\otimes t^k$ and 
$\bar x_i^+\otimes t^k$ by $\iota_{|L\B_+}$ are 
$h_i[k]^{\cV^L}$ and $t_i^k$, the images of 
$h_{i,k}^*$ and $e_{i,k}^*$ by $\iota^\star$ are the generators 
$\bar h_i\otimes t^k, 1\leq i\leq n, k\geq 0$ and 
$\bar x_i^-\otimes t^k, 1\leq i \leq n, k\in\ZZ$, of 
$(L\B_-)_{pol}$. 

The statement follows because $\iota^\star$ is a Lie algebra 
morphism. 
\hfill \qed\medskip 

Lemma \ref{surj} implies that the kernel of $\iota_{|L\B_+}$ is 
contained in contains the annihilator of $(L\B_-)_{pol}$ in $L\B_+$. 
Since this annihilator is zero, $\iota_{|L\B_+}$ is injective. 
It follows that $\iota_{|L\B_+}$ is an isomorphism. 

Therefore, $\iota: UL\B_+\to \cW_0$ is also an isomorphism. 
Recall that $\iota$ was obtained from the surjective 
$\CC[[\hbar]]$-modules morphism $\iota_\hbar = p\circ i_\hbar$, where 
$p$ is the projection of $\cV^L$ on $\cW$. 

We now use: 

\begin{lemma} \label{3.15}
Let $E$ and $F$ be $\CC[[\hbar]]$-modules, such that $F$ is 
torsion-free and $E$ is separated (i.e. $\cap_{N>0} \hbar^N E = 0$). 
Let $\pi : E\to F$ be a surjective morphism of $\CC[[\hbar]]$-modules, 
such that the induced morphism $\pi_0: E / \hbar E \to F / \hbar F$
is an isomorphism of vector spaces. Then $\pi$ is an isomorphism.  
\end{lemma}

{\em Proof.} Let $x$ belong to $\Ker p$. $\pi_0(x \on{\ mod\ }\hbar)$
is zero, therefore $x$ belongs to $\hbar E$. Set $x = \hbar x_1$. 
$\hbar \pi(x_1)$ is zero; since $F$ is torsion-free, $x_1$ belongs to 
$\Ker p$. Therefore, $\Ker p\subset \hbar \Ker p$. It follows that
$\Ker p\subset \cap_{N>0} \hbar^N E$, so that $\Ker p = 0$. It
follows that $\pi$ is an isomorphism.  
\hfill \qed\medskip

Recall that $U_\hbar L\N_+$ was defined as the quotient $\cA / 
(\cap_{N>0}\hbar^N \cA)$. It follows that $U_\hbar L\N_+$ is  separated.
The above Lemma therefore shows that $p\circ i_\hbar$ is an isomorphism.
Since $p$ and $i_\hbar$ are both surjective,  they are both
isomorphisms. Cor.\ \ref{cor:second} follows, together  with $\cap_{N>0}
\cI_N^{(\infty)} = 0$ (from where also follows that  $\cap_{N>0} \cI_N =
0$), and also,  by Lemma \ref{LV:free}, Thm.\ \ref{thm:third}, 1). 

It is then clear that the map $U_\hbar L\N_+  \to 
U_\hbar L\N_+^{top}$ is injective and that 
$U_\hbar L\N_+^{top}$ is the $\hbar$-adic completion of 
$U_\hbar L\N_+$. This proves Thm.\ \ref{thm:third}, 2). 

There is a unique algebra morphism $\varsigma$ from $\wt U_\hbar L\N_+$ 
to $U_\hbar L\N_+$, which sends each $e_i[k]^{\wt\cA}$ to $e_i[k]$. As 
we have seen in Prop.\ \ref{fargo}, $\varsigma$ induces an isomorphism 
between $\wt U_\hbar L\N_+ / \hbar\wt U_\hbar L\N_+$  and 
$U_\hbar L\N_+ / \hbar U_\hbar L\N_+$. Moreover,  $\wt U_\hbar L\N_+$  
is separated and by Thm.\ \ref{thm:third}, 1), $U_\hbar L\N_+$ is free. 
Lemma \ref{3.15} then implies that $\varsigma$ is an isomorphism. 
\hfill\qed\medskip 

\begin{remark} We have not been able to prove directly that  $\cA$ ot
$\wt \cA$ are themselves $\hbar$-adically separated. Since the
homogeneous components of $\cA$ are infinitely generated
$\CC[[\hbar]]$-modules,  it might happen that $\cA$ has components of
the type $\CC((\hbar))$ or $\CC((\hbar)) / \CC[[\hbar]]$. We cannot
exclude the existence of these components by the same argument as in the
proof of Thm.\  \ref{thm:first}, because their images by the map $\cA
\to \cA /  \hbar\cA$ are zero (whereas in the finitely generated case, 
the torsion submodules had a nonzero image by the same map). 
\end{remark}

\begin{remark} \label{open}
 Let $FO^{(0)}$ be the subspace of $FO$ formed of the 
functions satisfying $f(^{(i)}_\al) = 0$ when $t^{(i)}_1 
= q_{d_i}^2 t^{(i)}_2 = \cdots = q_{d_i}^{-2a_{ij}} t^{(i)}_{a_{ij}}
= q_{d_i}^{-a_{ij}}t^{(j)}_1$ for any $i,j$. We showed in \cite{Enr} that the image of
$U_\hbar L\N_+$ in $FO$ is contained in $FO^{(0)}[\hbar^{-1}]$. It is natural 
to expect that this image is actually the subspace of $FO^{(0)}$ consisting
in the functions such that $f(t_1,\cdots,t_N) = O(\hbar^k)$ whenever $k$
out of the $N$ variables $t_i$ coincide.
\end{remark}

\subsection{Nondegeneracy of the pairing $\langle , 
\rangle_{U_\hbar L\N_\pm}$ (proof of Thm.\ \ref{nondeg:L})}

Let us define $T(LV)$ as the tensor algebra $\oplus_{k\geq 0}
(LV)^{\otimes_{\CC[[\hbar]]}k}$, where $LV = \oplus_{i=1}^{n}
\CC[[\hbar]][t_i,t_i^{-1}]$.  Denote in this algebra, the element
$t_i^l$ of $LV$ as $f_i[l]^{(T)}$. 

Define a pairing 
$$
\langle , \rangle_{FO \times T(LV)} : FO \times T(LV) \to \CC((\hbar))
$$
as follows: if $P$ belongs to $FO_\kk$, 
\begin{align} \label{pairing:NR}
& \langle P, f_{i_1}[l_1]^{(T)} \cdots f_{i_{N'}}[l_{N}]^{(T)} 
\rangle_{FO \times T(LV)} 
\\ & \nonumber  
= \delta_{\kk,\sum_{j=1}^{N} \eps_{i_j}}
\res_{u_N = 0} \cdots \res_{u_1 = 0}
\left( P(t_1,\cdots,t_N) \prod_{l < l'} {{u_{l'} - u_l}\over
{q^{\langle \eps_{i_{l'}}, \eps_{i_l}\rangle} u_{l'} - u_l}} 
u_1^{l_1} \cdots u_N^{l_N} {{du_1}\over{u_1}}\cdots 
{{du_N}\over{u_N}} \right) , 
\end{align}
where we set as usual 
$(t_1, \ldots, t_{k_1}) = (t^{(1)}_1, \ldots, t^{(1)}_{k_1})$, etc., 
$(t_{k_1+ \cdots + k_{n-1} + 1}, \ldots, t_{k_1+ \cdots + k_n}) 
= (t^{(n)}_1, \ldots, t^{(n)}_{k_n})$, and 
$u_1 = t^{(i_1)}_1$, $u_2 = t^{(i_2)}_1$ if $i_2\neq i_1$
and in general 
$u_s = t^{(i_s)}_{\nu_s + 1}$, where $\nu_s$ is the number of indices $t$
such that $t<s$ and $i_t = i_s$. 

\begin{lemma} \label{volod}
The pairing $\langle , \rangle_{FO \times T(LV)}$ verifies $(T(LV))^{\perp} = 0$. 
\end{lemma}

{\em Proof.} Assume that the polynomial $P$ of $FO_\kk$ is such that (\ref{pairing:NR}) 
vanishes for any families of indices $(i_k)$ and $(l_k)$. Fix a family of indices 
$(i_k)$ such that $\kk = \sum_{j=1}^{N} \eps_{i_j}$. Since (\ref{pairing:NR}) 
vanishes for any family $(l_k)$, the rational function 
$P(t_1,\cdots,t_N) \prod_{l < l'} {{u_{l'} - u_l}\over
{q^{\langle \eps_{i_{l'}}, \eps_{i_l}\rangle} u_{l'} - u_l}}$ vanishes,
therefore $P$ is zero. \hfill \qed\medskip 

Let $\langle , \rangle_{\langle LV \rangle \times T(LV)}$ be the
restriction  of $\langle , \rangle_{FO \times T(LV)}$ to $\langle
LV\rangle \times T(LV)$.  Lemma \ref{volod} implies that $T(LV)^\perp =
0$ for this pairing. Using the isomorphism of Thm.\ \ref{thm:third}
between $\langle LV \rangle$ and $U_\hbar L\N_+$, we may view 
$\langle , \rangle_{\langle LV \rangle \times T(LV)}$ as a pairing 
$\langle , \rangle_{U_\hbar L\N_+ \times T(LV)}$ between 
$U_\hbar L\N_+$ and $T(LV)$. So again $T(LV)^\perp = 0$ for 
$\langle , \rangle_{U_\hbar L\N_+ \times T(LV)}$. 

Let $p$ be the quotient map from  $T(LV)$ to $U_\hbar L\N_+$. Composing
$\langle , \rangle_{U_\hbar L\N_+ \times T(LV)}$  with $p\otimes id$, we
get a pairing $\langle , \rangle_{T(LV) \times T(LV)}$  between $T(LV)$
and itself. It follows from (\ref{pairing:NR}) and (\ref{pdt:FO}) that 
$\langle , \rangle_{T(LV) \times T(LV)}$ is given by formula
(\ref{pairing:introd}).  Moreover, it follows  from \cite{Enr}, Prop.\
4.1 (relying on an identity of  \cite{Jing}) that $\langle ,
\rangle_{T(LV) \times T(LV)}$ induces a pairing $\langle ,
\rangle_{U_\hbar L\N_+ \times U_\hbar L\N_-}$ between $U_\hbar L\N_+$ 
and $U_\hbar L\N_-$. 

Since $\langle , \rangle_{U_\hbar L\N_+ \times U_\hbar L\N_-}$ is
induced by  the pairing $\langle , \rangle_{U_\hbar L\N_+ \times
T(LV)}$, and $T(LV)^\perp = 0$ for  this pairing, we get that $(U_\hbar
L\N_+)^\perp = 0$ for   $\langle , \rangle_{U_\hbar L\N_+ \times U_\hbar
L\N_-}$. Exchanging the roles of $U_\hbar L\N_+$ and $U_\hbar L\N_-$, 
we find that $(U_\hbar L\N_-)^\perp = 0$. 
Thm.\ \ref{nondeg:L}  follows. \hfill \qed \medskip 

\begin{remark} This argument is completely similar to the proof of 
Thm.\ \ref{thm:second}, the pairing between $U_\hbar L\N_+$ and $FO$ 
playing the role of the pairing beween $U_\hbar\N_+$ and $\Sh(V)$.
\end{remark}

\subsection{The form of the $R$-matrix (proof of Prop.\ 
\ref{R:mat:QC})}

Let us define $A_+^{a,b}$ as the subalgebra of $U_\hbar L\N_+$ generated by 
the $e_i[k]$, $i = 1, \ldots,n$, $a\leq k \leq b$.

\begin{lemma} \label{bar:mitswah} 
$A_+^{a,b}$ is a graded subalgebra of $U_\hbar L\N_+$. We have 
$A_+^{a,b} + I^+_{\leq a} + I^+_{\geq b}  = U_\hbar L\N_+$.  Moreover,
the graded components of $A_+^{a,b}$ are finite $\CC[[\hbar]]$-modules. 
\end{lemma}

{\em Proof.} Let us define $A_+^{\leq a}$ and $A_+^{\geq b}$ as the
subalgebras of $A_+$ generated by the $e_i[k]$, $k\leq a$ (resp.\ $k\geq
b$). It follows from  Thm.\ \ref{thm:third} that the product defines a
surjective  morphism from $A_+^{\leq a}\otimes A_+^{a,b} \otimes
A_+^{\geq b}$ to $A_+$.  The Lemma follows.  \hfill \qed\medskip

Since $I^+_{\geq a}[\al]  + I^+_{\leq b}[\al] \subset  (I^+_{\geq a}[\al] +
I^+_{\leq b}[\al])^{\perp\perp}$,  it follows from Lemma
\ref{bar:mitswah} that  $(I^+_{\geq a}[\al] + I^+_{\leq
b}[\al])^{\perp\perp}$ is a  submodule of $A_+$ with a complement of
finite type. Moreover, this module is also divisible, so that  $A_+
[\al]/ (I^+_{\geq a}[\al] + I^+_{\leq b}[\al])^{\perp\perp}$ is
torsion-free. Since it is finitely generated, it follows that  
$A_+[\al] / (I^+_{\geq a} [\al]+ I^+_{\leq b}[\al])^{\perp\perp}$ is
a free, finite-dimensional $\CC[[\hbar]]$-module. 

On the other hand, $(I^+_{\geq a}[\al] + I^+_{\leq b}[\al] )^{\perp}$ is a
submodule of $\Hom_{\CC[[\hbar]]} (A_+[\al] / (I^+_{\geq a}[\al] + I^+_{\leq b}[\al]) ,
\CC[[\hbar]])$, and is therefore a $\CC[[\hbar]]$-module of finite type.
It is a submodule of $A_-[-\al]$, so it is torsion-free. It follows that  
$(I^+_{\geq a} [\al]+ I^+_{\leq b}[\al] )^{\perp}$ is also a free,
finite-dimensional $\CC[[\hbar]]$-module.  

By construction, the pairing induced by $\langle , \rangle_{U_\hbar
L\N_\pm }$ between $(I^+_{\geq a}[\al] + I^+_{\leq b}[\al] )^{\perp}$
and  $A_+ [\al]/ (I^+_{\geq a}[\al] + I^+_{\leq b}[\al])^{\perp\perp}$
is nondegenerate. 

The fact that $P_{a,b}[\al]$ defines an element of 
$\limm_{\leftarrow a,b} A_+ / (I^+_{\leq a} +  I^+_{\geq
b})^{\perp\perp} \otimes_{\CC[[\hbar]]}  A_-[\hbar^{-1}]$ follows
from the following fact: if $F\subset G$ is an inclusion of finite
dimensional vector spaces, and $id_F$ and $id_G$ are the identity 
elements of $F\otimes F^*$ and $G\otimes G^*$, then their images
in $G\otimes F^*$ by the natural maps coincide.

Let $x$ be any product of the $f_i[k]$, with $c\leq k\leq d$. Then $x$
is orthogonal to $I^+_{\leq -c} + I^+_{\geq -d}$.  It follows that
$\cup_{a,b} (I^+_{\leq a} + I^+_{\geq b})^{\perp} = 0$,  therefore
$\cap_{a,b}(I^+_{\leq a} + I^+_{\geq b})^{\perp\perp} = 0$.

The proof of Prop.\ \ref{R:mat:QC} follows then the proof of Prop.\ \ref{R:mat}. 
\hfill \qed \medskip

\begin{remark} We have in the $\SL_{2}$ case
$$
P = \sum_r \sum_{i_1<\cdots<i_r, n_r\geq 0} {{\hbar^{n_1 + \cdots + n_r}}
\over{[n_1]^!_q \cdots [n_r]^!_q }
}e_{i_1}^{n_1} \cdots e_{i_r}^{n_r} \otimes  
f_{-i_1}^{n_1} \cdots f_{-i_r}^{n_r} ;  
$$
for $r\leq 2$, this formula is shown in \cite{Kh:D}, App.\ B. 
It would be interesting to obtain analogous explicit formulas in Yangian or
elliptic cases.   
\end{remark}

\section{Toroidal algebras (proofs of Props.\ \ref{oz}, \ref{nagila})} 
\label{toroidal}

\subsection{Proof Prop.\ \ref{oz}}

1) follows from the argument of the beginning of the proof 
of Prop.\ \ref{fargo}. The first statements of 2) are obvious. 

The proof of Prop.\ \ref{fargo} then implies that $j_+$ induces a surjective
Lie algebra morphism from  $\wt F$ to $\G\otimes\CC[t,t^{-1}]$, which
restricts to an isomorphim  between $\oplus_{\al\in\pm\Delta_+, \al\
\on{real}; k\in\ZZ} \wt F[(\al,k)]$  and  $(\oplus_{\al\in\pm\Delta_+,
\al\ \on{real}} \G[\al]) \otimes\CC[t,t^{-1}]$, which are the real roots
part of both Lie algebras, and that $\wt F_+[(\al,k)] = 0$ if $\al$ does not 
belong to $\Delta_+$. 

It follows that  $\Ker j_+$ is a graded subalgebra of $\wt F_+$, contained in 
$$
\oplus_{\al\in\Delta_+, \al\on{\ imaginary}, k\in\ZZ} \wt F_+[(\al,k)].
$$ 
\hfill \qed\medskip 

\subsection{Proof of Prop.\ \ref{nagila}}

Let us prove 1). Let us denote by $Z(\wt F_+)$ the center of $\wt F_+$. 
Let us first prove that $\Ker j_+\subset Z(\wt F_+)$.  Let $x$ belong to
$\Ker j_+$. We may assume that $x$ is  homogeneous of degree $n\delta$.
Then any  $n\delta + \al_i$, which is a real root. $[\wt e_i[k], x]$ is
homogeneous of degree $n\delta + \al_i$, which is a real root.  Since
the restriction of $j_+$ on the subspace of $\wt F$ of degree  $n\delta
+ \al_i$ is injective, $j_+([\wt e_i[k], x])$ is nonzero  unless $[\wt
e_i[k], x]$ is itself zero. But  $j_+([\wt e_i[k], x])$ is equal to
$[j_+(\wt e_i[k]), j_+(x)]$, which  is zero because $j_+(x) = 0$.
Therefore,  $[\wt e_i[k], x]$ is equal  to zero. 

On the other hand, since $j_+$ is surjective and the center of $L\N_+$ 
is zero, $\Ker j_+ = Z(\wt F_+)$. This proves 1).

Let us prove 2). The argument used in the proof of 1) implies that $\Ker
j_+$ is contained in the center $Z(\wt F)$ of $\wt F$. In the same way,
one proves that  $\Ker j_-$ is contained in $Z(\wt F)$, therefore $\Ker
j\subset Z(\wt F)$.  On the other hand, $\wt F$ is perfect. 

It follows that we have a surjective Lie algebra morphism $j' : 
\T \to \wt F$, such that the composition $\T\to \wt F \to L\G$
is the natural projection of $\T$ on $L\G$. Let $j'_+$ be the restriction of 
$j'$ to $\T_+$. For any $i,k$, we have 
$j'_+(e_i[k]^{\T}) = \wt e_i[k] + k_{i,k}$, with $k_{i,k}$ in $\Ker(j)$. 
Let $\la$ be any linear map from $\T_+$ to $\Ker(j)$, such that 
$\la(e_i[k]^{\T}) = k_{i,k}$. Set  $\wt j'_+ = j'_+ - \la$. Then 
$\wt j'_+$ is a Lie algebra map from $\T_+$ to $\wt F$. Since 
$\wt j'_+(e_i[k]^{\T}) = \wt e_i[k]$ and the $e_i[k]^{\T}$ generate
$\T_+$, the image of $\wt j'_+$ is contained in $\wt F_+$. 

Moreover, $\wt j'_+$ is graded, and it coincides with $j'_+$ on the
nonsimple roots subspace $[\T_+,\T_+] = \oplus_{\al\in\Delta_+ \setminus
\{\eps_i\}} \T_+[\al]$ of $\T_+$.  It follows that the restrictions of 
$j'$ on $[\T_+,\T_+]$ and $[\T_-,\T_-]$ are graded. 

Let us show that $\wt j'_+$ is surjective. 
Since the 
composition of $\wt j'_+$ with the projection $j:\wt F_+ \to L\N_+$ is 
the natural projection, it suffices to show that any element $x$ of 
$\Ker j_+$ is contained
 in $\wt j'_+(\T_+)$.  $x$ belongs to the image of $j'$, so let us set
$x = j'(y)$, with  $y = y_+ + y_- + y_0$, $y_\pm$ in $[\T_\pm,\T_\pm]$
and $y_0$ in $\HH_\T \oplus \oplus_{i = 1}^n\T_+[\eps_i]  \oplus
\oplus_{i = 1}^n\T_-[-\eps_i]$, where $\HH_\T$ is the Cartan subalgebra 
of $\T$ (defined as $\HH[\la^{\pm 1}] \oplus Z_0$, see 
Rem.\ \ref{rem:generalizations}). 
Then  $j'(y_\pm)$ belong to $[\wt F_\pm,\wt F_\pm]$ and
$j'(y_0)$ belongs to  $\wt H \oplus \oplus_{i = 1}^n \wt F_+[\eps_i] 
\oplus \oplus_{i = 1}^n \wt F_-[-\eps_i] \oplus \Ker(j)$. Moreover,  the
map from $\HH_\T \oplus \oplus_{i = 1}^n\T_+[\eps_i]  \oplus \oplus_{i =
1}^n\T_-[-\eps_i]$ to  $\wt H \oplus \oplus_{i = 1}^n \wt F_+[\eps_i] 
\oplus \oplus_{i = 1}^n \wt F_-[-\eps_i]$  induced by $j'$ is injective,
therefore $y_0 = 0$. It follows that $y_- = 0$ and $x = j'(y_+) = 
\wt j'_+(y_+)$, because $\wt j'_+$ coincides with $j'$ on $[\T_+,\T_+]$. 

\begin{lemma}
1) Assume that $A$ is not of type $A_1^{(1)}$. 
There is a unique Lie algebra map $j''$ from $\wt F_+$ to $\T_+$ such that 
$j''(\wt e_i[k]) = e_i[k]^\T$, for any $i = 0,\ldots, n$ and $k$ integer. 

2) Assume that $A$ is the Cartan matrix of type $A_1^{(1)}$. 
There is a unique Lie algebra map $j''$ from $\wt F_+$ to $\T_+ / 
\oplus_{l\in\ZZ} \CC K_\delta[l] $ such that 
$j''(\wt e_i[k]) = e_i[k]^\T$, for any $i = 0,1$ and $k$ integer. 
\end{lemma}

{\em Proof.} One should just check that the defining relations of $\wt
F_+$ are satisfied by the $e_i[k]^{\T}$  (in the $A_1^{(1)}$ case, by
the images of $e_i[k]^{\T}$  in  $\T_+ /  \oplus_{l\in\ZZ} \CC
K_\delta[l]$). This is the case when $A$ is not of type 
$A_1^{(1)}$,  because in that case we set $e_i = \bar x_i \otimes
\la^{\delta_{i0}}$ and we always have $\langle \bar x_i, \bar
x_j\rangle_{\bar\G} = 0$ for $i\neq j$. 

If $A$ is of type $A_1^{(1)}$, we have $x_0 = \bar f \otimes \la$, 
$x_1 = \bar e$, therefore 
$$
[x_0[l], x_1[m]] = ( - \bar h\la\otimes t^{l+m}, - m K_\delta[l+m]),  
$$
so that $[x_0[l+1], x_1[m]]  = [x_0[l], x_1[m+1]]$ holds in 
$\T_+ /  \oplus_{l\in\ZZ} \CC K_\delta[l]$. 
\hfill \qed\medskip 

Let us now prove Prop.\ \ref{nagila}, 2). The composition  $\wt j'_+
\circ j''$ are Lie algebra maps from  $\wt F_+ \to \T_+ \to \wt F_+$ 
($\wt F_+ \to \T_+ / \oplus_{k\in\ZZ} \CC K_\delta[k]  \to \wt F_+$  in
the $A_1^{(1)}$ case), which map the generators $\wt e_i[k]$ to 
themselves. Therefore, $\wt F_+$ can be viewed as a subalgebra of 
$\T_+$ (of  $\T_+ / \oplus_{k\in\ZZ} \CC K_\delta[k]$ in the 
$A_1^{(1)}$ case). This subalgebra contains the elements  $e_i[k]^\T$ of
$\T_+$ (resp.\ of $\T_+ / \oplus_{k\in\ZZ}  \CC K_\delta[k]$). Since the
Lie subalgebra of $\T_+$ generated by  the $e_i[k]^\T$ is $\T_+$ itself,
the image of $\wt F_+$ is equal to  $\T_+$ (resp.\ to $\T_+ /
\oplus_{k\in\ZZ}  \CC K_\delta[k]$). This proves Prop.\ \ref{nagila},
2).  \hfill \qed\medskip 

\begin{remark} \label{rem:india}  Prop.\ \ref{nagila}, 1), can also be
obtained using the presentation given in \cite{Moody} of $\T$. In this
paper, one shows that $\T$ is isomorphic to the algebra $\dot F$ with
generators $\dot e_i^\pm[k],\dot h_i[k]$ and  $\dot c$ and relations 
$\ad(\dot e_i^\pm[0])^{1 - a_{ij}}(\dot e_j^\pm[k]) = 0$, $[\dot
e^\pm_i[k] ,  \dot e^\pm_i[l]] = 0$, $[\dot h_i[k], \dot e^\pm_j[l]] =
\pm a_{ij} \dot e^\pm_j[k+l]$,  $[\dot e_i^+[k],  \dot e_j^-[l]] =
\delta_{ij} \dot h_i[k+l] + i \delta_{i+j,0} \langle  e_i, f_i
\rangle_{\bar\G} c$,  $[\dot h_i[k],\dot h_j[k]] = k \delta_{k+l,0} \langle h_i,
h_j\rangle_{\bar\G} c$,  $c$ central. It is then clear that there is a
Lie algebra map from $\T$ to $\wt F$. On the other hand, the system of
relations $[\dot e_i^+[k],  \dot e_i^+[l]] = \ad(\dot e_i^\pm[0])^{1 -
a_{ij}}(\dot e_j^\pm[k]) = 0$ is {\it not} a presentation of $\T_+$, 
because the ideal generated by these relations is not preserved by the
analogues of the $\Phi_{i,k}^\pm$ of the proof of Lemma \ref{india}. 
\hfill \qed\medskip 
\end{remark}

\begin{remark} \label{rem:generalizations} 
{\em Toroidal Manin triples.} 
It is easy to define an extension of the Lie algebra $\T$ with an 
invariant scalar product. Recall first (\cite{Moody,Kassel}) 
that if $\G$ is the central 
extension of the Lie algebra $\bar\G[\la,\la^{-1}]$, $\T$
is the universal central extension of 
$\bar\G[\la^{\pm 1},\mu^{\pm 1}]
$. We have therefore 
$$
\T = \bar\G[\la^{\pm 1},\mu^{\pm 1}] \oplus Z(\T). 
$$
$Z(\T)$ is isomorphic to $\Omega^1_\AAA / d\AAA$, where $\AAA = 
\CC[\la^{\pm 1},\mu^{\pm 1}]$. We have 
$$
Z(\T) = \oplus_{k,l\in\ZZ} K_{k\delta}[l] \oplus \CC c, 
$$
with 
$K_{k\delta}[l] = $ the class of ${1\over k} \la^k \mu^{l-1} d\mu$
if $k\neq 0$ , $K_0[l] = $ the class of $\mu^l {{d\la}\over{\la}}$, 
$c = $ the class of ${{d\mu}\over\mu}$.  

Define for $k,l$ in $\ZZ$, $\wt D_{k\delta}[l]$ as the
derivations of $\bar\G[\la^{\pm 1},\mu^{\pm 1}]$
equal to $\la^k \mu^l (l \la \pa_\la - k \mu \pa_\mu)$ if $k\neq 0$
and to   $\mu^l \la \pa_\la$ if $k = 0$, and $\wt d$ as the 
derivation $\mu \pa_\mu$. 

Endow $\CC^{\times 2}$ with the coordinates $(\la,\mu)$  and consider on
this space the Poisson structure defined  by $\{\la,\mu\} = \la\mu$. Let
$\Ham(\CC^{\times 2})$ be the  Lie algebra of Hamiltonian vector fields
on $\CC^{\times 2}$ generated by the functions $\la^k\mu^l$,
$k,l\in\ZZ^2$,   $\log\la$ and $\log\mu$. For any function $f$ on
$\CC^{\times 2}$, denote by $V_f$ the corresponding Hamiltonian vector
field.   Then $\Ham(\CC^{\times 2})$ is a Lie algebra, and the map 
$V_{\la^k\mu^l} \mapsto \wt D_{k\delta}[l]$, for $(k,l)\neq (0,0)$, 
$V_1 \mapsto 0$,  $V_{\log \la} \mapsto \wt D_{0}[0]$, $V_{\log \mu}
\mapsto  \wt d$, defines a Lie algebra map from $\Ham(\CC^{\times 2})$ 
to $\Der(\bar\G[\la^{\pm 1},\mu^{\pm 1}])$. 

The formula $V_f (\sum_i a_i db_i) = \sum_i \{f,a_i\} db_i + a_i
d\{f,b_i\}$  defines an action of $\Ham(\CC^{\times 2})$ on
$\Omega^1_{\AAA} / d\AAA$,  that is on $Z(\T)$. Define  $\bar
D_{k\delta}[l]$ and  $\bar d$ as the  following endomorphisms of $\T$: 
$\bar D_{k\delta}[l](x,0) = (\wt D_{k\delta}[l](x),0),   \bar d(x,0) =
(\wt d(x),0)$, and $\bar D_{k\delta}[l](0,\omega) =  (0,V_{\la^k
\mu^l}(\omega))$ for $(k,l)\neq (0,0)$,  $\bar D_{0}[0](0,\omega) = 
(0,V_{\log \la}(\omega))$,  $\bar d(0,\omega) =
(0,V_{\log\mu}(\omega))$.  These endomorphisms again define derivations
of $\T$, and we have now a Lie algebra map from $\Ham(\CC^{\times 2})$
to $\Der(\T)$. Let  $\wt\T$ be the corresponding crossed product Lie
algebra of $\T$  with $\Ham(\CC^{\times 2})$. We denote by
$D_{k\delta}[l]$ and $d$ the elements of $\wt t$ implementing the
extensions of the  derivations $\bar D_{k\delta}[l]$ and  $\bar d$ to
$\T$.

Define for $a,b$ integers, $x[a,b]$ as the element $(x\otimes \la^a \mu^b)$
of $\bar\G[\la^{\pm 1}, \mu^{\pm 1}]$. Define the bilinear form 
$\langle , \rangle_{\wt\T}$ by 
$$
\langle x[a,b], x'[a',b'] \rangle_{\wt\T} = \langle x,x'\rangle_{\bar\G}
\delta_{a+a',0} \delta_{b+b',0}, 
\langle D_{k\delta}[l] , K_{k'\delta}[l'] 
\rangle_{\wt\T} = \delta_{k+k',0} \delta_{l+l',0}, 
\langle d, c \rangle_{\wt\T}  = 1,  
$$
and all other pairings of elements $x[a,b],K_{k\delta}[l],D_{k'\delta}[l'],c$ 
and $d$ are zero. 

Then $\langle , \rangle_{\wt\T}$ is an invariant nondegenerate bilinear form 
on $\wt\T$.

Let us define $D$ as the image of $\Ham(\CC^{\times 2})$ in $\wt\T$. Let
us set  $D_>, D_<$ and $D_0$ as its subspaces  $\oplus_{k>0,l\in\ZZ} \CC
D_{k\delta}[l]$, $\oplus_{k<0,l\in\ZZ} \CC D_{k\delta}[l]$ and
$\oplus_{l\in\ZZ} \CC D_{0}[l] \oplus \CC d$. We have then $D = D_> \oplus 
D_0 \oplus D_<$. In the same way, define $Z_>, Z_<$ and $Z_0$ as the subspaces 
$\oplus_{k>0,l\in\ZZ} \CC
K_{k\delta}[l]$, $\oplus_{k<0,l\in\ZZ} \CC K_{k\delta}[l]$ and
$\oplus_{l\in\ZZ} \CC K_{0}[l] \oplus \CC c$ of $Z(\T)$. We have then 
$Z(\T) = Z_> \oplus Z_< \oplus Z_0$. 

Recall we defined $\HH_\T$ as the subalgebra  $\bar\HH[\la^{\pm 1}]
\oplus Z_0$ of $\T$. $\wt \HH_\T = \HH_\T\oplus D_0$ is then a Lie
subalgebra of $\wt\T$. In the spirit of the new realizations, we split
$\wt \HH_\T$ in two parts. 

Let us set $\HH_+ = \bar\HH[\la] \oplus Z_0$,  
$\HH_- = \bar\HH[\la^{-1}] \oplus D_0$; then $\HH_+ + \HH_- 
= \wt\HH_\T$, and $\HH_+ \cap \HH_-$ is $\bar\HH$. 

Define $L\N_+$  and $L\N_-$ as the linear spans of the  $x[a,b]$,
$a\in\ZZ, b>0$ ($b\geq 0$ if $x\in\bar\N_+$), resp.\ of the   $x[a,b]$,
$a\in\ZZ, b<0$ ($b\leq 0$ if $x\in\bar\N_-$). $L\N_\pm$ are Lie
subalgebras of $\bar\G[\la^{\pm 1},\mu^{\pm 1}]$.    $L\N_\pm \oplus
Z_\pm$ and $\wt{L\N}_\pm = L\N_\pm \oplus D_\pm \oplus Z_\pm$ are also Lie
subalgebras of $\wt\T$. Set $\wt\T_\pm = \wt{L\N}_\pm \oplus \HH_\pm$.  

Endow $\wt\T\times \bar\HH$ with the scalar product 
$\langle , \rangle_{\wt\T\times \bar\HH}$ defined by 
$$
\langle (x,h) , (x',h') \rangle_{\wt\T\times \bar\HH}  =
\langle x, x' \rangle_{\wt\T} -  \langle h, h' \rangle_{\bar\HH}.
$$  
Let $p_\pm$ be the natural projection of $\wt \T_\pm$ 
on $\bar \HH$. Identify $\wt\T_\pm$ as the Lie subalgebras
of $\{(x, \pm p_\pm(x)), x\in \wt\T_\pm\}$ of
$\wt\T\times \bar\HH[\la,\la^{-1}]$. 
$\wt\T_\pm$ are supplementary isotropic subspaces of 
$\wt\T\times \bar\HH[\la,\la^{-1}]$ and define therefore a 
Manin triple. This Manin triple is a central and cocentral 
extension (by $Z(\T)$ and $D$) of the Manin triple 
$$
(\bar\G[\la^{\pm1},\mu^{\pm1}] \times \bar\HH, 
L\N_+ \oplus \bar\HH[\la],  L\N_- \oplus \bar\HH[\la^{-1}])
$$ 
which is a part of the new realizations Manin triple
$(\G[\mu^{\pm1}]\times \HH, L\B_+,L\B_-)$.
  
One may also consider ``intermediate'' Manin triples, for example
$$
\left( \{\bar\G[\la^{\pm1},\mu^{\pm1}] \oplus Z_> \oplus D_< \} \times \bar\HH, 
L\N_+ \oplus \bar\HH[\la] \oplus Z_>,  L\N_- \oplus \bar\HH[\la^{-1}] \oplus D_<
\right).  
$$ 

It is a natural problem to quantize the corresponding Lie 
bialgebra structures on $\wt\T_\pm$. For this, one can think of the 
following program: 

1) to compute the centers of $U_\hbar L\N_+$ and (following
\cite{FO:new}) the center of $FO$. By duality, these central elements
should provide derivations of $U_\hbar L\N_-$ (and $FO$) of imaginary
degree. Compute these derivations and relations between them. One could
expect that the algebra generated by the derivations is some difference
analogue of the Lie algebra $\Ham(\CC^{\times 2})$.

2) it should then be easy, following Thm.\ \ref{thm:third},  to prove
that  the analogue of $i_\hbar$ is an isomorphism, and to derive from
there  the quantization of the Lie bialgebra $L\B_+$. 

We hope to return to these questions elsewhere. 
\end{remark}


\appendix

\section{Lemmas on $\CC[[\hbar]]$-modules}

\begin{lemma} \label{str:modules}
Let $E$ be a finitely generated $\CC[[\hbar]]$-module. Let $E_{tors} = 
\{x\in E | \hbar^k x = 0 $ for some $k >0\}$ be the torsion part of $E$. 
Then $E_{tors}$ is isomorphic to a direct sum $\oplus_{i=1}^p \CC[[\hbar]] / 
(\hbar^{n_i})$, where $n_i$ are positive integers, and $E$ is isomorphic to the
direct sum of $E_{tors}$ and a free module $\CC[[\hbar]]^{p'}$. 
\end{lemma}

{\em Proof.} As $E$ is finitely generated, we have a surjective
$\CC[[\hbar]]$-modules morphism $\CC[[\hbar]]^N \to E$. Let $K$ be the
kernel of this morphism. Then  $E$ is isomorphic to $\CC[[\hbar]]^N /
K$. 

Let us determine the form of $K$. Let us set $\bar K_i = K \cap \hbar^i
\CC[[\hbar]]^N$.  Then we have $\hbar\bar K_i\subset \bar K_{i+1}$. Let
us set $E_0 = \CC^N$, and $F_i =  \hbar^{-i}\bar K_i$ mod $\hbar$. Then
we have $F_0 \subset F_1 \subset\cdots \subset E_0$. Let $p$ the integer
such that  $F_k = F_p$ for $k\geq p$.  We can then find a basis 
$(v_i)_{1\leq i\leq N}$ of $E_0$ such that $(v_1,\ldots,v_{\dimm F_0})$
is a  basis of $F_0$, $(v_1,\ldots,v_{\dimm F_1})$ is a  basis of
$F_1$, etc., $(v_1,\ldots,v_{\dimm F_p})$ is a  basis of $F_p$. Then $K$
is the submodule $ \oplus_{i}  (\oplus_{k = \dimm F_{i-1} + 1}^{\dimm
F_i}\hbar^i \CC[[\hbar]] v_k)$ of $E_0[[\hbar]]$.

It follows that the quotient $E_0[[\hbar]] / K$ is isomorphic to a direct 
sum $\oplus_{i=1}^p \CC[[\hbar]] / (\hbar^{n_i}) \oplus \CC[[\hbar]]^{p'}$. 
The statement of the Lemma follows. 
\hfill \qed \medskip

\begin{cor} Any $\CC[[\hbar]]$-submodule of a finite-dimensional free 
$\CC[[\hbar]]$-module is free. 
\end{cor}

{\em Proof.} This follows from the fact that such a submodule has no torsion 
and from the above Lemma. \hfill \qed\medskip 

We have also 

\begin{lemma} \label{free:inf:dim}
Let $E$ be a free $\CC[[\hbar]]$-module with countable
basis $(v_i)_{i\geq 0}$. Any countably generated $\CC[[\hbar]]$-submodule of 
$E$ is free and has a countable basis. 
\end{lemma}

{\em Proof.} We repeat the reasoning of the proof of  Lemma
\ref{str:modules}. Let $(w_i)_{i\geq 0}$ be a countable family of $E$
and  let $F$ be the sub-$\CC[[\hbar]]$-module of $E$ generated by the 
$w_i$.
 
Set $\bar F_i = F\cap \hbar^i E$ and $F_i = \hbar^{-i}F$ mod $\hbar$. 
Generating families and bases for the $F_i$ can be constructed inductively 
as follows. 

A generating family for $F_0$ is $(w_i$ mod $\hbar)_{i\geq 0}$.  We can
then construct by induction a partition of $\NN$ in subsets $(i_k)$  and
$(j_k)$ such that $(w_{i_k}$ mod $\hbar)_{k\geq 0}$  is a basis of 
Span$(w_i$ mod $\hbar)_{i\geq 0}$. 

Let $\la_{kk'}$ be the scalars such that $w_{j_k} - \sum_{k'} \la_{kk'}
w_{i_{k'}}$ belongs to $\hbar E$.  Set $w_k^{(1)} = \hbar^{-1}
[w_{j_k} - \sum_{k'} \la_{kk'} w_{i_{k'}}]$. Then a generating family of 
$F_1$ is $(w_{i_k},w^{(1)}_k$ mod $\hbar)$. We then construct by 
induction a partition of $\NN$ in subsets $(i_k^{(1)})$ and $(j_k^{(1)})$
such that $(w_{i_k},w^{(1)}_{i_k^{(1)}}$ mod $\hbar)$ is a basis of $F_1$.

It is clear how to continue this procedure. Then 
$(w_{i_k},w^{(1)}_{i_k^{(1)}},w^{(2)}_{i_k^{(2)}} \ldots)$ forms a basis
of  $F$. \hfill \qed 
\end{document}